\documentclass[a4paper,twoside,12pt,english]{article}

\usepackage[applemac]{inputenc}
\usepackage[english]{babel}

\makeatletter
\renewcommand\section{\@startsection
 {section}{1}{0mm}%
 {-2\bigskipamount}%
 {\bigskipamount}%
 {\normalfont\normalsize\bfseries}%
}
\renewcommand\thesection{\arabic{section}}

\newcommand\dateymd{\number\year, \ifcase\month\or
 January\or February\or March\or April\or May\or June\or
 July\or August\or September\or October\or November\or
 December\fi, \number\day}
\newcommand\printtime{%
 \c@hours=\time \divide\c@hours by60
 \c@minutes=\c@hours \multiply\c@minutes by-60
 \advance \c@minutes by \time
 \ifnum\c@hours<10 0\fi\the\c@hours:%
 \ifnum\c@minutes<10 0\fi\the\c@minutes}
\makeatother

\usepackage[OT1]{fontenc}

\usepackage{fancyhdr}
\usepackage{relsize}
\usepackage{bm} 
\usepackage{amssymb}
\usepackage{amsmath}
\usepackage{amsthm}
\usepackage{enumerate}

\usepackage{verbatim} 

\usepackage[nodvipsnames]{color}
\newcommand\bfgr[1]{\textcolor[cmyk]{0,0,0,0.75}{\textbf{#1}}}

\usepackage[pdftex]{graphics} 

\newcommand\upla{\vskip0pt}

\newcommand\bnou{\textbf{\lower3.7pt\hbox{*}:~}}
\newcommand\enou{\unskip\textbf{~:\lower3.7pt\hbox{*}} }
\newcommand\bvell{\textbf{\lower3.7pt\hbox{*}:~}$\langle$}
\newcommand\evell{\unskip$\rangle$\textbf{~:\lower3.7pt\hbox{*}} }
\newcommand\bbnou{\textbf{\lower3.7pt\hbox{**}:~}}
\newcommand\eenou{\unskip\textbf{~:\lower3.7pt\hbox{**}} }

\newcommand\ie{i.\,e.~}
\newcommand\ifoi{\,\hbox{if\kern2.5pt and\kern2.5pt only\kern2.5pt if}\,{} }

\newcommand\df{\bfseries}
\newcommand\dfc[1]{\,{\df#1}\,}
\newcommand\dfd[1]{\,{\df#1}\hskip1pt}
\newcommand\secpar[1]{\S\,{#1}}

\newcommand\ensep{\unskip\hskip.65em\ignorespaces}

\newcommand\atilde{\lower3.5pt\hbox{\~{}}}
\newcommand\underl{\lower3.5pt\hbox{-}}

\newcommand\remarks{\bigskip\noindent\textit{Remarks}\par}
\newcommand\remark{\bigskip\noindent\textit{Remark}\par}

\newcommand\better{$\succ$}
\newcommand\tied{$\sim$}

\newcommand\pq[2]{\raise.25ex\hbox{\footnotesize${#1}\over{#2}$}%
\hskip-.35ex\null}
\newcommand\onehalf{\leavevmode\raise.5ex\hbox{\scriptsize$1$}\hskip-.25ex/\hskip-.25ex\lower.2ex\hbox{\scriptsize$2$}}

\newcommand\halfsmallskip{\vskip0.5\smallskipamount}

\newcommand\xxxx[1]{%
 \hangindent2.5\parindent
 \hangafter1
 \noindent\hskip.5\parindent
 \hbox to2\parindent{\hss#1\hss}}
\newcommand\condition[2]{\xxxx{#1}\textit{#2}.}
\newcommand\iim[1]{\xxxx{\textup{(#1)}}\ignorespaces}

\newcommand\brwrap[1]{[\textsl{#1}\kern1pt]}

\newcommand\refa[1]{\brwrap{#1}}
\newcommand*\drefa[2]{\brwrap{#1\kern.5pt:\,{\relscale{0.95}#2}}}

\makeatletter
\renewcommand\@biblabel[1]{\refa{#1}}
\renewcommand\@cite[2]{\hbox{\brwrap{#1\if@tempswa\/\upshape\,:\,{\relscale{0.95}#2}\fi}}} 
\makeatother

\newcommand\latop[2]{{#1\atop#2}}

\newcommand\sbset{\subset}
\newcommand\spset{\supset}
\newcommand\sbseteq{\subseteq}
\newcommand\spseteq{\supseteq}

\newcommand\cd[1]{\!#1\!}

\newcommand\stv{{}^\ast\kern-.25pt v}

\newcommand\tcl[1]{#1{}^\ast}
\newcommand\tclp[1]{\tcl{(#1)}}

\newcommand\isc{v^\ast}
\newcommand\mast{\nu} 
\newcommand\img{m^\mast}

\newcommand\rxi{\mathrel{\smash{\succ\kern-1.7ex\raise1.15ex\hbox{\mathsurround0pt$\scriptscriptstyle\xi$}\kern.4ex}}}
\newcommand\rxieq{\mathrel{\smash{%
 \vbox{\offinterlineskip\halign{\hfil##\hfil\cr
 \mathsurround0pt$\succ$\cr
 \noalign{\vskip-.5ex}%
 \mathsurround0pt$-$\cr
 \noalign{\vskip-1.15ex}%
 }\vss}\kern-1.05ex\raise1.15ex\hbox{\mathsurround0pt$\scriptscriptstyle\xi$}\kern.4ex}}}

\newcommand\ppmg{m^\sigma}
\newcommand\ppto{t^\sigma}

\newcommand\psc{v^\pi}
\newcommand\pmg{m^\pi}
\newcommand\pto{t^\pi}

\newcommand\vk{v^k}

\newcommand\uu{u}
\newcommand\ww{w}
\newcommand\pp{p}
\newcommand\qq{q}

\newcommand\plumpf{f}
\newcommand\nev{f}

\makeatletter
\gdef\centre#1{\smash{\vbox{\m@th\ialign{\hfil##\hfil\crcr
  \noalign{\kern3\p@}
  \hskip1pt{\tiny$\scriptscriptstyle\bullet$}\crcr
  \noalign{\kern2.5\p@\nointerlineskip}
  $\hfil\displaystyle{#1}\hfil$\crcr}}}}
\makeatother
\newcommand\pcentre[1]{\left(#1\right)^{\raise1pt\hbox{\tiny$\scriptscriptstyle\bullet$}}}

\newcommand\bbr{\mathbb{R}}

\newcommand\flr{\varphi}
\newcommand\phiset{Q}
\newcommand\phisetb{\hbox to1.75ex{\hss\hskip2pt$\overline{\hbox to1.5ex{\hss$Q$\hskip2pt\hss}}$\hss}}
\newcommand\lub{\hbox to1.75ex{\hss\hskip2pt$\overline{\hbox to1.5ex{\hss$F$\hskip2pt\hss}}$\hss}}
\newcommand\flrn{\flr^n}

\newcommand\sgn{\hbox{\textup{sgn}}}

\newcommand\ist{A}
\newcommand\xst{X}
\newcommand\yst{Y}
\newcommand\zst{Z}
\newcommand\wst{W}

\newcommand\tie{\hbox{\textit{\char"05}}}

\newcommand\istbis{\hbox to1.97ex{\hss\hskip3.5pt$\smash{\widetilde{%
 \hbox to1.9ex{\hss\vphantom{t}\smash{$\ist$}\hskip3.5pt\hss}}}$\hss}}
\newcommand\tiebis{\hbox to1.97ex{\hss\hskip3.5pt$\smash{\widetilde{%
 \hbox to1.9ex{\hss\vphantom{t}\smash{$\tie$}\hskip3.5pt\hss}}}$\hss}}
\newcommand\rhobis{\widetilde\rho}

\newcommand\cst{C}
\newcommand\clone{c}
\newcommand\clustit{\kern.1ex\widetilde c\kern.2ex}
\newcommand\contr[1]{\widetilde #1}

\newcommand\pre[1]{\mathsf{P}_{#1}}
\newcommand\suc[1]{\mathsf{S}_{#1}}
\newcommand\col[1]{\mathsf{C}_{#1}}

\newcommand\rank[1]{\kappa_{#1}}
\newcommand\ranka[2]{\kappa_{#1}(#2)}

\newcommand\rhoi{\eta}
\newcommand\last{\ell}
\newcommand\anteh{\,{}'\kern-.3ex h}
\newcommand\antef{\,{}'\kern-.45ex f}
\newcommand\antec{\,{}'\kern-.3ex c}

\newcommand\vbis{\widetilde v}
\newcommand\tbis{\widetilde t}
\newcommand\iscbis{\widetilde v{}^{\kern.5pt\ast}}
\newcommand\nubis{\widetilde\nu}
\newcommand\rbis{\widetilde r}
\newcommand\pscbis{\widetilde v{}^{\kern.75pt\pi}}

\newcommand\xibis{\smash{\widetilde{\hbox{\vphantom{t}\smash{$\xi$}}}}}
\newcommand\imgbis{\widetilde m^\nu}
\newcommand\ppmgbis{\widetilde m^\sigma}
\newcommand\pptobis{\widetilde t^\sigma}
\newcommand\ptobis{\widetilde t^\pi}
\newcommand\gammabis{\widetilde\gamma}
\newcommand\tref[1]{\smash{\hbox{$\widetilde{\hbox{\ref{#1}}}$}}}
\newcommand\trefbis[1]{\setbox0\hbox{\ref{#1}}%
\smash{\hbox to\wd0{\hss$\widetilde{\hbox to1em{\hss\ref{#1}\hss}}$\hss}}}
\newcommand\pmgbis{\widetilde m^\pi}

\newcommand\rev{\hbox{'}}
\newcommand\gpxy{\gamma\hbox{'}\kern-3pt{}_{xy}}
\newcommand\gbispxy{\widetilde\gamma\hbox{'}\kern-3pt{}_{xy}}

\newcommand\gto{\tau}

\newcommand\intg{\gamma}
\newcommand\inth{\eta}

\newcommand\uuu{\bm{u}}
\newcommand\aaa{\bm{a}}

\newcommand\vaa{\mathsf{V}}
\newcommand\vaasub{\vaa\kern-1pt}
\newcommand\vxx{\vaasub_{\scriptscriptstyle X\kern-1pt X}}

\newcommand\vxy{\vaasub_{\scriptscriptstyle X\kern-.25pt Y}}
\newcommand\vyx{\vaasub_{\scriptscriptstyle Y\kern-1pt X}}
\newcommand\vyy{\vaasub_{\scriptscriptstyle Y\kern-.25pt Y}}
\newcommand\vrs{\vaasub_{\scriptscriptstyle R\kern-.25pt S}}
\newcommand\vzz{\vaasub_{\scriptscriptstyle Z\kern-1pt Z}}
\newcommand\vxz{\vaasub_{\scriptscriptstyle X\kern-1pt Z}}
\newcommand\zeromatrix{\mathsf{O}}
\newcommand\irreopen{{\cal I}}
\newcommand\flrx{\flr_{\scriptscriptstyle X}}
\newcommand\flrax{\flr_{\scriptscriptstyle A\setminus X}}
\newcommand\flry{\flr_{\scriptscriptstyle Y}}
\newcommand\flrr{\flr_{\scriptscriptstyle R}}
\newcommand\gflr{\psi}
\newcommand\gflrx{\psi_{\scriptscriptstyle X}}
\newcommand\gflry{\psi_{\scriptscriptstyle Y}}

\newcommand\fxx{F_{\scriptscriptstyle X\kern-1pt X}}
\newcommand\fyy{F_{\scriptscriptstyle Y\kern-.25pt Y}}
\newcommand\fxy{F_{\scriptscriptstyle X\kern-.25pt Y}}
\newcommand\fzz{F_{\scriptscriptstyle Z\kern-1pt Z}}
\newcommand\frr{F_{\scriptscriptstyle R\kern-.25pt R}}
\newcommand\flrbis{\hbox to1.75ex{\hss\hskip1.25pt$\widetilde{\hbox to1.5ex{\hss$\varphi$\hskip1.25pt\hss}}$\hss}}
\newcommand\vaabis{\smash{\widetilde{\hbox{\vphantom{t}\smash{$\vaa$}}}}}
\newcommand\vaabissub{\vaabis\kern-1pt}
\newcommand\vaabisxx{\vaabissub_{\scriptscriptstyle X\kern-1pt X}}
\newcommand\vaabiszz{\vaabissub_{\scriptscriptstyle Z\kern-1pt Z}}
\newcommand\waa{\mathsf{W}}
\newcommand\waabis{\smash{\widetilde{\hbox{\vphantom{t}\smash{$\waa$}}}}}
\newcommand\xstbis{\hbox to1.97ex{\hss\hskip2.5pt$\smash{\widetilde{%
 \hbox to1.9ex{\hss\vphantom{t}\smash{$\xst$}\hskip2.5pt\hss}}}$\hss}}

\newcommand\xstbiss{\smash{\widetilde\xst}}
\newcommand\ystbiss{\smash{\widetilde\yst}}
\newcommand\zstbiss{\smash{\widetilde\zst}}
\newcommand\vxxbis{\vaabissub_{\scriptscriptstyle \xstbiss\kern-1pt\xstbiss}}
\newcommand\vyxbis{\vaabissub_{\scriptscriptstyle \ystbiss\kern-1pt\xstbiss}}
\newcommand\vzzbis{\vaabissub_{\scriptscriptstyle \zstbiss\kern-1pt\zstbiss}}
\newcommand\flrbissub{\flrbis\kern-1pt}
\newcommand\flrxbis{\flrbissub_{\scriptscriptstyle\xstbiss}}

\newcommand\flrbisX{\flrbissub_{\scriptscriptstyle\xst}}

\newcommand\flrbisx{\flrbissub_x}
\newcommand\flrbisy{\flrbissub_y}
\newcommand\flrbisz{\flrbissub_z}
\newcommand\flraxbis{\flrbissub_{\scriptscriptstyle A\setminus\xstbiss}}

\newcommand\xsth{\smash{\widehat X}}
\newcommand\ysth{\smash{\widehat Y}}
\newcommand\xh{\hat x}
\newcommand\yh{\hat y}

\newcommand\fbis{\hbox to1.97ex{\hss\hskip2.5pt$\smash{\widetilde{%
 \hbox to1.9ex{\hss\vphantom{t}\smash{$F$}\hskip2.5pt\hss}}}$\hss}}
\newcommand\fbisxx{\fbis_{\kern-2pt\scriptscriptstyle \xstbiss\kern-1pt\xstbiss}}

\newcommand\av[1]{\alpha_{#1}}

\newcommand\vbisbis{\smash{\widetilde{\hbox{\vphantom{t}\smash{$\vbis$}}}}\kern.5pt}
\newcommand\iscbisbis{\smash{\widetilde{\hbox{\vphantom{t}\smash{$\vbis$}}\kern.75pt}}^\ast}
\newcommand\nubisbis{\smash{\widetilde{\hbox{\vphantom{t}\smash{$\nubis$}}}}\kern.5pt}

\newcommand\chrel{\mathrel{\trianglerighteq}}
\newcommand\chrels{\mathrel{\equiv}}
\newcommand\chrela{\mathrel{\hbox{\relscale{1.25}$\kern1pt\triangleright\kern1pt$}}}
\newcommand\chrelbis{\mathrel{\smash{\widetilde{\hbox{\vrule width0pt height6.5pt\smash{$\trianglerighteq$}}}}}}
\newcommand\chrelabis{\mathrel{\smash{\widetilde{\hbox{\vrule width0pt height6.5pt\smash{\hbox{\relscale{1.25}$\kern1pt\triangleright\kern1pt$}}}}}}}

\newcommand\fun{{\cal F}}
\newcommand\funbis{{\cal F}^{\kern.5pt\prime}}
\newcommand\funter{{\cal G}}

\newtheorem{proposition}{Proposition}[section]
\newtheorem{lemma}[proposition]{Lemma}
\newtheorem{theorem}[proposition]{Theorem}
\newtheorem{corollary}[proposition]{Corollary}

\newlength\repskip 
\newenvironment{repeated}[1]%
{\null\hskip-\repskip\hbox to0.9\hsize\bgroup\hbox to0pt{\small$(\ref{#1})$\hss}\hfil\hskip.1\hsize$\displaystyle}%
{$\hfil\egroup}

\newenvironment{repeatedNN}[1]%
{\hbox to\hsize\bgroup\hbox to0pt{\small$(\ref{#1})$\hss}\hfil$\displaystyle}%
{$\hfil\egroup\nonumber}

\begin{document}


\thispagestyle{empty}

\null\vskip-5\baselineskip\null 
\renewcommand\upla{\vskip-2.5pt} 

\begin{center}
\hrule
\vskip7.5mm\upla
\textbf{\uppercase{A continuous rating method for~preferential~voting}}
\par\medskip
\textsc{Rosa Camps,\, Xavier Mora \textup{and} Laia Saumell}
\par
Departament de Matem{\`a}tiques,
Universitat Aut\`onoma de Barcelona,
Catalonia,
Spain
\par\medskip
\texttt{xmora\,@\,mat.uab.cat}
\par\medskip
29th July 2008,\ \  revised 26th November 2008
\vskip7.5mm\upla
\hrule
\end{center}

\upla\vskip-4mm\null
\begin{abstract}
A method is given for quantitatively rating the social acceptance of different options which are the matter of a preferential vote.\linebreak
The~proposed method is proved to satisfy certain desirable conditions, among which there is a majority principle, a property of clone consistency, and the continuity of the rates with respect to the data.\linebreak
One can view this method as a quantitative complement for a qualitative method introduced in 1997 by Markus Schulze. It~is also related to certain methods of one-dimensional scaling or cluster analysis.

\vskip2pt\upla

\upla
\bigskip\noindent
\textbf{AMS subject classifications:} 
\textit{%
05C20, 
91B12, 
91B14, 
91C15, 
91C20. 
}
\end{abstract}

\upla
\vskip5mm
\hrule

\upla
\vskip-12mm\null 
\section*{}
The outcome of a vote is often expected to entail a quantitative
rating of the candidate options according to their social acceptance.
Some voting\linebreak[3] methods are directly based upon such a rating.
\ensep
This is the case when each voter is asked to choose one option \,and each option is rated by the fraction of the vote in its favour. The resulting rates can be used for filling a single seat (first past the post) or for distributing a number of them (proportional representation).
\ensep
A more elaborate voting method based upon quantitative rates was introduced in 1433 by Nikolaus von Kues~\cite[\secpar{1.4.3}, \secpar{4}]{mu}
and again in 1770--1784 by Jean-Charles de~Borda~\cite[\secpar{1.5.2}, \secpar{5}]{mu}.
Here, each voter is asked to rank the different options in~order of preference \,and each option is rated by the average of its ranks, \ie the~ordinal numbers that give its position in these different rankings (this formulation differs from the traditional one by a linear function).
\ensep
For future reference in this paper, 
these two rating methods will be called respectively\,
the method of \dfc{first-choice fractions}\,
and\, the method of \dfc{average ranks}.

However, both of these methods have important drawbacks, which leads to the point of view of \dfc{paired comparisons}, where each option is confronted with every other by counting how many voters prefer the former to the latter and vice versa. From this point of view it is quite natural to~abide by the so-called Condorcet principle: an option should be deemed the winner whenever it defeats every other one in this sort of tournament. This approach was introduced as early as in the thirteenth century by Ramon Llull \cite[\secpar{1.4.2}, \secpar{3}]{mu}, and later on it was propounded again by the marquis of Condorcet in 1785--1794~\cite[\secpar{1.5.4}, \secpar{7}]{mu}, and by Charles Dodgson, alias Lewis Carroll, in 1873--1876
\hbox{\brwrap{
\,\textsl{2}\,\textup{;}\space
\textsl{23}\,\textup{:\,\relscale{0.95}\secpar{12}}
}}.
\ensep
Its~development gives rise to a variety of methods, 
some of them with remarkably good properties.
This is particularly the case of 
the method of \dfd{ranked~pairs}, proposed in 1986/87 by~Thomas M.~Zavist and T.~Nicolaus Tideman~\cite{ti, zt},
and the method introduced in~1997 by~Markus Schulze~\cite{sc, scbis}, which we will refer to as the method of \dfd{paths}.
In~spite of the fact that generally speaking they can produce different results, both of them comply with the Condorcet principle
and they share the remarkable property of clone consistency~\cite{t6, sc}.

Nevertheless, these methods do not immediately give
a quantitative rating of the candidate options.
Instead, they are defined only as algorithms for determining a winner or at most a purely ordinal ranking.
On the other hand, they are still based upon the quantitative information provided by the table of paired-comparison scores,
which raises the question of whether their qualitative results can be consistently converted into quantitative ratings.

In~\cite[\secpar{10}]{xm}
a quantitative rating algorithm was devised with the aim of complementing the method of ranked pairs.
Although a strong evidence was given for its fulfilling
certain desirable conditions
---like the ones stated below---,
it was also pointed out that it
fails a most natural one,
namely that the output, \ie the rating,
be a continuous function of the input,\linebreak 
\ie the frequency of each possible
content of an individual vote.
In~fact, such a lack of continuity seems
unavoidable when the method of ranked pairs is considered
and those other conditions are imposed.
\ensep
In~contrast, in this paper we will see that the method of paths
does admit such a continuous rating procedure.

Our method can be viewed as a projection of the matrix of paired-comparison scores onto a special set of such matrices. This projection is combined with a subsequent application of two standard rating methods, one of which the method of average ranks. The overall idea has some points in common with \cite{sm}.

\pagebreak 

We will refer to the method described in this paper
as the \dfd{CLC~rating method}, where the capital letters stand for ``Continuous Llull Condorcet''.

\bigskip
The paper is organized as follows: In section~1 we state the problem which is to be solved and we make some general remarks. Section~2 presents an heuristic outline of the proposed method.
Section~3 gives a summary of the procedure, after which certain variants are introduced. Section~4 presents some illustrative examples. Finally, sections~5--18 give detailed mathematical proofs of the claimed properties for the main variant.

\bigskip
The reader interested to try the CLC method can make use of the tool  which is available at\ensep {\small\texttt{http://mat.uab.cat/{\atilde}xmora/CLC\underline{ }calculator/}}.

\section{Statement of the problem and general remarks}

\paragraph{1.1}
We consider a set of $N$~options
which are the matter of a vote.
Although more general cases will be included later on~(\secpar{3.3}), for the moment we assume that each voter expresses his preferences in the form of a ranking; by~it we mean an ordering of the options in question by decreasing degree of preference, with the possibility of ties
and/or truncation (\ie expressing a top segment only).
We~want to aggregate these individual preferences into a \emph{social rating},
where each option is assigned a rate
that quantifies its social acceptance.

In some places we will restrict our attention to the case of complete votes. For ranking votes, we are in such a situation whenever we are dealing with non-truncated rankings.
As we will see, the incomplete case will give us much more work than the complete one.


We will consider two kinds of ratings, which will be referred to respectively as rank-like ratings and fraction-like ones. As it is suggested by these names, a rank-like rating will be reminiscent of a ranking, whereas a fraction-like one will evoke the notion of proportional representation. Our method will produce both a rank-like rating and a fraction-like one. They will agree with each other in the ordering of the candidate options, except that the ordering implied by the fraction-like rating may be restricted to a top segment of the other one.
Quantitatively speaking, the two ratings have different meanings. In particular, the fraction-like rates can be viewed as an estimate of the first-choice fractions based not only on the first choices of the voters, but also on the whole set of preferences expressed by them. In~contrast, the rank-like rates are not focused on choosing, but they aim simply at positioning all the candidate options on a certain scale.

More specifically, the two ratings are asked to satisfy the following conditions:

\smallskip
\newcommand\llsi{\textup{A}}
\condition{\llsi}{Scale invariance}
The rates depend only on the relative frequency of each possible content of an individual vote.
In~other words, if every individual vote is replaced by a fixed number of copies of it, the rates remain exactly the same.

\smallskip
\newcommand\llpe{\textup{B}}
\condition{\llpe}{Permutation equivariance}
Applying a certain permutation of the options to all of the individual votes has no other effect than getting the same permutation in the social rating.

\smallskip
\newcommand\llco{\textup{C}}
\condition{\llco}{Continuity}
The rates depend continuously on the relative frequency of each possible content of an individual vote.

\medskip
\noindent
The next conditions apply to the rank-like rating:

\newcommand\condr[1]{#1\rlap{\mathsurround0pt$_r$}}
\newcommand\condf[1]{#1\rlap{\mathsurround0pt$_f$}}

\smallskip
\newcommand\llrr{\textup{D}}
\condition{\llrr}{Rank-like range}
Each rank-like rate is a number, integer or fractional, between $1$ and~$N$. The~best possible value is~$1$ and the worst possible one is~$N$.

\smallskip
\newcommand\llrd{\textup{E}}
\condition{\llrd}{Rank-like decomposition}
Let us restrict the attention to the complete case.
Consider a splitting of the options into a~`top~class'~$\xst$ plus a `low~class'~$\yst$.
Assume that all of the voters have put each member of $\xst$ above every member of~$\yst$.\ensep
In~that case, and only in that~case, the rank-like rates can be obtained separately for each of these two classes according to the corresponding restriction of the ranking votes (with the proviso that the unassembled low-class rates differ from the assembled
ones by the number of top-class members).

\smallskip
\noindent
In its turn, the fraction-like rating is required to satisfy the following conditions:

\smallskip
\newcommand\llfc{\textup{F}}
\condition{\llfc}{Fraction-like character}
Each fraction-like rate is a number greater than or equal to $0$.
Their sum is equal to a fixed value.
More specifically, we will take this value to be the participation fraction, \ie the fraction of non-empty votes.

\smallskip
\newcommand\llfd{\textup{G}}
\label{ultimaprop}
\condition{\llfd}{Fraction-like decomposition}
Consider the same situation as in \llrd\ with the additional assumption that there is no proper subset of~$\xst$ with the same splitting property as~$\xst$ (namely, that all voters have put each option from that set above every one outside it).\ensep
In~that case, and only in that~case, the top-class fraction-like rates are all of them positive and they can be obtained according to the corresponding restriction of the ranking votes, whereas the low-class fraction-like rates are all of them equal to~$0$.

\smallskip
\newcommand\llpv{\textup{H}}
\condition{\llpv}{Case of plumping votes}
Assume that each voter plumps for a single option. In that case, the fraction-like rates coincide with the fractions of the vote obtained by each option.

\renewcommand\condr[1]{$\mathrm#1_r$}
\renewcommand\condf[1]{$\mathrm#1_f$}

\medskip
\noindent
Furthermore, we ask for some properties that concern only
the concomitant social ranking, \ie
the purely ordinal information contained in the social rating:

\smallskip
\newcommand\llmp{\textup{I}}
\condition{\llmp}{Majority principle}
Consider a splitting of the options into a~`top class'~$X$ plus a `low~class'~$Y$. Assume that for each member of $X$ and every member of $Y$ there are more than half of the individual votes where the former is preferred to the latter.
In that case, the social ranking also prefers each member of $X$ to every member of~$Y$.

\medskip
\newcommand\llcc{\textup{J}}
\condition{\llcc}{Clone consistency}\!
A set \!$\cst$\! of options is said to be a \emph{cluster} (of~\emph{clones}) for a given ranking \,when\, each element from outside~$\cst$ compares with all elements of~$\cst$ in the same way
(\ie either it lies above all of them, or it lies below all of them, or it ties with all of them).\linebreak[3]
In~this connection, it is required that if a set of options is a cluster for each of the individual votes, then: (a) it is a cluster for the social ranking;\, and\, (b)~contracting it to a single option in all of the individual votes has no other effect in the social ranking than getting the same contraction.

\paragraph{1.2}
Let us emphasize that the individual votes that we are dealing with do \textit{not} have a quantitative character (at least for the moment): each voter is allowed to express a preference for $x$ rather than $y$, or vice versa, or maybe a tie between them, but he is not allowed to quantify such a preference.

This contrasts with ‘range voting’ methods, where each individual vote is already a quantitative rating \cite{rv, bala}.
Such methods are free from many of the difficulties that lurk behind the present setting.
However, they make sense only as long as all voters mean the same by  each possible value of the rating variable.
This hypothesis may be reasonable in some cases,
but quite often it is hardly applicable
(a typical symptom of its not being appropriate
is a concentration of the rates in a small set
independently of which particular options are under consideration).
In such cases, it is quite natural that the individual votes 
express only qualitative comparisons between pairs of options.
If the issue is not too complicated, one can expect these comparisons to form a ranking.
In the own words of\, \hbox{\brwrap{\textsl{\kern.4pt 1\kern.4pt a\kern.4pt}}},
“When there is no common language, a judge's only meaningful input is the order of his grades”.
\ensep
Certainly, the judges will agree upon the qualitative comparison between two options much more often than they will agree upon their respective rates in a certain scale. Such a lack of quantitative agreement may be due to truly different opinions; but quite often it is rather  meaningless.
Of course, the rates will coincide more easily if a discrete scale of few grades is used. But then it may happen that the judges rate equally two options about which they all share
a definite preference for one over the other, in which case these discrete rates are throwing away genuine information.
\ensep
Anyway, voting is often used in connection with moral, psychological or aesthetic qualities, whose appreciation may be as little  quantifiable, but also as much “comparable”, as, for instance, the feelings of pleasure or pain.



\medskip
So, in our case the quantitative character of the output is not present in~the individual votes (unless we adopt the general setting considered at the end of~\secpar{3.3}), but it derives from the fact
of having a number of them.
The~larger this number, the more meaningful is the quantitative character of the social rating.
This is especially applicable to the continuity property~\llco,
according to which 
a small variation in the proportion of votes with a given content produces only small variations in the rates.
In fact, if all individual votes have the same weight,
a few votes will be 
a small proportion only
in~the measure that 
the total number of votes is large enough.

In this connection, it~should be noticed
that property~\llco\ differs from the continuity
property adopted in~\cite{bala} (axiom~6), which
does not refer to small variations in the proportion of votes with a given content, but to small variations in the quantitative content of each individual vote. In the general setting considered at the end of~\secpar{3.3}, the~CLC method satisfies not only the continuity property~\llco, but also the axiom~6 of~\cite{bala}; in~contrast, the “majority-grade” method considered in~\cite{bala} satisfies the latter but not the former.

\paragraph{1.3}
One can easily see that the method of average ranks satisfies conditions \,\llsi--\llrd. In principle that method assumes that all of the individual votes are complete rankings; however, one can extend it to the general case of rankings with ties and/or truncation while keeping those conditions (it suffices to use formula~(\ref{eq:avranksfromscores}) of~\secpar{2.5}).
\ensep
In~their turn, the first-choice fractions are easily seen to satisfy conditions\linebreak[3] \llsi--\llco\, and \llfc--\llpv.
\ensep
However, neither of these two methods satisfies conditions \llmp\ and \llcc.
In fact, these conditions were introduced precisely as
particularly desirable properties that are not satisfied by those methods \cite{mu, bl, t6}.

Of~course, one can go for a particular ranking method
that satisfies conditions~\llmp\ and \llcc\ and then look for an appropriate algorithm
to convert the ranking result into the desired rating
according to the quantitative information coming from the vote.
But this should be done in such a way that the final rating
be always
in agreement with the ranking method
as well as 
in compliance with conditions~\llsi--\llpv,
which is not so easy to achieve.
\ensep
From this point of view, our proposal can be viewed
as providing such a complement for 
one of the variants of the method of paths \cite{sc, scbis}.

\newcommand\uplapar{\vskip-6.5mm\null} 

\uplapar
\paragraph{1.4}
When the set $\xst$ consists of a single option, the~majority principle~\llmp\,\ takes the following form:

\renewcommand\upla{\vskip-2pt}

\upla
\smallskip
\newcommand\llmpw{\textup{I1}}
\condition{\llmpw}{Majority principle, winner form} If an option $x$ has the property that for every $y\neq x$ there are more than half of the individual votes where $x$ is preferred to $y$, then $x$ is the social winner.

\upla
\smallskip
\noindent
In the complete case the preceding condition is equivalent to the following one:

\upla
\smallskip
\newcommand\llcpw{\textup{I1$'$}}
\condition{\llcpw}{Condorcet principle} If an option $x$ has the property that\linebreak 
for~every $y\neq x$ there are more individual votes where $x$ is preferred to $y$ than vice versa, then $x$ is the social winner.

\smallskip
\noindent
However, we want to admit the possibility 
of individual votes where no information is given
about certain pairs of options.
For instance, in the case of a truncated ranking it makes sense to interpret that there is no information about two particular options which are not present in the list.
In that case condition~\llmpw\ is weaker than~\llcpw,
and the CLC~method will satisfy only the weaker version.

\smallskip
This~lack of~compliance with the Condorcet principle
and its being replaced by a weaker condition  
may be considered undesirable.
\ensep
However, other authors have already remarked that
such a weakening of the Condorcet principle
is necessary in order to be able to keep other properties \cite{wo} (see also~\cite{kl}).
In our case, Condorcet principle seems to conflict with the continuity property~\llco\ (see \secpar{3.3}).
\ensep
On the other hand, the Condorcet principle was originally proposed in connection with the complete case~\cite{mu}, its generalization in the form~\llcpw\ instead of~\llmpw\ being due to later authors.
Even so, nowadays it is a common practice to refer to~\llcpw\ by the name of “Condorcet principle’.

\uplapar
\paragraph{1.5}
As we mentioned in the preceding subsection, we want to admit the possibility of individual votes where no information is given about certain pairs of options.
In this connection, the CLC~method will carefully distinguish a~{\em definite indifference} about two or more options from a~{\em lack of information} about them (see~\cite{he}).
For~instance, if all of the individual votes are complete rankings but they balance into an exact social indifference ---in particular if~each individual vote expresses such a complete indifference---, the resulting rank-like rates will be all of them equal to~$(N\!+\!1)/2$ and the corresponding fraction-like rates will be equal to~$1/N$.
In contrast, in the case of a full abstention, \ie where no voter has expressed any opinion, the rank-like rates will be all of them equal to~$N$ and the corresponding fraction-like rates will be equal to~$0$.


\medskip
Although the decomposition conditions~\llrd\ and \llfd\ have been stated only for the complete case, some partial results of that sort will hold under more general conditions. In particular, the following condition will be satisfied for general, possibly incomplete, ranking votes: The winner will be rated exactly~$1$ (in~both the rank-like rating and the fraction-like one) if and \emph{only~if} all of~the voters have put that option into first place.

\begin{comment}
En canvi, el següent és fals: An option will get a rank-like rate  exactly equal to~$N$ \ifoi all of the votes are either a complete ranking with that option in the last place or a truncated ranking which leaves that option out. 
Exemple contra ‘if’: 1~$a\succ b\succ c\succ d$, 1~$a\succ c$, 1~$a$ (el resultat és $1.0, 3.25, 3.25, 3.8333$). ‘Only if’: no ho sé demostrar.\par
\end{comment}

\smallskip
Conditions \llrd\ and \llfd, as well as the preceding property, refer to cases where ``all of the voters'' proceed in a certain way.
Of course, it should be clear whether we mean all of the ``actual'' voters or maybe all of the ``potential'' ones (\ie actual voters plus abstainers).
We assume that one has made a choice in that connection,
thus defining a total number of voters~$V$.
Considering all potential voters instead of only the actual ones 
has no other effect than contracting
the final rating towards the point where all rates take the minimal value
(namely, $N$~for rank-like rates \,and\, $0$~for fraction-like ones).

\uplapar
\paragraph{1.6}
It is interesting to look at the results of the CLC~rating method when it is applied to the approval voting situation, \ie the case where each voter gives only a list of approved options, without any expression of preference between them. In~such a situation it is quite natural to rate each option by the number of received approvals; the resulting method has pretty good properties, not the least of which is its eminent simplicity~\cite{br}.

\begin{comment}
Notice that approval voting is the limit case of range voting where the “quantitative” rating scale reduces to two values only.
\end{comment}

Now, an individual vote of approval type can be viewed as a truncated ranking which ties up all of the options that appear in it.
So it makes sense to apply the CLC~rating method.
Quite remarkably, one of its variants turns out to order the options  in exactly the same way as the number of received approvals (see \secpar{17}). More specifically, the variant in~question corresponds to interpreting that the non-approved options of an~individual vote are tied to each other. However, the main variant, which acknowledges 
a lack of comparison between non-approved options, can lead to different results.

\uplapar
\section{Heuristic outline}

This section presents our proposal as the result of a quest 
for the desired properties. Hopefully, this will communicate the main ideas that lie behind the formulas.

\paragraph{2.1}
The aim of complying with conditions~\llmp\ and \llcc\ 
calls for the point of view of \dfd{paired comparisons}.
In accordance with it,
our procedure will be based upon considering every pair of options 
and counting how many voters prefer one to the other or vice versa.
\ensep
To that effect, we must adopt some rules for translating the ranking votes
(possibly truncated or with ties) into binary preferences.
In principle, these rules will be the following:

\iim{a}When $x$ and $y$ are both in the list\,
and $x$ is ranked above $y$ (without a tie),
we certainly interpret that $x$ is preferred to $y$.

\iim{b}When $x$ and $y$ are both in the list\,
and $x$ is ranked as good as $y$,\,
we interpret it as being equivalent to
half a vote preferring $x$ to $y$\,
plus another half a vote preferring~$y$~to~$x$.

\iim{c}When $x$ is in the list and $y$ is not in it,\,
we interpret that $x$ is preferred to $y$.

\iim{d}When neither $x$ nor $y$ are in the list,\,
we interpret nothing
about the preference of the voter between $x$ and $y$.

\noindent
Later on~(\secpar{3.2,\,3.3}) we will consider certain alternatives to rules~(d) and (c).

\medskip
The preceding rules allow us to count how many voters
support a given binary preference,
\ie a~particular statement of the form ``$x$~is~preferable~to~$y$''.
\ensep
By~doing so for each possible pair of options $x$~and~$y$,
the whole vote gets summarized into a set of $N(N-1)$~numbers
(since $x$~and $y$ must be different from each other).
We will denote these numbers by~$V_{xy}$
and we will call them the binary \dfc{scores} of the vote.
The collection of these numbers will be called the \dfc{Llull matrix} of the vote.
\ensep
Since we look for scale invariance,
it makes sense to divide all of these numbers 
by the total number of votes~$V$,
which normalizes them to range from~$0$ to~$1$.
In the following we will work mostly with these normalized scores,
which will be denoted by~$v_{xy}$.
In practice, however, the absolute scores~$V_{xy}$ have the advantage
that they are integer numbers, so we will use them in the examples.

\medskip
In general, the numbers $V_{xy}$ are bound to satisfy $V_{xy} + V_{yx} \le V$, or~equivalently $v_{xy} + v_{yx} \le 1$.
\ensep
The special case where the ranking votes are all of them complete, \ie without truncation,
is characterized by the condition that $V_{xy} + V_{yx} = V$,
or~equivalently $v_{xy} + v_{yx} = 1$.
From now on we will refer to such a situation as the case of \dfc{complete votes}.

\medskip
Besides the scores $v_{xy}$, in the sequel we will often deal with the
\dfc{margins} $m_{xy}$ and the \dfc{turnovers}~$t_{xy}$,
which are defined respectively by
\begin{equation}
m_{xy} \,=\, \hbox to24mm{$v_{xy} - v_{yx},$\hfil}\qquad
t_{xy} \,=\, \hbox to24mm{$v_{xy} + v_{yx}.$\hfil}
\end{equation}
Obviously, their dependence on the pair $xy$ is respectively antisymmetric and symmetric, that is
\begin{equation}
m_{yx} \,=\, \hbox to24mm{$- m_{xy},$\hfil}\qquad
t_{yx} \,=\, \hbox to24mm{$t_{xy}.$\hfil}
\end{equation}
It is clear also that the scores $v_{xy}$ and $v_{yx}$
can be recovered from $m_{xy}$ and $t_{xy}$ by means of the formulas
\begin{equation}
v_{xy} \,=\, \hbox to24mm{$(t_{xy} + m_{xy})/2,$\hfil}\qquad
v_{yx} \,=\, \hbox to24mm{$(t_{xy} - m_{xy})/2.$\hfil}
\end{equation}

\paragraph{2.2}
A natural candidate for defining the social preference is the following: $x$~is socially preferred to~$y$
whenever $v_{xy} > v_{yx}$.
Of course, it can happen that $v_{xy} = v_{yx}$,
in~which case one would consider that $x$~is socially equivalent to~$y$.
The~binary relation that includes all pairs $xy$ for which $v_{xy} > 
v_{yx}$ will be denoted by $\mu(v)$ and will be called the
\dfc{comparison relation};
together with it, we will consider also the
\dfc{adjoint comparison relation} $\hat\mu(v)$ which is defined by the condition
$v_{xy} \ge v_{yx}$.

\medskip
As it is well-known, the main problem with paired comparisons
is that the comparison relations $\mu(v)$ and~$\hat\mu(v)$
may lack transitivity even if the individual
preferences are all of them transitive~\cite{mu, bl}.
More specifically, $\mu(v)$ can contain a ‘Condorcet~cycle’, \ie a sequence $x_0 x_1 \dots x_n$ such that $x_n=x_0$ and $x_ix_{i+1}\in\mu(v)$ for all $i$.

A most natural reaction to it is going for the transitive closure
of~$\hat\mu(v)$, which we will denote by~$\hat\mu^\ast(v)$.
By definition, $\hat\mu^\ast(v)$ includes all (ordered) pairs $xy$
for which there is a path $x_0 x_1 \dots x_n$
from~$x_0 = x$ to~$x_n = y$ whose links $x_ix_{i+1}$
are all of them in $\hat\mu(v)$.
In~other words, we can say that $\hat\mu^\ast(v)$ includes all pairs that are
``indirectly related'' through $\hat\mu(v)$.
\ensep
However, this operation replaces each cycle of intransivity by an equivalence between its members. Instead of that, we would rather break these equivalences
according to the quantitative information provided by the scores $v_{xy}$.
This is what is done in such methods as ranked pairs or paths.
However, these methods use that quantitative information
to reach only qualitative results.
In~contrast, our results will keep a quantitative character until the end.

\paragraph{2.3}
The next developments rely upon an operation $(v_{xy})\rightarrow (\isc_{xy})$
that transforms the original binary scores into a new one.
This operation is defined in the following way:
for every pair~$xy$, one considers all possible paths $x_0 x_1 \dots x_n$
going from~$x_0 = x$ to~$x_n = y$; every such path is associated with the
score of its weakest link, \ie the smallest value of $v_{x_ix_{i+1}}$;
finally, $\isc_{xy}$ is defined as the maximum value of this associated score
over all paths from $x$ to $y$.
In other words,
\begin{equation}
\isc_{xy} \hskip.75em = \hskip.75em
\max_{\vtop{\scriptsize\halign{\hfil#\hfil\cr\noalign{\vskip.5pt}$x_0=x$\cr$x_n=y$\cr}}}
\hskip.75em
\min_{\vtop{\scriptsize\halign{\hfil#\hfil\cr\noalign{\vskip-1.25pt}$i\ge0$\cr$i<n$\cr}}}
\hskip.75em v_{x_ix_{i+1}},
\label{eq:paths}
\end{equation}
where the \,$\max$\, operator considers all possible paths from $x$ to $y$,
and the \,$\min$\, operator considers all the links of a particular path.
The scores $\isc_{xy}$ will be called the \dfc{indirect scores} 
associated with the (direct) scores $v_{xy}$.

If $(v_{xy})$ is the table of 0's and 1's associated with a binary relation $\rho$
(by putting $v_{xy}=1$ \ifoi $xy\in\rho$), then $(\isc_{xy})$
is exactly the table associated with $\rho^\ast$, the transitive closure of $\rho$.
So, the operation $(v_{xy})\mapsto (\isc_{xy})$ can be viewed as a
quantitative analog of the notion of transitive
closure (see~\cite[Ch.\,25]{co}).

The main point, remarked in 1998 by Markus Schulze~%
\hbox{\brwrap{\textsl{34\,b}}}, 
is that the comparison relation
associated with a table of indirect scores is always transitive (Theorem~\ref{st:transSchulze}).
So, $\mu(\isc)$ is
always transitive, no matter what the case is for~$\mu(v)$. This is true in spite of the fact that
$\mu(\isc)$ can easily differ from~$\mu^\ast(v)$.
In the following we will refer to~$\mu(\isc)$ as the \dfd{indirect comparison relation}.

\remark

Somewhat surprisingly, in the case of incomplete votes the transitive relation $\mu(\isc)$ may differ from~$\mu(v)$ even when the latter is already transitive. 
\ensep
An example is given by the following profile, where each indicated preference is preceded by the number of people who voted in that way: 17~$a$, 24~$c$, 16~$a\succ b\succ c$, 16~$b\succ a\succ c$, 8~$b\succ c\succ a$, 8~$c\succ b\succ a$; in this case  the direct comparison relation $\mu(v)$ is the ranking $a\succ b\succ c$, whereas the indirect comparison relation $\mu(\isc)$ is the ranking $b\succ a\succ c$.
\ensep
More specifically, we have $V_{ab} = 33 > 32 = V_{ba}$ but $V^\ast_{ab} = 33 < 40 = V^\ast_{ba}$.

The agreement with $\mu(v)$ can be forced by suitably redefining the indirect scores; more specifically, formula (\ref{eq:paths}) can be replaced by an analogous one where the \,$\max$\, operator is not concerned with all possible paths from $x$ to $y$ but only those contained in $\mu(v)$. This idea is put forward in~\cite{scbis}. Generally speaking, however, such a method cannot be made into a continuous rating procedure since one does quite different things depending on whether $v_{xy} > v_{yx}$ or $v_{xy} < v_{yx}$. On the other hand, we will see that in the complete case the indirect comparison relation does not change when the paths are restricted to be contained in $\mu(v)$ (\secpar{7}).

\paragraph{2.4}
In the following we put
\begin{equation}
\nu=\mu(\isc),\qquad \hat\nu=\hat\mu(\isc),\qquad
\img_{xy} = \isc_{xy} - \isc_{yx}.
\end{equation}
So, $xy\in\nu$ \ifoi $\isc_{xy} > \isc_{yx}$, \ie $\img_{xy}>0$,
and $xy\in\hat\nu$ \ifoi $\isc_{xy} \ge \isc_{yx}$, \ie $\img_{xy}\ge0$.
From now on we will refer to $\img_{xy}$ as the {\df indirect margin} associated with the pair $xy$.

As it has been stated above,
the relation~$\nu$ is transitive.\ensep
Besides that, it~is clearly antisymmetric
(one cannot have 
both $\isc_{xy} > \isc_{yx}$ and vice versa).
On~the other hand, it may be not complete
(one can have $\isc_{xy} = \isc_{yx}$).
When it differs from~$\nu$, the complete relation $\hat\nu$
is not antisymmetric and ---somewhat surprisingly---
it may be not transitive either.
\ensep
For~instance, consider the profile given by 4~$b\succ a\succ c$, 3~$a\succ c\succ b$, 2~$c\succ b\succ a$, 1~$c\succ a\succ b$; in this case the indirect comparison relation $\nu = \mu(\isc)$ contains only the pair $ac$; as a consequence, $\hat\nu$ contains $cb$ and $ba$ but not $ca$.
\ensep
However, one can always find a total order~$\xi$ which
satisfies $\nu \sbseteq \xi \sbseteq \hat\nu$ (Theorem~\ref{st:existenceXiThm}).
From now on, any total order~$\xi$ that satisfies this condition 
will be called an \dfc{admissible order}.

The rating that we are looking for will be based on such an order~$\xi$. 
More specifically, it will be compatible with $\xi$ in the sense that the rates~$r_x$
will satisfy the inequality $r_x \le r_y$ whenever $xy\in\xi$.
If $\nu$ is already a total order, 
so that $\xi=\nu$,
the preceding inequality will be satisfied in the strict form $r_x<r_y$,
and this will happen \ifoi $xy\in\nu$ (Theorem~\ref{st:RvsNu}).

If there is more than one admissible order 
then some options will have equal rates.
In fact, we will have $r_x = r_y$ whenever $xy\in\xi_1$ and $yx\in\xi_2$,
where $\xi_1,\xi_2$ are two admissible orders. 
This will be so because we want the rating to be independent of the choice of~$\xi$.
This independence with respect to~$\xi$ seems essential for achieving the continuity property~C;
in fact, each possible choice of~$\xi$ for a given profile of vote frequencies may easily become the only one for a slight perturbation of that profile
(but not necessarily,
as it is illustrated by example~10 of \cite[\secpar{4.6}]{scbis}).

The following steps assume that one has fixed an admissible order~$\xi$.
From now on the situation~$xy\in\xi$ will be expressed also by~$x\rxi y$.
According to the definitions, the inclusions $\nu \cd\sbseteq \xi \cd\sbseteq \hat\nu$
are~equivalent to saying that $\isc_{xy}>\isc_{yx}$ implies $x\rxi y$
and that the latter implies $\isc_{xy} \ge \isc_{yx}$.
In other words, if the different options are ordered according to $\rxi$,
the~matrix $\isc_{xy}$ has then the property that each element above the diagonal
is larger than or equal to its symmetric over the diagonal.

\paragraph{2.5}
Rating the different options means positioning them on a line.
Besides\linebreak[3] complying with the qualitative restriction of
being compatible with~$\xi$ in the sense above,
we want that the distances between items
reflect the quantitative information provided by the binary scores.
However, a rating is expressed by $N$~numbers,
whereas the binary scores are $N(N-1)$ numbers.
So we are bound to do some sort of projection.
Problems of this kind have a certain tradition in 
combinatorial data analysis and cluster analysis
\cite{js, mi, ham}. 
In~fact, some of the operations that will be used below
can be viewed from that point of view.

\medskip
Let us assume for a while that we are dealing with complete ranking votes,
so that it makes sense to talk about the average ranks.
It is well-known~\cite[Ch.\,9]{bl} 
that their values, which we will denote by~$\bar r_x$, 
can be obtained from the Llull matrix by means of the following formula:
\begin{equation}
\bar r_x \,=\, N - \sum_{y\neq x}\,v_{xy}.
\label{eq:avranksfromscores}
\end{equation}
Equivalently, we can write
\begin{equation}
\bar r_x \,=\, (N+1 - \sum_{y\neq x}\,m_{xy}\,) \,/\, 2,
\label{eq:avranksfrommargins}
\end{equation}
where the $m_{xy}$ are the margins of the original scores $v_{xy}$,
\ie $m_{xy} = v_{xy}-v_{yx}$.
In fact, the hypothesis of complete votes means that $v_{xy} + v_{yx}\,=\, 1$,
so that $m_{xy} \,=\, 2 v_{xy} - 1$, which gives the equivalence between (\ref{eq:avranksfromscores}) and (\ref{eq:avranksfrommargins}).

Let us look at the meaning of the margins~$m_{xy}$ in connection with the idea of projecting
the Llull matrix into a rating: If there are no other items than $x$~and~$y$,
we can certainly view the sign and magnitude of~$m_{xy}$ as giving respectively the qualitative and quantitative aspects of the relative positions of $x$ and $y$ on the rating line,
that~is, the order and the distance between them.
\ensep
When there are more than two items, however, 
we have several pieces of information of this kind, one for every pair,
and these different pieces may be incompatible with each other,
quantitatively or even qualitatively,
which motivates indeed the problem that we are dealing with.
In particular, the~average ranks often violate
the desired compatibility with the relation~$\xi$.

\medskip
In order to construct a rating compatible with~$\xi$, we will
use a formula analogous to~(\ref{eq:avranksfromscores}) where 
the scores~$v_{xy}$ are replaced by certain \dfc{projected scores}~$\psc_{xy}$ to be defined in the following paragraphs.
Together with them, we will make use of the corresponding \dfc{projected margins}~$\pmg_{xy} = \psc_{xy} - \psc_{yx}$
and the corresponding \dfc{projected turnovers}~$\pto_{xy} = \psc_{xy} + \psc_{yx}$.
So,~the rank-like rates that we are looking for will be obtained in the following way:
\begin{equation}
r_x \,=\, N - \sum_{y\neq x}\,\psc_{xy}.
\label{eq:rrates}
\end{equation}
This formula will be used not only in the case of complete ranking votes,
but also in the general case where the votes are allowed to be incomplete
and/or intransitive binary relations.

\paragraph{2.6}
Let us begin by the case of complete votes, \ie $t_{xy} = v_{xy} + v_{yx} = 1$.
In~this case, we put also $\pto_{xy} = 1$. Analogously to (\ref{eq:avranksfromscores}) and (\ref{eq:avranksfrommargins}), 
formula~(\ref{eq:rrates}) is then equivalent to the following one:
\begin{equation}
r_x \,=\, (N+1 - \sum_{y\neq x}\,\pmg_{xy}\,) \,/\, 2.
\label{eq:rratesfrommargins}
\end{equation}

We want to define the projected margins $\pmg_{xy}$ so that
the rating defined by~(\ref{eq:rratesfrommargins}) be compatible with the ranking~$\xi$, \ie $xy\in\xi$ implies $r_x\le r_y$.
Now,~this ranking derives from the relation~$\nu$,
which is concerned with the sign of $\img_{xy} = \isc_{xy}-\isc_{yx}$.
This clearly points towards taking $\pmg_{xy} = \img_{xy}$.
However, this is still not enough for ensuring the compatibility with $\xi$.\linebreak[3]
In~order to ensure this property, it suffices that the projected margins, which we assume antisymmetric, behave in the following way:
\setbox0=\vbox{\hsize105mm\noindent
$x\rxi y \,\,\Longrightarrow\,\, \pmg_{xy} \ge 0$\,
\,and\, \,$\pmg_{xz} \ge \pmg_{yz}$\,\,
for any $z\not\in\{x,y\}$.}
\begin{equation}
\box0
\label{eq:mrobinson}
\end{equation}

On the other hand, we also want the rates to be independent of~$\xi$ when there are several possibilities for it.
To this effect, we will require the projected margins to have already such an independence.

The next operation will transform the indirect margins so as to
satisfy these conditions. It is defined in the following way, where we assume $x\rxi y$
and $x'$~denotes the item that immediately follows~$x$ in the total order~$\xi$:
\begin{gather}
\img_{xy} \,=\, \isc_{xy} - \isc_{yx},
\label{eq:cprojection1}
\\[2.5pt]
\ppmg_{xy} \,=\, \min\,\{\, \img_{pq} \;\vert\; p\rxieq x,\; y\rxieq q\,\},
\label{eq:cprojection2}
\\[2.5pt]
\pmg_{xy} \,=\, \max\,\{\, \ppmg_{pp'} \;\vert\; x\rxieq p\rxi y\,\},
\label{eq:cprojection3}
\\[2.5pt]
\pmg_{yx} \,=\, -\pmg_{xy}.
\label{eq:cprojection4}
\end{gather}
\noindent
One can easily see that the $\ppmg_{xy}$ obtained in (\ref{eq:cprojection2}) already satisfy a condition analogous to~(\ref{eq:mrobinson}). However, the independence of~$\xi$ is not be ensured until steps~(\ref{eq:cprojection3}--\ref{eq:cprojection4}).
This property is a consequence of the fact that the projected margins given by the preceding formulas satisfy not only condition~(\ref{eq:mrobinson}) but also the following one:
\begin{equation}
\pmg_{xy} = 0 \,\,\Longrightarrow\,\,
\pmg_{xz} = \pmg_{yz}\,\,\  
\hbox{for any $z\not\in\{x,y\}$}.
\label{eq:msemidefinite}
\end{equation}
In particular, this will happen whenever $\img_{xy}=0$
(since this implies $\ppmg_{pp'}=0$ for all $p$ such that $x\rxieq p\rxi y$). More particularly, in the event of having two admissible orders that interchange two consecutive elements $p$ and $p'$ we will have $\pmg_{pp'} = \ppmg_{pp'} = \img_{pp'} = 0$ and consequently $\pmg_{pz} = \pmg_{p'z}$ for any $z\not\in\{p,p'\}$,
as it is required by the desired independence of~$\xi$.

\medskip
Anyway, the projected margins are
finally introduced in~(\ref{eq:rratesfrommargins}), which determines the rank-like rates~$r_x$.
The corresponding fraction-like rates will be introduced in~\secpar{2.9}.

\remarks
1. Condition~(\ref{eq:mrobinson}) gives the pattern of growth of the projected margins~$\pmg_{pq}$ when $p$ and $q$ vary according to an admissible order~$\xi$. This pattern is illustrated in figure~1 below, where the square represents the matrix~$(\pmg_{pq})$ with $p$ and $q$ ordered according to~$\xi$, from better to worse. As~usual, the first index labels the rows, and the second one labels the columns. The diagonal corresponds to the case $p=q$, which we systematically leave out of consideration. Having said that, here it would be appropriate to put $\pmg_{pp}=0$.
\ensep
Anyway, the projected margins are greater than or equal to zero above the diagonal and smaller than or equal to zero below it, and they increase or remain the same as one moves along the indicated arrows.
The right-hand side of the figure follows from the left-hand one because the projected margins are antisymmetric.

\vbox to2.2in{%
\vskip-.32in
\hskip-.35in\hbox to2.4in{%
\resizebox{\paperwidth}{!}{\includegraphics{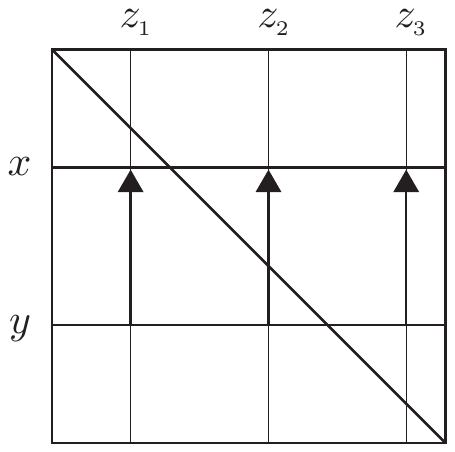}}
\hss}%
\hskip.5in
\hskip-.4in\hbox to2.4in{%
\resizebox{\paperwidth}{!}{\includegraphics{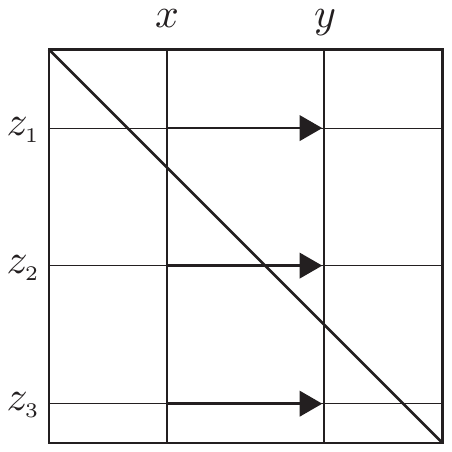}}
\hss}%
\vss
}
\vskip-.05in
\centerline{\label{figNE}Figure~1. Directions of growth of the projected margins.}
\vskip.15in

Of course, the absolute values $d_{xy} = |\pmg_{xy}|$ keep this pattern in the upper triangle but they behave in the reverse way in the lower one. Such a behaviour is often considered in combinatorial data analysis, where it is associated with the name of W.~S.~Robinson, a statistician who in 1951 introduced a condition of this kind as the cornerstone of a method for seriating archaeological deposits
(\ie placing them in chronological order)
\hbox{\brwrap{
\textsl{31}\,\textup{;} 
\textsl{14}\,\textup{:\,\relscale{0.95}\secpar{4.1.1,\,4.1.2,\,4.1.4}\,;} 
\textsl{33}
}}.

Condition~(\ref{eq:msemidefinite}), more precisely its expression in terms of the $d_{xy}$, is also considered in cluster analysis, where it is referred to by saying that the `dissimilarities' 
$d_{xy}$ are `even' \cite[\secpar{9.1}]{js} (`semidefinite' according to other authors).

\medskip
In the present case of complete votes, 
the projected margins $\pmg_{xy}$ defined by~(\ref{eq:cprojection1}--\ref{eq:cprojection4})
satisfy not only~(\ref{eq:mrobinson})
and (\ref{eq:msemidefinite}), but also the stronger condition
\begin{equation}
\pmg_{xz} \,=\, \max\,(\pmg_{xy},\pmg_{yz}),\qquad \hbox{whenever $x\rxi y\rxi z$}.
\label{eq:equaltomax}
\end{equation}
Besides~(\ref{eq:mrobinson}) and (\ref{eq:msemidefinite}), this property implies also that the~$d_{xy}$ satisfy the following inequality,
which makes no reference to the relation $\xi$:
\begin{equation}
d_{xz} \,\le\, \max\,(d_{xy},d_{yz}),\qquad \hbox{for any $x,y,z$.}
\label{eq:ultrametric}
\end{equation}
This condition, called the ultrametric inequality,
is also well known in cluster analysis,
where it appears as a necessary and sufficient condition
for the dissimilarities $d_{xy}$ to define a
hyerarchical classification of the set under consideration
\hbox{\brwrap{
\textsl{14}\,\textup{:\,\relscale{0.95}\secpar{3.2.1}\,;} 
\textsl{29}
}}.

\medskip
Our problem differs from the standard one of combinatorial data analysis
in that our dissimilarities, namely the margins, are antisymmetric,
whereas the standard problem
considers symmetric dissimilarities. In other words, our dissimilarities
have both magnitude and direction,
whereas the standard ones have magnitude only.
This makes an important difference in connection with the seriation problem,
\ie positioning the items on a line.
Let us remark that the case of directed dissimilarities is considered in \cite[\secpar{4.1.2}]{ham}.

\smallskip
2. The operation $(\img_{xy})\rightarrow(\pmg_{xy})$ defined by (\ref{eq:cprojection2}--\ref{eq:cprojection3}) is akin to the single-link method of cluster analysis, which can be viewed as a continuous method for projecting a matrix of dissimilarities onto the set of ultrametric distances; such a continuous projection is achieved by taking the maximal ultrametric distance which is bounded by the given matrix of dissimilarities \cite[\secpar{7.3,\,7.4,\,8.3,\,9.3}]{js}.
The operation $(\img_{xy})\rightarrow(\pmg_{xy})$ does the same kind of job under the constraint that the clusters ---in the sense of cluster analysis--- be intervals of the total order~$\xi$.

\paragraph{2.7}
In order to get more insight into the case of incomplete votes,
it is interesting to look at the case of plumping votes,
\ie the case where each vote plumps for a single option.
In this case, and assuming interpretation~(d),
the binary scores of the vote have the form $v_{xy} = \plumpf_x$ for~every $y\neq  x$,
where $\plumpf_x$ is the fraction of voters who choose~$x$.

In the spirit of condition~\llpv, in this case we expect the projected scores $\psc_{xy}$ to coincide with the original ones~$v_{xy} = \plumpf_x$. So, both the projected margins~$\pmg_{xy}$ and the projected turnovers
$\pto_{xy}$ should also coincide with the original ones, namely~$\plumpf_x-\plumpf_y$ and~$\plumpf_x+\plumpf_y$.
\ensep
In~this connection, one easily sees that the indirect scores $\isc_{xy}$ coincide with $v_{xy}$ (see Proposition~\ref{st:plumpps}).
As~a~consequence, $\xi$ is any total order for which the $\plumpf_x$ are non-increasing.

If we apply formulas~(\ref{eq:cprojection1}--\ref{eq:cprojection4}), we first get 
$\img_{xy} = \isc_{xy} - \isc_{yx} = v_{xy} - v_{yx} = \plumpf_x - \plumpf_y$, and then~$\ppmg_{xx'} = \plumpf_x - \plumpf_{x'}$, but the projected margins resulting from~(\ref{eq:cprojection3}) cease to coincide with the original ones.
\ensep
Most interestingly, such a coincidence would hold
if the \,$\max$\, operator of formula~(\ref{eq:cprojection3}) was replaced by a sum.

Now, these two apparently different operations ---maximum and addition--- can be viewed as particular cases of a general procedure which involves taking the union of certain intervals, namely $\gamma_{xx'} = [\,(t_{xx'}-\ppmg_{xx'})/2\,,\, (t_{xx'}+\ppmg_{xx'})/2\,]$.
\ensep
In~fact, in~the case of complete votes, all the turnovers are equal to~$1$, so these intervals are all of them centred at~$1/2$ and the union operation is equivalent to looking for the maximum of the widths.
\ensep
In~the case of plumping votes, we~know that $t_{xx'} = \plumpf_x + \plumpf_{x'}$ and we have just seen that $\ppmg_{xx'} = \plumpf_x - \plumpf_{x'}$, which implies that $\gamma_{xx'} = [\plumpf_{x'},\plumpf_x]$; so,~the intervals $\gamma_{xx'}$ and $\gamma_{x'x''}$ are  then adjacent to each other (the right end of the latter coincides with the left end of the former) and their union involves adding up the widths.

This remark strongly suggests that the general method should rely on such intervals.
In the following we will refer to them as {\df score intervals}.
A score interval can be viewed as giving a pair of scores about two options, the two scores being respectively in favour and against a specified preference relation about the two options. Alternatively, it can be viewed as giving a certain margin
together with a certain turnover.

More specifically, one is immediately tempted to replace the minimum and maximum operations of~(\ref{eq:cprojection2}--\ref{eq:cprojection3}) by the intersection and union of score intervals.
The starting point would be the score intervals that combine the original turnovers~$t_{xy}$
with the indirect margins~$\img_{xy}$.
Such a procedure works as desired both in the case of complete votes and that of plumping ones. Unfortunately, however, it breaks down in other cases of incomplete votes which produce empty intersections or disjoint unions.
So,~a~more elaborate method is required.

\paragraph{2.8}
In this subsection we will finally describe a rank-like rating procedure which is able to cope with the general case.
This procedure will use score intervals.
However, these intervals will not be based directly on the original turnovers, but on certain transformed ones.
This prior transformation of the turnovers will have the virtue of avoiding the problems pointed out at the end of the preceding paragraph.

So, we are given as input from one side the indirect margins~$\img_{xy}$, and from the other side the original turnovers~$t_{xy}$. The output to be produced is a set of projected scores~$\psc_{xy}$. They should have the virtue that the associated rank-like rating given by~(\ref{eq:rrates}) has the following properties:\ensep
(a)~it is the exactly the same for all admissible orders~$\xi$;
\ensep and
(b)~it is compatible with any such order~$\xi$,
\ie $xy\in\xi$ implies $r_x\le r_y$.

As we did in the complete case, we will require the projected scores~$\psc_{xy}$ to satisfy the condition of independence with respect to~$\xi$.

On the other hand, in order to ensure the compatibility condition~(b), it~suffices that the projected scores behave in the following way:
\setbox0=\vbox{\hsize105mm\noindent
$x\rxi y \,\,\Longrightarrow\,\, \psc_{xy} \ge \psc_{yx}$\,
\,and\, \,$\psc_{xz} \ge \psc_{yz}$\,\,
for any $z\not\in\{x,y\}$.}
\begin{equation}
\box0
\label{eq:vcompatibility}
\end{equation}
If~we think in terms of the associated margins~$\pmg_{xy}$ and turnovers~$\pto_{xy}$ \,\hbox{---which} add up to $2\psc_{xy}$---\, it suffices that both of them satisfy conditions analogous to~(\ref{eq:vcompatibility}).
More, specifically, it suffices that the projected margins be antisymmetric and satisfy condition~(\ref{eq:mrobinson}) of \secpar{2.6} \,and\, that the projected turnovers be symmetric and satisfy
\setbox0=\vbox{\hsize75mm\noindent
$x\rxi y \,\,\Longrightarrow\,\, \pto_{xz} \ge \pto_{yz}$\,\,
for any $z\not\in\{x,y\}$.}
\begin{equation}
\box0
\label{eq:tgreenberg}
\end{equation}

\medskip
So, we want the projected scores to be independent of~$\xi$, and their associated margins and turnovers to satisfy conditions (\ref{eq:mrobinson}) and (\ref{eq:tgreenberg}).
These requirements are fulfilled by the procedure formulated in (\ref{eq:projection1}--\ref{eq:projection7}) below. These formulas use the following notations:
\ensep
$\Psi$~is an operator to be described in a while;
\ensep
$[a,b]$ means the closed interval $\{\,x\in\bbr\mid a\le x\le b\,\}$;
\ensep
\,$|\gamma|$\,~means the length of such an interval \,$\gamma = [a,b]$,\, \ie the number $b-a$;
\ensep
and \,$\centre\gamma$\,~means its barycentre, or centroid,\,
\ie the~number \,$(a+b)/2$.
\ensep
As~in (\ref{eq:cprojection1}--\ref{eq:cprojection4}), the following formulas assume that $x\rxi y$,\, and $x'$~denotes the option that immediately follows $x$ in the total order~$\xi$.
\newcommand\sep{\hskip2em}
\begin{gather}
\img_{xy} \,=\, \hbox to44mm{$\isc_{xy} - \isc_{yx},$\hfil}\sep
t_{xy} \,=\, \hbox to40mm{$v_{xy} + v_{yx},$\hfil}
\label{eq:projection1}
\\[2.5pt]
\ppmg_{xy} \,=\, \hbox to44mm{$
\min\,\{\, \img_{pq} \;\vert\; p\rxieq x,\; y\rxieq q\,\},$\hfil}\sep
\ppto_{xy} =\, \hbox to40mm{$
\Psi\kern.5pt[\kern.5pt(t_{pq}),(\ppmg_{pp'})\kern.5pt]\kern1.5pt_{xy},$\hfil}
\label{eq:projection2}
\\[2.5pt]
\gamma_{xx'} \,=\, \hbox to84mm{$[\,(\ppto_{xx'}-\ppmg_{xx'})/2\,,\, (\ppto_{xx'}+\ppmg_{xx'})/2\,],$\hfil}
\label{eq:projection3}
\\[2.5pt]
\,\gamma_{xy}\, \,=\, \hbox to84mm{$\bigcup\,\{\, \gamma_{pp'} \;\vert\; x\rxieq p\rxi y\,\},$\hfil}
\label{eq:projection4}
\\[2.5pt]
\pmg_{xy} \,=\, \hbox to44mm{$
|\gamma_{xy}|,$\hfil}\sep
\pto_{xy} \,=\, \hbox to40mm{$
2\, \centre\gamma_{xy},$\hfil}
\label{eq:projection5}
\\[2.5pt]
\pmg_{yx} \,=\, \hbox to44mm{$
-\pmg_{xy},$\hfil}\sep
\pto_{yx} \,=\, \hbox to40mm{$
\pto_{xy},$\hfil}
\label{eq:projection6}
\\[2.5pt]
\null\;\psc_{xy}\; \,=\, \hbox to44mm{$\max\gamma_{xy}
 = (\pto_{xy}+\pmg_{xy})/2,
$\hfil}\sep\ 
\null\hskip-5pt 
\psc_{yx} \,=\, \hbox to40mm{$\min\gamma_{xy} 
 = (\pto_{xy}-\pmg_{xy})/2.
$\hfil}
\label{eq:projection7}
\end{gather}

\medskip
Like~(\ref{eq:cprojection1}--\ref{eq:cprojection4}), the preceding  procedure can be viewed as a two-step transformation. The first step is given by~(\ref{eq:projection2}) and it transforms the input margins and turnovers $(\img_{xy}, t_{xy})$ into certain intermediate projections $(\ppmg_{xy}, \ppto_{xy})$ which already satisfy conditions analogous to (\ref{eq:mrobinson}) and (\ref{eq:tgreenberg}) but are not independent of~$\xi$.
The condition of independence requires a second step which is described by (\ref{eq:projection3}--\ref{eq:projection7}). As in~\secpar{2.6}, the superdiagonal final projections coincide with the intermediate ones, \ie $\pmg_{xx'}=\ppmg_{xx'}$ and $\pto_{xx'}=\ppto_{xx'}$. Notice also that the 
intermediate margins $\ppmg_{xy}$ are constructed exactly as in (\ref{eq:cprojection2}).

\newcommand\incrementone{eq:pincrement}
\newcommand\incrementtwo{eq:qincrement}

\medskip
The main difficulty lies in constructing the intermediate turnovers
$\ppto_{xy}$ so that they do not depend on~$\xi$. The reason is that this condition involves the admissible orders, which depend on the relation $\nu$ associated with the indirect margins $\img_{xy}$. So, that construction must take into account not only the original turnovers but also the indirect margins. This connection with the $\img_{xy}$ will be controlled indirectly through the~$\ppmg_{xx'}$. In~fact, we~will look for the $\ppto_{xy}$ so as to satisfy the following conditions:
\begin{gather}
\ppmg_{xx'} \,\le\, \ppto_{xx'}\,\le\, 1,
\label{eq:bounded}
\\[2.5pt]
0 \,\le\, \ppto_{py} - \ppto_{p'y} \,\le\, \ppmg_{pp'},
\label{eq:pincrement}
\\[2.5pt]
0 \,\le\, \ppto_{xq} - \ppto_{xq'} \,\le\, \ppmg_{qq'}\strut.
\label{eq:qincrement}
\end{gather}
Notice that (\ref{eq:pincrement}) ensures that $\ppto_{py}$ and $\ppto_{p'y}$ will coincide with each other whenever $\ppmg_{pp'}=0$. Since $\ppmg_{pp'}=\img_{pp'}$, we are in the case of having two admissible orders that interchange $p$ with $p'$. The fact that this implies $\ppto_{py}=\ppto_{p'\kern-.75pt y}$ eventually ensures the independence of $\xi$ (Theorem~\ref{st:independenceOfXi}; we say ‘eventually’ because the full proof is quite long).

In the case of complete votes we will have $\ppto_{xy}=1$, so that condition~(\ref{eq:bounded}) will be satisfied with an equality sign in the right-hand inequality, whereas (\ref{\incrementone}) and (\ref{\incrementtwo})~will be satisfied with an equality sign in the left-hand inequality.
\ensep
In the case of plumping votes, where we know that 
$\ppmg_{xx'} = m_{xx'} = \plumpf_x - \plumpf_{x'}$ (\secpar{2.7}),
we will have 
$\ppto_{xy} = t_{xy} = \plumpf_x + \plumpf_y$, so that 
(\ref{\incrementone}) and (\ref{\incrementtwo}) will be satisfied with an equality sign in the right-hand inequalities (equation (\ref{eq:bounded}) is satisfied too, but in this case the inequalities can be strict).

Notice also that conditions (\ref{\incrementone}--\ref{\incrementtwo}) imply the following one:
\begin{equation}
0 \,\le\, \ppto_{xx'} - \ppto_{x'x''} \,\le\, \ppmg_{xx'} + \ppmg_{x'x''}.
\label{eq:overlap}
\end{equation}

In geometrical terms, the inequalities in~(\ref{eq:bounded}) mean that\ensep (a)~the interval~$\gamma_{xx'}$ is contained in $[0,1]$. On the other hand, the inequalities in~(\ref{eq:overlap}) mean that the intervals $\gamma_{xx'}$ and $\gamma_{x'x''}$ are related to each other in the following way:\ensep (b)~the barycentre of the first one lies to the right of that of the second one;\ensep (c)~the two intervals overlap each other.

Conditions (\ref{eq:bounded}--\ref{\incrementtwo}) can be easily achieved by taking simply 
$\ppto_{xy}=1$. However, this choice goes against our aim of distinguishing between definite indifference and lack of information; in~particular, condition~\llpv\ requires that in the case of plumping votes the projected turnovers should coincide with the original ones (which are then less than $1$).
\ensep
Now, conditions~(\ref{eq:bounded}--\ref{\incrementtwo}) are convex with respect to the~$\ppto_{xy}$ (the~whole set of them), \ie if~they are satisfied by two different choices of these numbers, they are satisfied also by any convex combination of them.
This~implies that for any given set of original turnovers $t_{xy}$ there is a unique set of values $\ppto_{xy}$ which  minimizes the euclidean distance to the given one while satisfying those conditions.

So,~the operator $\Psi$ can defined in the following way:
\label{defPsi}
$\ppto_{xy}$ is the set of turnovers which is determined by conditions (\ref{eq:bounded}--\ref{\incrementtwo}) together with that of minimizing the following measure of deviation with respect to the~$t_{xy}$:
\begin{equation}
\Phi \,=\, \sum_x \, \sum_y \, (\ppto_{xy} - t_{xy})^2.
\label{eq:phi}
\end{equation}
The actual computation of the~$\ppto_{xy}$ can be carried out in a finite number of steps by means of a quadratic programming algorithm \cite[\secpar{14.1} (2nd ed.)]{lu}.

Anyway, the preceding operations have the virtue of ensuring the desired properties.


\renewcommand\uplapar{\vskip-6.5mm\null} 
\uplapar
\remarks
1. Condition~(\ref{eq:tgreenberg}) is illustrated in figure~2, where the arrows indicate the directions of growth of the projected turnovers. The right-hand side of the figure follows from the left-hand one because the projected turnovers are symmetric.

\vbox to2.2in{%
\vskip-.32in
\hskip-.35in\hbox to2.4in{%
\resizebox{\paperwidth}{!}{\includegraphics{cratingfigs/N}}
\hss}%
\hskip.5in
\hskip-.4in\hbox to2.4in{%
\resizebox{\paperwidth}{!}{\includegraphics{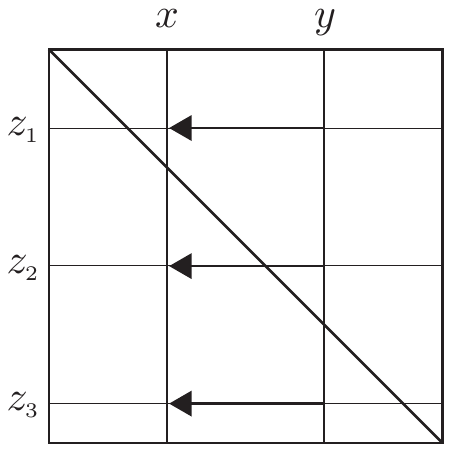}}
\hss}%
\vss
}
\vskip-.05in
\centerline{\label{figNW}Figure~2. Directions of growth of the projected turnovers.}
\vskip.15in

This condition can be associated with the name of Marshall~G.~Greenberg, a mathematical psychologist who in 1965 considered a condition of this form ---at the suggestion of Clyde H.\ Coombs--- in connection with the problem of producing a rating after paired-comparison data, specially in the incomplete case
\hbox{\brwrap{
\textsl{12}\,\textup{;} 
\textsl{14}\,\textup{:\,\relscale{0.95}\secpar{4.1.2}}
}}.
Having said that, we strongly differ
from that author in that he applies a property like~(\ref{eq:tgreenberg}) to the scores, whereas we consider more
appropriate to apply it to the turnovers.

In fact, under the general assumption that each vote is a ranking,
possibly incomplete,
and that each ranking is translated into a set of binary preferences according to rules (a--d) of~\secpar{2.1},
it is fairly reasonable to expect that the turnover for a pair~$xy$, \ie the number of voters who expressed an opinion about $x$ in comparison with~$y$, should increase as $x$ and/or $y$ are higher in the social ranking. In~practice, the original turnovers can deviate to a certain extent from this ideal behaviour. In~contrast, our projected turnovers are always in agreement with it (with respect to the total order~$\xi$).

\smallskip
2. The projected scores turn out to satisfy not only~(\ref{eq:vcompatibility}), but also the following stronger property: 
\setbox0=\vbox{\hsize105mm\noindent
if\, \,$x\rxi y$\, \,then\, \,$\psc_{xy} \ge \psc_{yx}$\hfil\break
\hphantom{if\, \,$x\rxi y$\, \,\,}%
and\,\, \,$\psc_{xz} \ge \psc_{yz},\ \psc_{zx} \le \psc_{zy}$\,\,
for any $z\not\in\{x,y\}$.}
\begin{equation}
\box0
\label{eq:vrobinson}
\end{equation}
So, the projected scores increase or remain constant in the directions shown in figure~1.\ensep Furthermore, we will see that the quotients $\pmg_{xy}/\pto_{xy}$ have also the same property.

\smallskip
3. In contrast to the case of complete votes, in this case the projected margins do not satisfy~(\ref{eq:equaltomax}) but only
\begin{equation}
\pmg_{xz} \,\le\, \pmg_{xy} + \pmg_{yz},\qquad \hbox{whenever $x\rxi y\rxi z$}.
\label{eq:mtriangular}
\end{equation}
As a consequence, the absolute values $d_{xy} = |\pmg_{xy}|$
satisfy the triangular inequality:
\begin{equation}
d_{xz} \,\le\, d_{xy} + d_{yz},\qquad \hbox{for any $x,y,z$.}
\label{eq:triangular}
\end{equation}

\smallskip
4. Under the assumption of ranking votes (but not necessarily in a more general setting)
one can see that the original turnovers already  satisfy (\ref{eq:bounded}).
In this case, the preceding definition of $\Psi$ turns out to be  equivalent to an analogous one where conditions~(\ref{eq:bounded}--\ref{\incrementtwo}) are replaced simply by~(\ref{\incrementone}--\ref{\incrementtwo}). From this it follows that the intermediate turnovers have then the same sum as the original ones:
\begin{equation}
\sum_x \, \sum_y \, \ppto_{xy} \,=\, \sum_x \, \sum_y \, t_{xy}.
\label{eq:mateixasuma}
\end{equation}

\uplapar
\paragraph{2.9}
Finally, let us see how shall we define the fraction-like rates~$\flr_x$. As~in the case of the rank-like rates, we will use a classical method which would usually be applied to the original scores, but here we will apply it to the projected scores.
This method was introduced in~1929 by Ernst Zermelo~\cite{ze} 
and it was rediscovered by other authors in the 1950s~\cite{bt,fo}
(see also~\cite{ke}).
Zermelo's work was motivated by chess tournaments, whereas the other   authors were considering comparative judgments.
Anyway, all of them were especially interested in the incomplete case, \ie the case where turnovers may depend on the pair~$xy$. 

\medskip
More specifically,
the fraction-like rates~$\flr_x$ will be determined by the following system of equations (together with the condition that $\flr_x\ge0$ for every~$x$):
\begin{align}
\sum_{y\neq x}\,\pto_{xy}\,\flr_x/(\flr_x+\flr_y) \,&=\, \sum_{y\neq x}\,\psc_{xy}\ \ (= N-r_x),
\label{eq:frates}
\\[2.5pt]
\sum_x\,\flr_x \,&=\, \nev,
\label{eq:fratesa}
\end{align}
where (\ref{eq:frates})~contains
one equation for every $x$, and $\nev$~stands for the fraction of non-empty votes (\ie $\nev = F/V$ where $F$~is the number of non-empty votes and $V$ is the total number of votes).
In spite of having $N+1$ equations, the $N$~equations contained in (\ref{eq:frates}) are not independent, since their sum results in the identity
$(\sum_x\sum_{y\neq x}\pto_{xy})/2 = \sum_x\sum_{y\neq x}\psc_{xy}$. On the other hand, it is clear that (\ref{eq:frates})~is insensitive to all of the~$\flr_x$ being multiplied by a~constant factor. This indeterminacy disappears once (\ref{eq:frates})~is supplemented with  equation~(\ref{eq:fratesa}).

\medskip
In the case of plumping votes, where we know that $\psc_{xy} = \plumpf_x$ and $\pto_{xy} = \plumpf_x + \plumpf_y$, the solution of (\ref{eq:frates}--\ref{eq:fratesa}) is easily seen to be $\flr_x = \plumpf_x$, as required by condition~\llpv.


\medskip
The problem of solving the system~(\ref{eq:frates}--\ref{eq:fratesa}) is well posed when the projected Llull matrix~$(\psc_{xy})$ is irreducible. This means that there is no splitting of the options into a~`top class'~$X$ plus a `low~class'~$Y$ so that $\psc_{yx}=0$ for any $x\in\xst$ and $y\in\yst$. When such a splitting exists, one is forced to put $\flr_y=0$ for all~$y\in\yst$.
For more details, the reader is referred to section~11.



\medskip
Zermelo (and the other authors) dealed also with the problem of numerically solving a non-linear system of the form~(\ref{eq:frates}). In this connection, he~showed that in the irreducible case its solution (up to a multiplicative constant) can be approximated to an arbitrary degree of accuracy by means of an iterative scheme of the form
\begin{equation}
\flr^{n+1}_x \,=\,
\bigg(\sum_{y\neq x}\,\psc_{xy}\bigg) \bigg/ \bigg(\sum_{y\neq x}\,\pto_{xy}/(\flr^{n}_x+\flr^{n}_y)\bigg),
\label{eq:iterates}
\end{equation}
starting from an arbitrary set of values~$\flr^0_x > 0$.

\medskip
The fraction-like rates~$\flr_x$ determined by~(\ref{eq:frates}--\ref{eq:fratesa}) can be viewed as an estimate of the first-choice fractions using not only the first choices but the whole rankings. 
Properly speaking, Zermelo's method (with the original scores and turnovers) corresponds to a maximum likelihood estimate of the parameters of a certain probabilistic model for the outcomes of a tournament between several players, or, more in the lines of our applications, for the outcomes of comparative judgments.
This model will be briefly described in section~11.
Although we are far from its hypotheses, 
we will see that Zermelo's method is quite suitable for translating our rank-like rates into fraction-like ones.

\renewcommand\uplapar{\vskip-6mm\null} 

\uplapar
\section{Summary of the method. Variants. General forms of vote}

\uplapar
\paragraph{3.1}Let us summarize the whole procedure. In the \emph{general case}, where the votes are not necessarily complete, it consists of the following steps:

\begin{enumerate}
\setlength\itemsep{0pt}
\item Form the Llull matrix~$(v_{xy})$~(\secpar{2.1}). Work out the turnovers~$t_{xy}=v_{xy}+v_{yx}$.
\item Compute the indirect scores~$\isc_{xy}$ defined by~(\ref{eq:paths}). An efficient way to do it is the Floyd-Warshall algorithm \cite[\secpar{25.2}]{co}. Work out the indirect margins $\img_{xy} \cd= \isc_{xy}\cd-\isc_{yx}$ and the associated indirect comparison relation $\nu = \{xy\mid\img_{xy} > 0\}$.
\item Find an admissible order~$\xi$~(\secpar{2.4}) and arrange the options according to it. For instance, it suffices to arrange the options by non-decreasing values of the ‘tie-splitting’
Copeland scores $\rank x = 1 + |\{\,y\mid y\cd\neq x,\ \img_{yx}\cd>0\}|\linebreak + \frac12\,|\{\,y\mid y\cd\neq x,\ \img_{yx}\cd=0\}|$ (Proposition~\ref{st:Copeland}).
\item Starting from the indirect margins~$\img_{xy}$, work~out the superdiagonal intermediate projected margins $\ppmg_{xx'}$ as defined in (\ref{eq:projection2}.1).
\item Starting from the original turnovers~$t_{xy}$, and taking into\linebreak account the superdiagonal intermediate projected margins $\ppmg_{xx'}$,\linebreak determine the intermediate projected turnovers $\ppto_{xy}$ so as to minimize (\ref{eq:phi}) under the constraints (\ref{eq:bounded}--\ref{\incrementtwo}). This can be carried out in a finite number of steps by means of a quadratic programming algorithm \cite[\secpar{14.1} (2nd ed.)]{lu}.
\item Form the intervals $\gamma_{xx'}$ defined by (\ref{eq:projection3}), derive their unions $\gamma_{xy}$ as defi\-ned by~(\ref{eq:projection4}),\, and read off the projected scores~$\psc_{xy}$ (\ref{eq:projection7}).
\item Compute the rank-like ranks~$r_x$ according to (\ref{eq:rrates}).
\item Determine the fraction-like rates $\flr_x$ by solving the system (\ref{eq:frates}--\ref{eq:fratesa}). This can be done numerically by means of the iterative scheme~(\ref{eq:iterates}).
\end{enumerate}

\smallskip
In the \emph{complete case}, the scores~$v_{xy}$ and the margins~$m_{xy}$ are~related to each other by the monotone increasing transformation $v_{xy} = (1+m_{xy})/2$.\linebreak[3] Because of this fact, the preceding procedure can then be simplified in the following way:

\begin{itemize}
\setlength\itemsep{0pt}
\item Step~2 computes $m^*_{xy}$ instead of $v^*_{xy}$ and takes
$\img_{xy} = (m^*_{xy}- m^*_{yx})/2$.
\item Step~5 is not needed.
\item Step~6 reduces to (\ref{eq:cprojection3}--\ref{eq:cprojection4}).
\item Step~7 makes use of formula~(\ref{eq:rratesfrommargins}).
\end{itemize}

\paragraph{3.2}
The preceding procedure admits of certain variants which might be appropriate to some special situations. Next we will distinguish four of them, namely
\vskip-6pt

\newcommand\iitm[2]{%
\halfsmallskip\hskip2\parindent\hbox to1.5em{#1.\hss}#2}

\iitm{1}{Main}
\iitm{2}{Dual}
\iitm{3}{Balanced}
\iitm{4}{Margin-based}
\smallskip

\noindent
The above-described procedure is included in this list as the main variant. The four variants are exactly equivalent to each other in the complete case, but in the incomplete case they can produce different results. In spite of this, they all share the main properties.

\smallskip
The dual variant is analogous to the main one except that the \hbox{max-min} indirect scores $\isc_{xy}$ are replaced by the following \hbox{min-max} ones:
\begin{equation}
\stv_{xy} \hskip.75em = \hskip.75em
\min_{\vtop{\scriptsize\halign{\hfil#\hfil\cr\noalign{\vskip.5pt}$x_0=x$\cr$x_n=y$\cr}}}
\hskip.75em
\max_{\vtop{\scriptsize\halign{\hfil#\hfil\cr\noalign{\vskip-1.25pt}$i\ge0$\cr$i<n$\cr}}}
\hskip.75em v_{x_ix_{i+1}}.
\label{eq:dual}
\end{equation}
Equivalently, $\stv_{xy} = 1 - {\hat v}^\ast_{yx}$ where ${\hat v}_{xy} = 1 - v_{yx}$. In the complete case one has $\stv_{xy} = 1 - \isc_{yx}$, so that $\stv_{xy} - \stv_{yx} = \isc_{xy} - \isc_{yx}$ and $\mu(\stv) = \mu(\isc)$; as a consequence, the dual variant is then equivalent to the main one.

\smallskip
The balanced variant takes $\nu = \mu(\isc) \cap \mu(\stv)$ together with
\begin{equation}
\img_{xy} = \begin{cases}
\min\,(\isc_{xy} - \isc_{yx}, \stv_{xy} - \stv_{yx}),
& \text{if $xy\in\mu(\isc) \cap \mu(\stv)$,}\\
-\img_{yx},
& \text{if $yx\in\mu(\isc) \cap \mu(\stv)$,}\\
0,
& \text{otherwise.}\\
\end{cases}
\label{eq:balanced}
\end{equation}
The remarks made in connection with the dual variant show that in the complete case the balanced variant is also equivalent to the preceding ones.

\pagebreak\null\vskip-14mm\null 
The margin-based variant follows the simplified procedure of the end of~\secpar{3.1} even if one is not originally in the complete case. 
Equivalently, it corresponds to replacing the original scores $v_{xy}$ by the following ones: $v'_{xy} = (1+m_{xy})/2$.
This amounts to replacing any lack of information about a pair of options by a definite indifference between them, which brings the problem into the complete case. 
So, the specific character of this variant lies only in its  interpretation of incomplete votes. Although this interpretation goes against the general principle stated in \secpar{1.5}, it may be suitable to certain situations where the voters are well acquainted with all of the options.
\ensep
In the case of ranking votes, it amounts to replace rule~(d) of \secpar{2.1} by the following one:

\iim{d$'$} When neither $x$ nor $y$ are in the list,\,
we interpret that they are considered equally good
(or equally bad),\, so we proceed as in~(b).

\noindent
In other words, each truncated vote is completed by appending to it all the missing options tied to each other.

\vskip-6mm\null 
\remark
Other variants ---in the incomplete case--- arise when equation (\ref{eq:rrates}) is replaced by the following one:
\begin{equation}
r_x \,=\, 1 + \sum_{y\neq x}\,\psc_{yx}.
\label{eq:rratesVar}
\end{equation}

\vskip-6mm\null 
\paragraph{3.3}
Most of our results will hold if the ``votes'' are not required to be rankings, but they are allowed to be general binary relations.
In particular, this allows to deal with certain situations where it makes sense to replace rule~(c) of \secpar{2.1} by the following one:

\iim{c$'$} When $x$ is in the list and $y$ is not in it,\,
we interpret nothing about the preference of the voter between $x$ and $y$.

\noindent
One could even allow the votes to be non-transitive binary relations;
such a lack of transitivity in the individual preferences may 
arise when individuals are aggregating a variety of criteria~\cite{he}.

A vote in the form of a binary relation $\rho$ contributes to the binary scores with the following amounts:
\begin{equation}
\label{eq:binrelmatrix}
v_{xy} =
\begin{cases}
1, &\text{if } xy\in\rho \text{ and } yx\notin\rho\\
1/2, &\text{if }  xy\in\rho \text{ and } yx\in\rho\\
0, &\text{if } xy\notin\rho.
\end{cases}
\end{equation}

\bigskip
Even more generally, a vote could be any set of normalized binary scores, \ie an element of the set~$\Omega = \{\,v\in[0,1]^{\textit{\char"05}} \mid v_{xy}+v_{yx}\le1\,\}$, where $\tie$~denotes the set of pairs $xy\in\ist\times\ist$ with $x\neq  y$.

\bigskip
Anyway, the collective Llull matrix is simply the center of gravity of a~distribution of individual votes:
\begin{equation}
\label{eq:cog}
v_{xy} = \sum_k \alpha_k\, \vk_{xy},
\end{equation}
where $\alpha_k$ are the relative frequencies or weights of the individual votes~$\vk$.

\newbox\Strutbox
\newdimen\Strutheight

\def\initsize{%
 \Strutheight=.7\baselineskip \advance\Strutheight by2pt
 \setbox\Strutbox=\hbox{\vrule height\Strutheight depth .3\baselineskip width0pt}%
 \def\Strut{\relax\unhcopy\Strutbox}%
}
\initsize

\def\initsizesmall{%
 \Strutheight=.7\baselineskip \advance\Strutheight by2pt
 \setbox\Strutbox=\hbox{\vrule height\Strutheight depth .3\baselineskip width0pt}%
 \def\Strut{\relax\unhcopy\Strutbox}%
}

\edef\bv{|} \catcode`\|=\active \let|\bv
\edef\bi{/} \catcode`\/=\active \let/\bi
\def\htskip{\hskip.4em plus.25em minus.25em}
\newcount\nspan

\def\ls#1{\hbox{{\xipt\sl@\/}#1}}
\def\Ls#1{\hbox{{\xipt\sl@\/}#1}}
\def≈{\,=\,}

\newif\ifhborders

\def\taula#1{
\vbox\bgroup
\def|{&&}
\ifhborders
\def/{\unskip&\cr\noalign{\hrule height.2pt}\Strut&&\ignorespaces}
\else
\def/{\unskip&\cr\strut&&\ignorespaces}
\fi
\def\hrf{\unskip&\cr\noalign{\hrule height.2pt}\Strut&&\ignorespaces}
\def\hrg{\unskip&\cr\noalign{\hrule height.8pt}\Strut&&\ignorespaces}
\halign\bgroup##\hfil\cr 
\hbox{\sc\ #1\hss}\cr 
\noalign{\vskip5pt}
\hbox\bgroup\vrule
}

\def\fitaula{\vrule\egroup\cr\egroup\egroup}

\def\bloc #1#2#3#4#5#6#7◊{%
\nspan=#1\multiply\nspan by2\advance\nspan by-1%
\vrule width.2pt%
\hbox{\vbox{\offinterlineskip\halign{%
##&\vrule## width.2pt%
&&\htskip#2{#4##\/}#3\htskip &\vrule## width.2pt\cr
\noalign{\hrule height.8pt}
\Strut&&\multispan\nspan\it#5\/\hfil&\cr
\strut&&#6&\cr
\noalign{\hrule height.8pt}
\Strut&&\ignorespaces#7&\cr
\noalign{\hrule height.8pt}
}}}%
\vrule width.2pt}

\def\multicol(#1)#2:#3◊{\bloc{#1}\hfil\hfil\rm\null{#2}{#3}◊}
\def\hmulticol(#1)#2:#3:#4◊{\bloc{#1}\hfil\hfil\rm{#2}{#3}{#4}◊}
\def\gap#1{\vrule\egroup{#1}\hbox\bgroup\vrule}

\def\numeros#1◊{\bloc1\hfil\hfil\sf\null{\it Id\/}#1◊} 
\def\numx#1◊{\bloc1\hfil\hfil\sf\null{$x$}#1◊} 
\def\jutges(#1)#2:#3◊{\bloc{#1}\hfil\hfil\rm{Judges}{#2}{#3}◊}
\def\calculs(#1)#2:#3◊{\bloc{#1}\hfil\hfil\rm{Computations\ }{#2}{#3}◊}
\def\balls(#1)#2:#3◊{\bloc{#1}\hfil\hfil\rm{Dances}{#2}{#3}◊}
\def\rxxi(#1)#2:#3◊{\bloc{#1}\null\hfil\rm{Computations}{#2}{#3}◊}
\def\balljutges#1(#2)#3:#4◊{\bloc{#2}\hfil\hfil\rm{#1}{#3}{#4}◊}
\def\fballs(#1)#2:#3◊{\bloc{#1}\null\hfil\rm{Dances}{#2}{#3}◊}
\def\rate#1◊{\bloc1\hfil\hfil\bf\null{$X$}{#1}◊}

\def\fl{\it} 
\def\al{\sl} 

\def\llocs#1◊{\bloc1\hfil\hfil\fl\null{\it R\/}#1◊}
\def\suma#1◊{\bloc1\hfil\hfil\bf\null{\it S\/}#1◊}
\def\promig#1◊{\bloc1\hfil\hfil\bf\null{\it A\/}{#1}◊}
\def\reordenats(#1):#2◊{\bloc{#1}\hfil\hfil\rm\null%
{\multispan\nspan\hfil{\it Rearranged\/}\hskip3.5pt\hfil}{#2}◊}
\def\mediana#1◊{\bloc1\hfil\hfil\bf\null{\it M\/}#1◊}
\def\coco#1◊{\bloc1\hfil\hfil\bf\null{\it C\/}#1◊}

\def\dmedas(#1)#2:#3◊{\bloc{#1}\hfil\hfil\rm{Dances}{#2}{#3}◊}
\def\jmedas(#1)#2:#3◊{\bloc{#1}\hfil\hfil\rm{Judges}{#2}{#3}◊}
\def\bmedas(#1)#2:#3◊{\bloc{#1}\hfil\hfil\rm{Balance}{#2}{#3}◊}

\def\mballs(#1)#2:#3◊{\bloc{#1}\hfil\hfil\rm{Dances \rm($M$)}{#2}{#3}◊}
\def\md#1◊{\bloc1\hfil\hfil\bf\null{$M^D$}#1◊}
\def\mj#1◊{\bloc1\hfil\hfil\bf\null{$M^J$}#1◊}
\def\mb#1◊{\bloc1\hfil\hfil\bf\null{$M^B$}#1◊}
\def\liballs(#1)#2:#3◊{\bloc{#1}\hfil\hfil\rm{Dances \rm($L_1$)}{#2}{#3}◊}
\def\lid#1◊{\bloc1\hfil\hfil\rm\null{$L_1^D$}#1◊}

\def\xd#1◊{\bloc1\null\hfil\bf\null{$X^D$}{#1}◊}
\def\xj#1◊{\bloc1\null\hfil\bf\null{$X^J$}#1◊}
\def\xdj#1◊{\bloc{2}\hfil\hfil\rm\null{$X^D$|$X^J$}{#1}◊}
\def\xb#1◊{\bloc1\null\hfil\bf\null{$X^B$}#1◊}

\def\rates#1◊{\bloc1\hfil\hfil\rm\null{\it X\/}#1◊}
\def\ratesbf#1◊{\bloc1\hfil\hfil\bf\null{\it X\/}#1◊}

\def\peu{\unskip&\cr\noalign{\hrule height.8pt}\Strut&&\ignorespaces}
\def\pballs#1:#2◊{\bloc1\hfil\hfil\rm{\sf #1 }\null
{#2\peu $\scriptstyle D$/$\scriptstyle J$\peu$\scriptstyle B$}◊}

\def\hbi#1{\hbox to.55em{\hss #1\hss}}
\def\hbij#1{\hbox to.85em{\hss #1\hss}}
\def\hbii#1{\hbox to1em{\hss #1\hss}}
\def\hbiii#1{\hbox to1.3em{\hss #1\hss}}

\def\iiijutges{\hbi{\sf A}|\hbi{\sf B}|\hbi{\sf C}}
\def\ivjutges{\hbi{\sf A}|\hbi{\sf B}|\hbi{\sf C}|\hbi{\sf D}}
\def\vjutges{\hbi{\sf A}|\hbi{\sf B}|\hbi{\sf C}|\hbi{\sf D}|\hbi{\sf E}}
\def\viijutges{\hbi{\sf A}|\hbi{\sf B}|\hbi{\sf C}|\hbi{\sf D}|\hbi{\sf E}|\hbi{\sf F}|\hbi{\sf G}}
\def\ixjutges{\hbi{\sf A}|\hbi{\sf B}|\hbi{\sf C}|\hbi{\sf D}|\hbi{\sf E}|\hbi{\sf F}|\hbi{\sf G}|\hbi{\sf H}|\hbi{\sf I}}
\def\xijutges{\hbi{\sf A}|\hbi{\sf B}|\hbi{\sf C}|\hbi{\sf D}|\hbi{\sf E}|\hbi{\sf F}|\hbi{\sf G}|\hbi{\sf H}|\hbi{\sf I}|\hbi{\sf J}|\hbi{\sf K}}
\def\standard{\hbi{W}|\hbi{T}|\hbi{V}|\hbi{F}|\hbi{Q}}
\def\standardUK{\hbi{W}|\hbi{T}|\hbi{F}|\hbi{Q}|\hbi{V}}
\def\llatins{\hbi{C}|\hbi{S}|\hbi{R}|\hbi{P}|\hbi{J}}
\def\llatinsIDSF{\hbi{S}|\hbi{C}|\hbi{R}|\hbi{P}|\hbi{J}}
\def\wtq{\hbi{W}|\hbi{T}|\hbi{Q}}
\def\wtqii{\hbii{W}|\hbii{T}|\hbii{Q}}
\def\iiijutgesii{\hbii{A}|\hbii{B}|\hbii{C}}

\def\rkf#1#2{\hbi{\hbox{\vbox to0pt{\hsize1ex\vss\halign{\hss##\hss\cr
 #2\cr\noalign{\vskip5pt}#1\cr}}}}}

\def\onehalf{\leavevmode\xiipt
 \raise.4ex\hbox{\xpt 1}\hskip-.2ex/\hskip-.15ex\lower.4ex\hbox{\xpt 2}\hskip.2ex}

\def¶#1{{\,\sf#1\,}}   

\def\qbloc #1#2#3#4#5#6◊{%
\nspan=#1\multiply\nspan by2\advance\nspan by-1%
\vrule width.2pt%
\hbox{\vbox{\offinterlineskip\halign{%
##&\vrule## width.2pt%
&&\htskip#2{#4##\/}#3\htskip &\vrule## width.2pt\cr
\noalign{\hrule height.8pt}
\Strut&&#5&\cr
\noalign{\hrule height.8pt}
\Strut&&\ignorespaces#6&\cr
\noalign{\hrule height.8pt}
}}}%
\vrule width.2pt}

\def\qnumx#1◊{\qbloc1\hfil\hfil\sf{$x$}#1◊} 
\def\qmulticol(#1)#2:#3◊{\qbloc{#1}\hfil\hfil\rm{\multispan\nspan\it#2\/\hfil}{#3}◊}

\def\xnumx#1◊{\bloc1\hfil\hfil\sf{\lower1pt\hbox{\Strut}}{$x$}#1◊} 

\hborderstrue

\let\xx\bf
\newcommand\hsk{\hskip5pt}
\newcommand\hssk{\hskip2.5pt}
\newcommand\st{$\ast$}
\newcommand\hph{\hphantom{2}}
\renewcommand\onehalf{\leavevmode\raise.5ex\hbox{\scriptsize$1$}\hskip-.2ex/\hskip-.2ex\lower.2ex\hbox{\scriptsize$2$}}
\newcommand\hq{\hskip.25em}

\section{Examples}

\paragraph{4.1} As a first example of a vote which involved truncated rankings, we look at an election which took place the 16th of February of 1652 in the Spanish royal household. This election is quoted in~\cite{ri}, but we use the slightly different data which are given in~\cite[vol.\,2, p.\,263--264]{cv}. The office under election was that of ``aposentador mayor de palacio'', and the king was assessed
by six noblemen, who expressed the following preferences:

\newcommand\vote[2]{\leftline{\hskip2\parindent
 \hbox to10cm{#1\ \dotfill\ \hbox to5em{#2\hfill}}}}

\medskip
\vote{Marqu\'es de Ari\c{c}a}{¶{b}\better ¶{e}\better ¶{d}\better ¶{a}}
\vote{Conde de Barajas}{¶{b}\better ¶{a}\better ¶{f}}
\vote{Conde de Montalb\'an}{¶{a}\better ¶{f}\better ¶{b}\better ¶{d}}
\vote{Marqu\'es de Povar}{¶{e}\better ¶{b}\better ¶{f}\better ¶{c}}
\vote{Conde de Pu\~nonrostro}{¶{e}\better ¶{a}\better ¶{b}\better ¶{f}}
\vote{Conde de Ysingui\'en}{¶{b}\better ¶{d}\better ¶{a}\better ¶{f}}

\medskip
\noindent
The candidates ¶{a}--¶{f} nominated in these preferences were:

\newcommand\candidate[2]{\leftline{\hskip2.3\parindent
 \hbox to.7em{\hss¶{#1}\hss}\quad #2}}

\medskip
\candidate{a}{Alonso Carbonel \,(architect, 1583--1660)}
\candidate{b}{Gaspar de Fuensalida \,(died 1664)} 
\candidate{c}{Joseph Nieto}
\candidate{d}{Sim\'on Rodr\'{\i}guez}
\candidate{e}{Francisco de Rojas \,(1583--1659)} 
\candidate{f}{Diego Vel\'azquez \,(painter, 1599--1660)}
\medskip

\bigskip
The CLC computations are as follows:

\bgroup\small
\initsizesmall

\bigskip
\leftline{\hskip2\parindent
\hbox to60mm{\ignorespaces
\taula{}%
\numx a/b/c/d/e/f◊%
\hmulticol(6)\hskip3pt$V_{xy}$:
\sf a|\sf b|\sf c|\sf d|\sf e|\sf f:
\st|2|5|3|3|5/
4|\st|6|6|4|5/
1|0|\st|1|0|0/
2|0|3|\st|2|2/
3|2|3|3|\st|3/
1|1|5|4|3|\st◊%
\fitaula
\hfill}%
\taula{}%
\numx a/b/c/d/e/f◊%
\hmulticol(6)\hskip3pt\smash{\hbox{$V^\ast_{xy}$}}:
\sf a|\sf b|\sf c|\sf d|\sf e|\sf f:
\st|2|\bf5|\bf4|\bfgr{3}|\bf5/
\bf4|\st|\bf6|\bf6|\bf4|\bf5/
1|1|\st|1|1|1/
2|2|\bf3|\st|2|2/
\bfgr{3}|2|\bf3|\bf3|\st|\bfgr{3}/
3|2|\bf5|\bf4|\bfgr{3}|\st◊%
\hmulticol(1):$\kappa$: 2\pq12/ 1/ 6/ 5/ 3/ 3\pq12◊
\fitaula
}

\bigskip
\leftline{\hskip2\parindent
\hbox to60mm{\ignorespaces
\taula{}%
\numx b/a/e/f/d/c◊%
\hmulticol(6)\hskip3pt\smash{\hbox{$M^\mast_{xy}$}}:
\sf b|\sf a|\sf e|\sf f|\sf d|\sf c:
\st|2|2|3|4|5/
\st|\st|0|2|2|4/
\st|\st|\st|0|1|2/
\st|\st|\st|\st|2|4/
\st|\st|\st|\st|\st|2/
\st|\st|\st|\st|\st|\st◊%
\fitaula
\hfill}%
\taula{}%
\numx b/a/e/f/d/c◊%
\hmulticol(6)\hskip3pt$T_{xy}$:
\sf b|\sf a|\sf e|\sf f|\sf d|\sf c:
\st|6|6|6|6|6/
\st|\st|6|6|5|6/
\st|\st|\st|6|5|3/
\st|\st|\st|\st|6|5/
\st|\st|\st|\st|\st|4/
\st|\st|\st|\st|\st|\st◊%
\fitaula
}

\bigskip
\leftline{\hskip2\parindent
\hbox to60mm{\ignorespaces
\taula{}%
\numx b/a/e/f/d/c◊%
\hmulticol(6)\hskip3pt\smash{\hbox{$M^\sigma_{xy}$}}:
\sf b|\sf a|\sf e|\sf f|\sf d|\sf c:
\st|2|2|3|4|5/
\st|\st|0|2|2|4/
\st|\st|\st|0|1|2/
\st|\st|\st|\st|1|2/
\st|\st|\st|\st|\st|2/
\st|\st|\st|\st|\st|\st◊%
\fitaula
\hfill}%
\taula{}%
\numx b/a/e/f/d/c◊%
\hmulticol(6)\hskip3pt\smash{\hbox{$T^\sigma_{xy}$}}:
\sf b|\sf a|\sf e|\sf f|\sf d|\sf c:
\st|6|6|6|6|6/
\st|\st|6|6|5\pq13|4\pq23/
\st|\st|\st|6|5\pq13|4\pq23/
\st|\st|\st|\st|5\pq13|4\pq23/
\st|\st|\st|\st|\st|4/
\st|\st|\st|\st|\st|\st◊%
\fitaula
}

\bigskip
\leftline{\hskip2\parindent
\hbox to60mm{\ignorespaces
\taula{}%
\numx b/a/e/f/d/c◊%
\hmulticol(6)\hskip3pt\smash{\hbox{$M^\pi_{xy}$}}:
\sf b|\sf a|\sf e|\sf f|\sf d|\sf c:
\st|2|2|2|2|3/
\st|\st|0|0|1|2\pq16/
\st|\st|\st|0|1|2\pq16/
\st|\st|\st|\st|1|2\pq16/
\st|\st|\st|\st|\st|2/
\st|\st|\st|\st|\st|\st◊%
\fitaula
\hfill}%
\taula{}%
\numx b/a/e/f/d/c◊%
\hmulticol(6)\hskip3pt\smash{\hbox{$T^\pi_{xy}$}}:
\sf b|\sf a|\sf e|\sf f|\sf d|\sf c:
\st|6|6|6|6|5/
\st|\st|6|6|5\pq13|4\pq16/
\st|\st|\st|6|5\pq13|4\pq16/
\st|\st|\st|\st|5\pq13|4\pq16/
\st|\st|\st|\st|\st|4/
\st|\st|\st|\st|\st|\st◊%
\fitaula
}

\bigskip
\leftline{\hskip2\parindent
\hbox to60mm{\ignorespaces
\taula{}%
\numx b/a/e/f/d/c◊%
\hmulticol(6)\hskip3pt\smash{\hbox{$V^\pi_{xy}$}}:
\sf b|\sf a|\sf e|\sf f|\sf d|\sf c:
\st|4|4|4|4|4/
2|\st|3|3|3\pq16|3\pq16/
2|3|\st|3|3\pq16|3\pq16/
2|3|3|\st|3\pq16|3\pq16/
2|2\pq16|2\pq16|2\pq16|\st|3/
1|1|1|1|1|\st◊%
\fitaula
\hfill}%
\taula{}%
\numx b/a/e/f/d/c◊%
\hmulticol(2):$r_x$|$\flr_x$:
2.6667|0.3049/
3.6111|0.1703/
3.6111|0.1703/
3.6111|0.1703/
4.0833|0.1293/
5.1667|0.0549◊%
\fitaula
}

\egroup

\smallskip\null\pagebreak 

\bigskip
According to these results, the office should have been given to candidate~¶{b}, who is also the winner by most other methods.
In the CLC method, this candidate is followed by three runners-up tied to each other, namely candidates ¶{a},¶{e} and ¶{f}. In spite of the clear advantage of candidate~¶{b}, the king appointed candidate ¶{f}, namely, the celebrated painter Diego Vel\'azquez.

\paragraph{4.2}
As an example where the votes are complete strict rankings, we will consider the final round of a dancesport competition. Specifically, we will take the final round of the Professional Latin Rising Star section of the 2007 Blackpool Dance Festival (Blackpool, England, 25th May 2007). The data were taken from  
\texttt{\small http://www.scrutelle.info/results/estelle/2007/blackpool\linebreak \underl2007/}.

As usual, the final was contested by six couples, whose numbers were ¶{3},¶{4},¶{31},¶{122},¶{264},¶{238}. Eleven adjudicators ranked their simultaneous performances in four equivalent dances.

The all-round official result was ¶{3}\better ¶{122}\better ¶{264}\better ¶{4}\better ¶{31}\better ¶{238}.
This\linebreak[3] result comes from the so-called ``Skating System'', whose name reflects a prior use in figure-skating. The Skating System has a first part which produces a separate result for each dance. This~is done mainly on the basis of the median rank obtained by each couple, a~criterion which Condorcet proposed as a ``practical'' method in 1792/93~\cite[ch.\,8]{mu}. However, the fine properties of this criterion are lost in the second part of the Skating System, where the~all-round result is obtained by adding the up the final ranks obtained in the~different dances.

From the point of view of paired comparisons, it makes sense to base the all-round result on the Llull matrix which collects the 44~rankings produced by the 11~adjudicators over the 4~dances~\cite[\secpar{11}]{xm}. As~one can see below, in the present case this matrix exhibits several Condorcet cycles, like for instance ¶{3}\better ¶{4}\better ¶{264}\better ¶{3} and ¶{3}\better ¶{122}\better ¶{264}\better ¶{3}, which means that the competition was closely contested. In such close contests, the Skating System often has to resort to certain tie-breaking rules which are virtually equivalent to throwing the dice. In contrast, the all-round Llull matrix has the virtue of being a more accurate quantitative aggregate over the different dances. On~the basis of this more accurate aggregate, in this case the indirect scores reveal an all-round ranking which is quite different from the one produced by the Skating System (but it coincides with the one produced by other paired-comparison methods, like ranked pairs). In consonance with all this, the CLC rates obtained below are quite close to each other, particularly for the couples ¶{3},¶{4},¶{122} and ¶{264}.



Since we are dealing with complete votes, in this case the CLC computations can be carried out entirely in terms of the margins. In the following we have chosen to pass to margins after computing the indirect scores, but we could have done it before that step. 

\renewcommand\upla{\vskip-5pt} 

\bgroup\small
\initsizesmall

\bigskip\upla
\leftline{\hskip2\parindent
\taula{}%
\numx 3/4/31/122/238/264◊%
\hmulticol(6)\hskip3.5pt$V_{xy}$:
\hsk\sf 3\hsk|\hsk\sf 4\hsk|\hssk\sf 31\hssk|\sf 122|\sf 238|\sf 264:
\st|23|28|23|28|20/
21|\st|23|20|30|24/
16|21|\st|15|25|18/
21|24|29|\st|28|23/
16|14|19|16|\st|19/
24|20|26|21|25|\st◊%
\fitaula
}

\bigskip\upla
\leftline{\hskip2\parindent
\taula{}%
\numx 3/4/31/122/238/264◊%
\hmulticol(6)\hskip3.5pt\smash{\hbox{$V^\ast_{xy}$}}:
\hsk\sf 3\hsk|\hsk\sf 4\hsk|\hssk\sf 31\hssk|\sf 122|\sf 238|\sf 264:
\st|23|\xx28|23|\xx28|23/
\xx24|\st|\xx24|23|\xx30|\xx24/
21|21|\st|21|\xx25|21/
\xx24|\xx24|\xx29|\st|\xx28|\xx24/
19|19|19|19|\st|19/
\xx24|23|\xx26|23|\xx25|\st◊%
\hmulticol(1):$\kappa$: 4/ 2/ 5/ 1/ 6/ 3◊
\fitaula
}

\bigskip\upla
\leftline{\hskip2\parindent
\taula{}%
\numx 122/4/264/3/31/238◊%
\hmulticol(6)\hskip2pt\smash{\hbox{$M^\mast_{xy}$}}:
\sf 122|\hsk\sf 4\hsk|\sf 264|\hsk\sf 3\hsk|\hssk\sf 31\hssk|\sf 238:
\st|1|1|1|8|9/
\st|\st|1|1|3|11/
\st|\st|\st|1|5|6/
\st|\st|\st|\st|7|9/
\st|\st|\st|\st|\st|6/
\st|\st|\st|\st|\st|\st◊%
\fitaula
}

\bigskip\upla
\leftline{\hskip2\parindent
\taula{}%
\numx 122/4/264/3/31/238◊%
\hmulticol(6)\hskip2pt\smash{\hbox{$M^\pi_{xy}$}}:
\sf 122|\hsk\sf 4\hsk|\sf 264|\hsk\sf 3\hsk|\hssk\sf 31\hssk|\sf 238:
\st|1|1|1|3|6/
\st|\st|1|1|3|6/
\st|\st|\st|1|3|6/
\st|\st|\st|\st|3|6/
\st|\st|\st|\st|\st|6/
\st|\st|\st|\st|\st|\st◊%
\fitaula
\hskip2\parindent
\taula{}%
\numx 122/4/264/3/31/238◊%
\hmulticol(2):$r_x$|$\flr_x$:
3.3636|0.1815/
3.3864|0.1788/
3.4091|0.1761/
3.4318|0.1734/
3.5682|0.1583/
3.8409|0.1318◊%
\fitaula
}

\egroup

\paragraph{4.3} As a second example of an election involving truncated rankings we take the Debian Project leader election, which is using the method of paths since 2003.
So far, the winners of these elections have been clear enough. 
However, a quantitative measure of this clearness was lacking.
In the~following we consider the 2006 election, which had a participation of $V=421$ actual voters out of a total population of $972$ members. The individual votes are available in \texttt{\small http://www.debian.org/vote/2006/vote\underl002}.

\bigskip
That election resulted in the following Llull matrix:

\bgroup\small 
\initsizesmall

\bigskip
\leftline{\hskip2\parindent
\taula{}%
\numx 1/2/3/4/5/6/7/8◊%
\hmulticol(8)\hskip6pt\smash{\hbox{$V_{xy}$}}:
\sf 1|\sf 2|\sf 3|\sf 4|\sf 5|\sf 6|\sf 7|\sf 8:
\st|\hsk\xx321\hsk|144|159\pq12|\bfgr{193\pq12}|\xx347\pq12|\xx246|\xx320/
51|\st|42|53|50|\xx262|65|163/
\xx251|\xx340|\st|198\pq12|\xx253|\xx362|\xx300|\xx345/
\xx245\pq12|\xx341|\xx204\pq12|\st|\xx256|\xx371\pq12|\xx291\pq12|\xx339\pq12/
\bfgr{193\pq12}|\xx325|144|149|\st|\xx357|\xx254|\xx321\pq12/
26\pq12|77|24|22\pq12|21|\st|30|74\pq12/
137|\xx292|90|109\pq12|131|\xx330|\st|\xx296/
76|\xx207|54|71\pq12|75\pq12|\xx302\pq12|89|\st◊%
\fitaula
}

\egroup

\bigskip
\noindent
Notice that candidate ¶{4} is the winner according to the Condorcet principle (but not according to the majority principle, since $V_{\textsf{43}}$ does not reach $V/2$). Notice also that there is no Condorcet cycle. However, candidates ¶{1} and ¶{5} are in a tie for third place: both of them defeat all other candidates except ¶{4} and ¶{3}, and $V_{\textsf{15}}$ coincides exactly with $V_{\textsf{51}}$.

\bigskip
The ensuing CLC computations are as follows:

\bgroup\small 
\initsizesmall

\bigskip
\leftline{\hskip2\parindent
\taula{}%
\numx 1/2/3/4/5/6/7/8◊%
\hmulticol(8)\hskip6pt\smash{\hbox{$V^\ast_{xy}$}}:
\sf 1|\sf 2|\sf 3|\sf 4|\sf 5|\sf 6|\sf 7|\sf 8:
\st|\hsk\xx321\hsk|159\pq12|159\pq12|\bfgr{193\pq12}|\xx347\pq12|\xx246|\xx320/
89|\st|89|89|89|\xx262|89|163/
\xx251|\xx340|\st|198\pq12|\xx253|\xx362|\xx300|\xx345/
\xx245\pq12|\xx341|\xx204\pq12|\st|\xx256|\xx371\pq12|\xx291\pq12|\xx339\pq12/
\bfgr{193\pq12}|\xx325|159\pq12|159\pq12|\st|\xx357|\xx254|\xx321\pq12/
77|77|77|77|77|\st|77|77/
137|\xx292|137|137|137|\xx330|\st|\xx296/
89|\xx207|89|89|89|\xx302\pq12|89|\st◊%
\hmulticol(1):$\kappa$: 3\pq12/ 7/ 2/ 1/ 3\pq12/ 8/ 5/ 6◊
\fitaula
}

\bigskip
\leftline{\hskip2\parindent
\taula{}%
\numx 4/3/1/5/7/8/2/6◊%
\hmulticol(8)\hskip3pt\smash{\hbox{$M^\mast_{xy}$}}:
\sf 4|\sf 3|\sf 1|\sf 5|\sf 7|\sf 8|\sf 2|\sf 6:
\hph\st\hph|\hph6\hph|86|96\pq12|154\pq12|250\pq12|252|294\pq12/
\st|\st|91\pq12|93\pq12|163|256|251|285/
\st|\st|\st|0|109|231|232|270\pq12/
\st|\st|\st|\st|117|232\pq12|236|280/
\st|\st|\st|\st|\st|207|203|253/
\st|\st|\st|\st|\st|\st|44|225\pq12/
\st|\st|\st|\st|\st|\st|\st|185/
\st|\st|\st|\st|\st|\st|\st|\st◊%
\fitaula
}

\bigskip
\leftline{\hskip2\parindent
\taula{}%
\numx 4/3/1/5/7/8/2/6◊%
\hmulticol(8)\hskip6pt\smash{\hbox{$T_{xy}$}}:
\sf 4|\sf 3|\sf 1|\sf 5|\sf 7|\sf 8|\sf 2|\sf 6:
\hph\st\hph|403|405|405|401|411|394|394/
\st|\st|395|397|390|399|382|386/
\st|\st|\st|387|383|396|372|374/
\st|\st|\st|\st|385|397|375|378/
\st|\st|\st|\st|\st|385|357|360/
\st|\st|\st|\st|\st|\st|370|377/
\st|\st|\st|\st|\st|\st|\st|339/
\st|\st|\st|\st|\st|\st|\st|\st◊%
\fitaula
}

\bigskip
\leftline{\hskip2\parindent
\taula{}%
\numx 4/3/1/5/7/8/2/6◊%
\hmulticol(8)\hskip3pt $M^\sigma_{xy}$:
\sf 4|\sf 3|\sf 1|\sf 5|\sf 7|\sf 8|\sf 2|\sf 6:
\st|6|86|96\pq12|154\pq12|250\pq12|252|294\pq12/
\st|\st|86|93\pq12|154\pq12|250\pq12|251|285/
\st|\st|\st|0|109|231|232|270\pq12/
\st|\st|\st|\st|109|231|232|270\pq12/
\st|\st|\st|\st|\st|203|203|253/
\st|\st|\st|\st|\st|\st|44|225\pq12/
\st|\st|\st|\st|\st|\st|\st|185/
\hph\st\hph|\hph\st\hph|\hph\st\hph|\hph\st\hph|\st|\st|\st|\st◊%
\fitaula
}

\newpage 

\bigskip
\leftline{\hskip2\parindent
\taula{}%
\numx 4/3/1/5/7/8/2/6◊%
\hmulticol(8)\hskip6pt\smash{\hbox{$T^\sigma_{xy}$}}:
\sf 4|\sf 3|\sf 1|\sf 5|\sf 7|\sf 8|\sf 2|\sf 6:
\st|403.4|403.4|403.4|403.25|403.25|392|392/
\st|\st|397.4|397.4|397.25|397.25|386|386/
\st|\st|\st|389.6|389.6|389.6|374.75|374.75/
\st|\st|\st|\st|389.6|389.6|374.75|374.75/
\st|\st|\st|\st|\st|385|366|366/
\st|\st|\st|\st|\st|\st|366|366/
\st|\st|\st|\st|\st|\st|\st|339/
\hph\st\hph|\st|\st|\st|\st|\st|\st|\st◊%
\fitaula
}

\bigskip
\leftline{\hskip2\parindent
\taula{}%
\numx 4/3/1/5/7/8/2/6◊%
\hmulticol(8)\hskip3pt$M^\pi_{xy}$:
\sf 4|\sf 3|\sf 1|\sf 5|\sf 7|\sf 8|\sf 2|\sf 6:
\st|6|86|86|109|203|203|217/
\st|\st|86|86|109|203|203|217/
\st|\st|\st|0|109|203|203|217/
\st|\st|\st|\st|109|203|203|217/
\st|\st|\st|\st|\st|203|203|217/
\st|\st|\st|\st|\st|\st|44|185/
\st|\st|\st|\st|\st|\st|\st|185/
\hph\st\hph|\hph\st\hph|\hph\st\hph|\hph\st\hph|\st|\st|\st|\st◊%
\fitaula
}

\bigskip
\leftline{\hskip2\parindent
\taula{}%
\numx 4/3/1/5/7/8/2/6◊%
\hmulticol(8)\hskip6pt\smash{\hbox{$T^\pi_{xy}$}}:
\sf 4|\sf 3|\sf 1|\sf 5|\sf 7|\sf 8|\sf 2|\sf 6:
\st|403.4|397.4|397.4|389.6|385|385|371/
\st|\st|397.4|397.4|389.6|385|385|371/
\st|\st|\st|389.6|389.6|385|385|371/
\st|\st|\st|\st|389.6|385|385|371/
\st|\st|\st|\st|\st|385|385|371/
\st|\st|\st|\st|\st|\st|366|339/
\st|\st|\st|\st|\st|\st|\st|339/
\hph\st\hph|\st|\st|\st|\st|\st|\st|\st◊%
\fitaula
}

\newpage 

\bigskip
\leftline{\hskip1\parindent 
\taula{}%
\numx 4/3/1/5/7/8/2/6◊%
\hmulticol(8)\hskip6pt\smash{\hbox{$V^\pi_{xy}$}}:
\sf 4|\sf 3|\sf 1|\sf 5|\sf 7|\sf 8|\sf 2|\sf 6:
\st|204.7|241.7|241.7|249.3|294|294|294/
198.7|\st|241.7|241.7|249.3|294|294|294/
155.7|155.7|\st|194.8|249.3|294|294|294/
155.7|155.7|194.8|\st|249.3|294|294|294/
140.3|140.3|140.3|140.3|\st|294|294|294/
91|91|91|91|91|\st|205|262/
91|91|91|91|91|161|\st|262/
77|77|77|77|77|77|77|\st◊%
\fitaula
\hskip.5\parindent
\taula{}%
\hmulticol(2):$r_x$|$\flr_x$:
3.6784|0.2067/
3.6926|0.2048/
4.1105|0.1596/
4.1105|0.1596/
4.5720|0.1218/
5.8100|0.0599/
5.9145|0.0559/
6.7197|0.0317◊%
\fitaula
}

\egroup 

\bigskip
As one can see, the CLC results are in full agreement with the Copeland scores of the original Llull matrix. In particular, they still give an exact tie between candidates~¶{1} and ¶{5}. Even so,
the CLC rates yield a quantitative information which is not present in the Copeland scores. In particular, they show that the victory of candidate~¶{4} over candidate~¶{3} was relatively narrow.


\smallskip
For the computation of the rates we have taken $V=421$ (the actual number of votes) instead of $V=972$ (the number of people with the right to vote); in particular, the fraction-like rates $\flr_x$ have been computed so that they add up to~$\nev=1$ instead of the true participation ratio $\nev=421/972$.\linebreak 
This is especially justified in Debian elections since they systematically\linebreak 
include ``none~of~the~above'' as one of the alternatives, so it is reasonable to~interpret that abstention does not have a critical character. In~the present case, ``none~of~the~above'' was alternative~¶{8}, which 
obtained a~better result than two of the real candidates.

\paragraph{4.4} Finally, we look at an example of approval voting. Specifically, we consider the 2006 Public Choice Society election~\cite{pcs06}. Besides an approval vote, here the voters were also asked for a preferential vote ``in the spirit of research on public choice''. However,  here we will limit ourselves to the approval vote, which was the official one. 
The vote had a participation of $V = 37$ voters, most of which approved more than one candidate.

The actual votes are listed in the following table,\footnote{We are grateful to Prof.\ Steven J.\ Brams, who was the president of the Public Choice Society when that election took place, for his kind permission to reproduce these data.} where we give not only the approval voting data but also the associated preferential votes. The approved candidates are the ones which lie at the left of the slash.

\newcommand\fiav{/}

\smallskip
\bigskip
\vbox{%
\leavevmode
\hbox to125mm{\hrulefill}\par

\bgroup\small
\centerline{%
\hfil
\vtop{\noindent
\halign{\strut#\hfill\cr
¶{A}\better ¶{B}\fiav\cr 
¶{A}\better ¶{C}\better ¶{B}\fiav\cr 
¶{D}\fiav\better ¶{A}\better ¶{B}\better ¶{E}\better ¶{C}\cr 
¶{B}\better ¶{A}\fiav\better ¶{D}\better ¶{C}\better ¶{E}\cr 
¶{D}\better ¶{A}\better ¶{B}\better ¶{C}\fiav\better ¶{E}\cr 
¶{C}\better ¶{B}\better ¶{A}\fiav\cr 
¶{E}\fiav\better ¶{D}\cr 
¶{C}\better ¶{A}\better ¶{B}\better ¶{E}\fiav\cr 
¶{D}\better ¶{E}\fiav\better ¶{C}\better ¶{A}\better ¶{B}\cr 
¶{E}\fiav\cr 
¶{B}\better ¶{C}\fiav\cr 
¶{D}\better ¶{C}\fiav\better ¶{B}\better ¶{E}\better ¶{A}\cr 
¶{B}\fiav\cr 
}}
\hfil
\vtop{\noindent
\halign{\strut#\hfill\cr
¶{A}\fiav\cr 
¶{A}\fiav\cr 
¶{D}\fiav\better ¶{A}\tied ¶{B}\tied ¶{C}\tied ¶{E}\cr 
¶{A}\tied ¶{C}\fiav\cr 
\fiav¶{B}\better ¶{E}\better ¶{A}\better ¶{D}\better ¶{C}\cr 
¶{A}\tied ¶{B}\tied ¶{E}\fiav\cr 
¶{A}\tied ¶{B}\tied ¶{C}\tied ¶{D}\tied ¶{E}\fiav\cr 
¶{D}\better ¶{A}\better ¶{B}\fiav\cr 
¶{B}\better ¶{D}\better ¶{A}\fiav\better ¶{C}\better ¶{E}\cr 
¶{A}\fiav\better ¶{B}\better ¶{E}\better ¶{C}\better ¶{D}\cr 
¶{D}\fiav\cr 
¶{A}\tied ¶{C}\better ¶{B}\fiav\better ¶{D}\better ¶{E}\cr 
¶{A}\fiav\better ¶{D}\better ¶{B}\better ¶{C}\better ¶{E}\cr 
}}
\hfil
\vtop{\noindent
\halign{\strut#\hfill\cr
¶{C}\fiav\better ¶{B}\better ¶{D}\better ¶{A}\better ¶{E}\cr 
¶{C}\fiav\cr 
¶{D}\tied ¶{E}\fiav\better ¶{A}\better ¶{B}\tied ¶{C}\cr 
¶{B}\fiav\better ¶{C}\better ¶{A}\better ¶{D}\better ¶{E}\cr 
¶{D}\better ¶{C}\better ¶{E}\fiav\cr 
¶{C}\fiav\better ¶{A}\better ¶{B}\tied ¶{D}\tied ¶{E}\cr 
¶{C}\fiav\cr 
¶{B}\better ¶{D}\fiav\better ¶{E}\better ¶{C}\better ¶{A}\cr 
¶{B}\better ¶{C}\fiav\better ¶{A}\better ¶{E}\better ¶{D}\cr 
¶{D}\better ¶{A}\better ¶{C}\better ¶{B}\fiav\cr 
¶{D}\better ¶{E}\fiav\better ¶{A}\better ¶{B}\cr 
}}
\hfil}
\egroup

\leavevmode
\hbox to125mm{\hrulefill}
}

\smallskip
\bigskip
The approval voting scores are the following:
¶{A}:~17, ¶{B}:~16, ¶{C}:~17, ¶{D}:~14, ¶{D}:~9.
So according to approval voting there was a tie between candidates~¶{A} and ¶{C}, which were followed at a minimum distance by candidate~¶{B}.


\bigskip
The CLC computations are as follows:

\bgroup\small 
\initsizesmall

\bigskip
\leftline{\hskip2\parindent
\hbox to60mm{\ignorespaces
\taula{}%
\numx A/B/C/D/E◊%
\hmulticol(5)\hskip3pt\hbox{$V_{xy}$}:
\sf A|\sf B|\sf C|\sf D|\sf E:
\st|12\pq12|11|14|15\pq12/
11\pq12|\st|12|13\pq12|14\pq12/
11|13|\st|14\pq12|15\pq12/
11|11\pq12|11\pq12|\st|11\pq12/
7\pq12|7\pq12|7\pq12|6\pq12|\st◊%
\fitaula
\hfill}%
\taula{}%
\numx A/B/C/D/E◊%
\hmulticol(5)\hskip3pt\smash{\hbox{$V^\ast_{xy}$}}:
\sf A|\sf B|\sf C|\sf D|\sf E:
\st|\xx12\pq12|\xx12|\xx14|\xx15\pq12/
11\pq12|\st|12|\xx13\pq12|\xx14\pq12/
11\pq12|\xx13|\st|\xx14\pq12|\xx15\pq12/
11\pq12|11\pq12|11\pq12|\st|\xx11\pq12/
7\pq12|7\pq12|7\pq12|7\pq12|\st◊%
\hmulticol(1):$\kappa$: 1/ 3/ 2/ 4/ 5◊
\fitaula
}

\bigskip
\leftline{\hskip2\parindent
\hbox to60mm{\ignorespaces
\taula{}%
\numx A/C/B/D/E◊%
\hmulticol(5)\hskip3pt\smash{\hbox{$M^\mast_{xy}$}}:
\sf A|\sf C|\sf B|\sf D|\sf E:
\st|\pq12|1|2\pq12|8/
\st|\st|1|3|8/
\st|\st|\st|2|7/
\st|\st|\st|\st|4/
\hq\st\hq|\hq\st\hq|\hq\st\hq|\hq\st\hq|\hq\st\hq◊%
\fitaula
\hfill}%
\taula{}%
\numx A/C/B/D/E◊%
\hmulticol(5)\hskip3pt$T_{xy}$:
\sf A|\sf C|\sf B|\sf D|\sf E:
\st|22|24|25|23/
\st|\st|25|26|23/
\st|\st|\st|25|22/
\st|\st|\st|\st|18/
\hq\st\hq|\hq\st\hq|\hq\st\hq|\hq\st\hq|\hq\st\hq◊%
\fitaula
}

\newpage 

\bigskip
\leftline{\hskip2\parindent
\hbox to60mm{\ignorespaces
\taula{}%
\numx A/C/B/D/E◊%
\hmulticol(5)\hskip3pt\smash{\hbox{$M^\sigma_{xy}$}}:
\sf A|\sf C|\sf B|\sf D|\sf E:
\st|\pq12|1|2\pq12|8/
\st|\st|1|2\pq12|8/
\st|\st|\st|2|7/
\st|\st|\st|\st|4/
\hq\st\hq|\hq\st\hq|\hq\st\hq|\hq\st\hq|\hq\st\hq◊%
\fitaula
\hfill}%
\taula{}%
\numx A/C/B/D/E◊%
\hmulticol(5)\hskip3pt\smash{\hbox{$T^\sigma_{xy}$}}:
\sf A|\sf C|\sf B|\sf D|\sf E:
\st|24\pq12|24\pq12|24\pq12|22\pq78/
\st|\st|24\pq12|24\pq12|22\pq38/
\st|\st|\st|24\pq12|21\pq38/
\st|\st|\st|\st|19\pq38/
\hq\st\hq|\st|\st|\st|\st◊%
\fitaula
}

\bigskip
\leftline{\hskip2\parindent
\hbox to60mm{\ignorespaces
\taula{}%
\numx A/C/B/D/E◊%
\hmulticol(5)\hskip3pt\smash{\hbox{$M^\pi_{xy}$}}:
\sf A|\sf C|\sf B|\sf D|\sf E:
\st|\pq12|1|2|5.56/
\st|\st|1|2|5.56/
\st|\st|\st|2|5.56/
\st|\st|\st|\st|4/
\hq\st\hq|\hq\st\hq|\hq\st\hq|\hq\st\hq|\hq\st\hq◊%
\fitaula
\hfill}%
\taula{}%
\numx A/C/B/D/E◊%
\hmulticol(5)\hskip3pt\smash{\hbox{$T^\pi_{xy}$}}:
\sf A|\sf C|\sf B|\sf D|\sf E:
\st|24\pq12|24\pq12|24\pq12|20.94/
\st|\st|24\pq12|24\pq12|20.94/
\st|\st|\st|24\pq12|20.94/
\st|\st|\st|\st|19\pq38/
\hq\st\hq|\st|\st|\st|\st◊%
\fitaula
}

\bigskip
\leftline{\hskip2\parindent
\hbox to60mm{\ignorespaces
\taula{}%
\numx A/C/B/D/E◊%
\hmulticol(5)\hskip3pt\smash{\hbox{$V^\pi_{xy}$}}:
\sf A|\sf C|\sf B|\sf D|\sf E:
\st|12\pq12|12\pq34|13\pq14|13\pq14/
12|\st|12\pq34|13\pq14|13\pq14/
11\pq34|11\pq34|\st|13\pq14|13\pq14/
11\pq14|11\pq14|11\pq14|\st|11.69/
7.69|7.69|7.69|7.69|\st◊%
\fitaula
\hfill}%
\hskip4em
\taula{}%
\numx A/C/B/D/E◊%
\hmulticol(2):$r_x$|$\flr_x$:
3.6014|0.2315/
3.6149|0.2276/
3.6486|0.2181/
3.7720|0.1928/
4.1689|0.1299◊%
\fitaula
}

\egroup 

\bigskip
So, the winner by the CLC method is candidate~¶{A}. However, this is true only for the main variant. For the other three variants (dual, balanced and margin-based) the result is a tie between ¶{A} and ¶{C}, in full agreement with the approval voting scores. In~\secpar{17} we will see that the margin-based variant always gives such a full agreement.

\medskip
\noindent
\textit{Remark}. In all of the preceding examples, the matrix of the indirect scores has a constant row which corresponds to the loser. However, it  is not always~so.

\section{Some terminology and notation}

We consider a finite set~$\ist$. Its elements represent the options which are the matter of a vote. 
The number of elements of $\ist$ will be denoted by~$N$.
\ensep
We~will be particularly concerned with (binary) \dfc{relations} on~$\ist$.
Stating that two elements $a$ and $b$ are in a certain relation~$\rho$ is equivalent to
saying that the (ordered) {\df pair} formed by these two elements is a member of a certain set~$\rho$.
The pair formed by $a$ and $b$, in this order, will be denoted simply as~$ab$.

\smallskip
The pairs that consist of two copies of the same element, \ie those of the form~$aa$, are not relevant for our purposes. So, we will systematically exclude them from our relations. This can be viewed as a sort of normalization. The~set of all proper pairs, \ie the pairs $ab$ with $a\neq  b$, will be denoted as~$\tie$, or if necessary as $\tie(\ist)$. So, we will restrict our attention to relations contained in~$\tie$ (such relations are sometimes called ``strict'', or\linebreak[3] ``irreflexive'').
\ensep
In particular, the relation that includes the whole of~$\tie$ will be called {\df complete tie}.

\newcommand\ddd[1]{\halfsmallskip\noindent\hbox to 2\parindent{\hss\footnotesize$\bullet$\ \ }{#1}\ensep}

\smallskip
A relation $\rho\sbseteq\tie$ will be called\,:

\halfsmallskip
\ddd{\df total\rm, or complete,} when
at least one of $ab\in\rho$ and $ba\in\rho$ holds for every pair $ab$.
\ddd{\df antisymmetric} when
$ab\in\rho$ and $ba\in\rho$ cannot occur simultaneously.
\ddd{\df transitive} when
the simultaneous occurrence of $ab\in\rho$ and $bc\in\rho$
implies $ac\in\rho$.
\ddd{a {\df partial order}, when
it is at the same time transitive and antisymmetric.}
\ddd{a {\df total order}, or {\df strict ranking},} when
it is at the same time transitive, antisymmetric and total.
\ddd{a {\df complete ranking}} when
it is at the same time transitive and total.
\ddd{a {\df truncated ranking}} when it consists of a complete ranking on a subset $\xst$ of $\ist$ together with all pairs $ab$ with $a\in\xst$ and $b\not\in\xst$.


\medskip
For every relation $\rho\sbseteq\tie$, we will denote by $\rho'$ the relation that
consists of all pairs of the form $ab$ where $ba\in\rho$; $\rho'$~will be called
the {\df converse} of~$\rho$.
\ensep
On the other hand, we will denote by $\bar\rho$ the relation that consists of all
pairs $ab$ for which $ab\notin\rho$; $\bar\rho$~will be called the {\df complement}
of~$\rho$.
\ensep
For~certain purposes, it will be useful to consider also the relation $\hat\rho$ given by the complement of the converse of~$\rho$, or equivalently by the converse of its complement.
So, $ab\in\hat\rho$ \ifoi $ba\not\in\rho$.
This relation will be called the {\df adjoint} of $\rho$.
This operation will be used mainly in \secpar{8}, in connection with the indirect comparison relation~$\nu = \mu(\isc)$.
The following proposition collects several properties which are immediate consequences of the definitions:

\smallskip
\begin{lemma}\hskip.5em
\label{st:adjoint}

\iim{a}$\hat{\hat\rho} = \rho$.

\iim{b}$\rho\sbset\sigma
\ensep\Longleftrightarrow\ensep \hat\sigma\sbset\hat\rho$.

\iim{c}$\rho \text{ is antisymmetric}
\ensep\Longleftrightarrow\ensep \rho\sbseteq\hat\rho
\ensep\Longleftrightarrow\ensep \hat\rho \text{ is total}$.

\iim{d}$\rho \text{ is total}
\ensep\Longleftrightarrow\ensep \hat\rho\sbseteq\rho
\ensep\Longleftrightarrow\ensep \hat\rho \text{ is antisymmetric}$.
\end{lemma}

\begin{comment}
\noindent
$(\rho\cup\sigma)\,\hat{} = \hat\rho \cap \hat\sigma$;\quad
$(\rho\cap\sigma)\,\hat{} = \hat\rho \cup \hat\sigma$;\quad
$\rho\cup\hat\rho$ is total;\quad
$\rho\cap\hat\rho$ is antisymmetric;\hfil\break
$\rho$ transitive $\Rightarrow$ the partial order $\rho\cap\hat\rho$ is compatible with the equivalence relation $\rho\cap\rho'$.\par
\end{comment}

\medskip
Besides pairs, we will be concerned also with longer sequences $a_0 a_1
\dots a_n$. They will be referred to as {\df paths}, and in the case
$a_n = a_0$ they are called {\df cycles}. When $a_ia_{i+1}\in\rho$ for every $i$,
we will say that the path $a_0 a_1 \dots a_n$ is {\df contained in~$\rho$},
and also that $a_0$ and $a_n$ are {\df indirectly related through~$\rho$}.
\ensep
When $\rho$ is transitive, the condition ``$a$~is indirectly related to $b$
through~$\rho$'' implies $ab\in\rho$. In general, however, it defines a new
relation, which is called the {\df transitive closure} of~$\rho$, and will be
denoted by $\tcl\rho$;
this is 
the minimum transitive relation that contains $\rho$.
\ensep
The transitive-closure operator is easily seen to have the following
properties:\ensep
$\tcl\rho\sbseteq\tcl\sigma$ whenever $\rho\sbseteq\sigma$;\ensep
$\tclp{\rho\cap\sigma}\sbseteq(\tcl\rho)\cap(\tcl\sigma)$;\ensep
$(\tcl\rho)\cup(\tcl\sigma)\sbseteq\tclp{\rho\cup\sigma}$;\ensep
$\tclp{\tcl\rho}=\tcl\rho$.
\ensep
On the other hand, one can easily check that

\smallskip
\begin{lemma} 
\label{st:trancycle}
The transitive closure $\tcl\rho$ is antisymmetric \ifoi $\rho$ contains no cycle.
More specifically, $ab,ba\in\tcl\rho$ \ifoi
$\rho$ contains a cycle that includes both~$a$ and~$b$.
\end{lemma}

\medskip
A subset $\cst\sbseteq\ist$ will be said to be a {\df cluster} for a relation~$\rho$ when, for any~$x\not\in\cst$, having $ax\in\rho$ for some $a\in\cst$ implies $bx\in\rho$ for any $b\in\cst$, and similarly, having $xa\in\rho$ for some $a\in\cst$ implies $xb\in\rho$ for any $b\in\cst$.
\ensep
On~the other hand, $\cst\sbseteq\ist$ will be said to be an {\df interval} for a relation~$\rho$ when the simultaneous occurrence of $ax\in\rho$ and $xb\in\rho$ with $a,b\in\cst$ implies $x\in\cst$.
\ensep
The following facts are easy consequences of the definitions:
\ensep
If~$\rho$ is antisymmetric and $\cst$ is a cluster for~$\rho$ then $\cst$ is also an interval for~$\rho$.
\ensep
If~$\rho$ is total and $\cst$ is an interval for~$\rho$ then $\cst$ is also a cluster for~$\rho$.
\ensep
As a corollary, if~$\rho$ is total and antisymmetric, then $\cst$ is a cluster for~$\rho$ \ifoi \,it is an interval for that relation.
\ensep
Later on we will make use of the following fact, which is also an easy consequence of the definitions:
\begin{lemma}\hskip.5em
\label{st:clusterLemma}The following conditions are equivalent to each other:

\iim{a} $\cst$~is a cluster for~$\rho$.

\iim{b} $\cst$~is a cluster for~$\hat\rho$.

\iim{c}The simultaneous occurrence of $ax\in\rho$ and $xb\in\hat\rho$ with $a,b\in\cst$ implies $x\in\cst$, and similarly, the simultaneous occurrence of $ax\in\hat\rho$ and $xb\in\rho$ with $a,b\in\cst$ implies also $x\in\cst$.
\end{lemma}

\medskip
When $\cst$ is a cluster for~$\rho$, it will be useful to consider a new set~$\istbis$ and a new relation $\rhobis$ defined in the following way:
\ensep
$\istbis$ is obtained from $\ist$
by replacing the set $\cst$ by a single element~$\clustit$,
\ie $\istbis = (\ist\setminus\cst)\cup\{\clustit\}$;
\ensep
for~$x,y\in\ist\setminus\cst$,
\ensep
$x\clustit\in\rhobis$ \ifoi there exists $\clone\in\cst$ such that
$x\clone\in\rho$,
\ensep
$\clustit y\in\rhobis$ \ifoi there exists $\clone\in\cst$ such that
$\clone y\in\rho$,
\ensep
and finally,\linebreak 
$xy\in\rhobis$ \ifoi $xy\in\rho$.
\ensep
We~will refer to this operation as the 
{\df con\-traction} of $\rho$ by the cluster $\cst$.
\ensep
If~$\rho$ is a strict ranking (resp.~a~complete ranking) on $\ist$,
then $\rhobis$ is a strict ranking (resp.~a~complete ranking) on~$\istbis$.

\medskip
Given a relation~$\rho$, we will associate every element~$x$ with the following sets:
\ddd{the {\df set of predecessors}, $\pre{x}$,} \ie the set of $y\in\ist$ such that $yx\in\rho$.
\ddd{the {\df set of successors}, $\suc{x}$,} \ie the set of $y\in\ist$ such that $xy\in\rho$.
\ddd{the {\df set of collaterals}, $\col{x}$,} \ie the set of $y\in\ist\setminus\{x\}$ which are neither predecessors nor successors of $x$ in $\rho$.

\medskip
The sets $\pre{x}$, $\suc{x}$ and $\col{x}$ are especially meaningful  when the relation $\rho$ is a partial order.
In that case, and it is quite natural to rank the elements of $\ist$ by their number of predecessors, or by the number of elements which are not their successors. More precisely, it makes sense to define the {\df rank} of $x$ in $\rho$ as
\vskip-8mm 
\begin{equation}
\rank x \,=\, 1 \,+\, |\pre{x}| \,+\, \vartheta\, |\col{x}|
\,=\, 1 \,+\, (1-\vartheta)\,|\pre{x}| \,+\, \vartheta\, (N - 1 - |\suc{x}|),
\label{eq:rank}
\end{equation}
\noindent
where $\vartheta$ is a fixed number in the interval $0\le \vartheta\le 1$.
If we do not say otherwise, we will take $\vartheta=\onehalf$.
\ensep
The following facts are easy consequences of the definitions:

\begin{lemma}\hskip.5em
\label{st:rank}
Assume that $\rho$ is a partial order.
In that case, having $xy\in\rho$ implies the following facts: $\pre{x}\sbset\pre{y}$, $\suc{x}\spset\suc{y}$ (both inclusions are strict), and $\rank x < \rank y$ (for any $\vartheta$ in the interval $0\le \vartheta\le 1$). For $\vartheta=\onehalf$, the~average of the numbers $\rank x$ is equal to $(N+1)/2$. If $\rho$~is a total order, then $\rank x$ does not depend on $\vartheta$; furthermore, having $xy\in\rho$ is then equivalent to $\rank x < \rank y$.
\end{lemma}

\medskip
As in \secpar{2.2}, given a set of binary scores $s_{xy}$, we denote by $\mu(s)$ the corresponding {\df comparison relation}:
\vskip-5mm 
\begin{equation}
xy\in\mu(s) \quad\equiv\quad s_{xy} > s_{yx}.
\end{equation}
\noindent
For such a relation, the adjoint $\hat\mu(s)$ corresponds to replacing the strict inequality by the non-strict one.

\section{The indirect scores and its comparison relation}

Let us recall that the indirect scores~$\isc_{xy}$ are defined in the following way:
$$
\isc_{xy} \,=\, \max\,\{v_\alpha\mid \text{$\alpha$ is a path \,$x_0 x_1 \dots x_n$\, from $x_0=x$ to $x_n=y$}\,\},
$$
where the score~$v_\alpha$ of a path $\alpha=x_0 x_1 \dots x_n$ is defined as
$$
v_\alpha \,=\, \min\,\{v_{x_ix_{i+1}}\mid 0 \le i < n\,\}.
$$
In the following statements, and the similar ones which appear elsewhere, “any $x,y,z$” should be understood as meaning “any $x,y,z$ which are pairwise different from each other”.

\bigskip\noindent\textit{Remark}.\hskip.5em
The matrix of indirect scores $\isc$ can be viewed as a power of $v$ (supplemented with $v_{xx} = 1$) for a matrix product defined in the following way: $(vw)_{xz} = \max_y \min(v_{xy},w_{yz})$. More precisely, 
$\isc$ coincides with such a power for any exponent greater than or equal to~$N-1$.

\medskip
\begin{lemma}\hskip.5em
\label{st:minInequality}
The indirect scores satisfy the following inequalities:
\begin{equation}
\isc_{xz} \ge \min\,(\isc_{xy}, \isc_{yz})\quad\hbox{for any $x,y,z$.}
\end{equation}
\end{lemma}

\begin{proof}\hskip.5em
Let $\alpha$ be a path from $x$ to $y$ such that $\isc_{xy} = v_\alpha$; let $\beta$ be a path from $y$ to $z$ such that $\isc_{yz} = v_\beta$. Consider now their concatenation~$\alpha\beta$. Since $\alpha\beta$~goes from $x$ to $z$, one has $\isc_{xz} \ge v_{\alpha\beta}$. On the other hand, the definition of the score of a path ensures that $v_{\alpha\beta} = \min\,(v_\alpha, v_\beta)$. Putting these things together gives the desired result.
\end{proof}

\medskip
The following lemma is somehow a converse of the preceding one:
\begin{lemma}\hskip.5em
\label{st:indirectEqualToDirect}
Assume that the original scores satisfy the following inequalities:
\begin{equation}
v_{xz} \ge \min\,(v_{xy}, v_{yz})\quad\hbox{for any $x,y,z$.}
\label{eq:minInequalityIndirect}
\end{equation}
In that case, the indirect scores coincide with the original ones.
\end{lemma}

\begin{proof}\hskip.5em
The inequality $\isc_{xz}\ge v_{xz}$ is an immediate consequence of the definition of $\isc_{xz}$.
The converse inequality can be obtained in the following way: Let $\gamma=x_0x_1x_2 \dots x_n$ be a path from~$x$ to~$z$ such that $\isc_{xz} = v_\gamma$. By virtue of~(\ref{eq:minInequalityIndirect}), we have
$$
\min\,\left(\,v_{x_0x_1}, v_{x_1x_2}, v_{x_2x_3}, \dots, v_{x_{n-1}x_n}\right)
\,\,\le\,\,
\min\,\left(\,v_{x_0x_2}, v_{x_2x_3}, \dots, v_{x_{n-1}x_n}\right).
$$
So, $\isc_{xz} \le v_{\gamma'}$ where $\gamma'=x_0x_2 \dots x_n$. By iteration, one eventually gets $\isc_{xz} \le v_{xz}$.
\end{proof}

\begin{comment}
\noindent
Remark: In order to obtain $\isc_{xz}=v_{xz}$ it suffices that inequality~(\ref{eq:minInequalityIndirect}) holds for all pairs of the form $xy$ or alternatively for all those of the form $yz$.\par
\end{comment}

\medskip
\begin{theorem}[Schulze, 1998
\hbox{\brwrap{\textsl{34\,b}}}
]\hskip.5em
\label{st:transSchulze}
$\mu(\isc)$ is a transitive relation.
\end{theorem}

\begin{proof}\hskip.5em
We will argue by contradiction. Let us assume that\linebreak $xy\in\mu(\isc)$ and $yz\in\mu(\isc)$, but $xz\notin\mu(\isc)$. This means respectively that\ensep
(a)~$\isc_{xy} > \isc_{yx}$\ensep
and\, (b)~$\isc_{yz} > \isc_{zy}$,\ensep
but\, (c)~$\isc_{zx} \ge \isc_{xz}$.
On the other hand, Lemma~\ref{st:minInequality} ensures also that\ensep (d)~$\isc_{xz} \ge \min\,(\isc_{xy}, \isc_{yz})$. 
We~will\linebreak distinguish two cases depending on which of the two last quantities is smaller:\ensep
(i)~$\isc_{yz} \ge \isc_{xy}$;\ensep
(ii)~$\isc_{xy} \ge \isc_{yz}$.

\halfsmallskip
Case~(i)\,: $\isc_{yz} \ge \isc_{xy}$.\ensep
We will see that in this case (c) and (d) entail a contradiction with~(a).
In fact, we have the following chain of inequalities:
$\isc_{yx} \ge
\min\,(\isc_{yz}, \isc_{zx}) \ge
\min\,(\isc_{yz}, \isc_{xz}) \ge
\min\,(\isc_{yz}, \isc_{xy}) = 
\isc_{xy}$,\ensep
where we are using successively:
Lemma~\ref{st:minInequality}, (c), (d) and (i).

\halfsmallskip
Case~(ii)\,: $\isc_{xy} \ge \isc_{yz}$.\ensep
An entirely analogous argument shows that in this case (c) and (d) entail a contradiction with~(b). In fact, we have 
$\isc_{zy} \ge
\min\,(\isc_{zx}, \isc_{xy}) \ge
\min\,(\isc_{xz}, \isc_{xy}) \ge
\min\,(\isc_{yz}, \isc_{xy}) = 
\isc_{yz}$,\ensep
where we are using successively:
Lemma~\ref{st:minInequality}, (c), (d) and (ii).
\end{proof}

\section{Restricted paths}

In this section we consider paths restricted to either $\mu(v)$ or $\hat\mu(v)$.
Such restricted paths allow to achieve not only the majority principle \llmpw, but also the Condorcet principle~\llcpw. In exchange, however, this idea can hardly be made into a continuous rating method, since one is doing quite different things depending on whether $v_{xy} > v_{yx}$ or $v_{xy} < v_{yx}$. Even so, we will see that in the complete case ---where \llmpw\ is equivalent to~\llcpw--- the indirect comparison relations which are obtained under such restrictions coincide with the one which is obtained when arbitrary paths are used. More specifically, we will look at the comparison relations associated with 
$\uu^\ast_
{xy}$ and $\ww^\ast_
{xy}$, where $\uu_
{xy}$ and $\ww_
{xy}$ are defined as
\begin{equation}
\uu_
{xy} = \begin{cases}
v_{xy}, & \text{if } v_{xy} > v_{yx},\\
0, & \text{otherwise};\\
\end{cases}
\qquad
\ww_
{xy} = \begin{cases}
v_{xy}, & \text{if } v_{xy} \ge v_{yx},\\
0, & \text{otherwise}.\\
\end{cases}
\label{eq:decomposition}
\end{equation}

\smallskip
\begin{proposition}\hskip.5em
\label{st:restricted}

\iim{a}$\mu(\uu^\ast) \,\sbseteq\, \mu^\ast(v).$

\iim{b}$\mu(\ww^\ast) \,\sbseteq\, \hat\mu^\ast(v).$
\end{proposition}

\begin{proof}\hskip.5em
Part~(a). Let us begin by recalling that $\mu^\ast(v)$ means the transitive closure of~$\mu(v)$. Let~us assume that $xy\in\mu(\uu^\ast)$, \ie  $\uu^\ast_
{xy} > \uu^\ast_
{yx}$. Since we are dealing with non-negative numbers, this ensures that $\uu^\ast_
{xy} > 0$. By the definition of~$\uu^\ast_
{xy}$, this implies the existence of a path $x_0x_1\dots x_n$ from $x_0=x$ to $x_n=y$ such that $\uu_
{x_ix_{i+1}} > 0$ for all $i$. According to~(\ref{eq:decomposition}.1), this ensures that $v_{x_ix_{i+1}} > v_{x_{i+1}x_i}$, \ie $x_ix_{i+1}\in\mu(v)$, for~all $i$. Therefore, $xy\in\mu^\ast(v)$.
\ensep
An entirely analogous argument proves part~(b).
\end{proof}

\smallskip
\begin{lemma}\hskip.5em
\label{st:laiarosa}

\iim{a}$\uu^\ast_
{xy} \le \ww^\ast_
{xy} \le \isc_{xy}.$

\iim{b}$\isc_{xy} > 1/2
\ensep\Longrightarrow\ensep
\uu^\ast_{xy} = \ww^\ast_{xy} = \isc_{xy}.$ 

\iim{c}$\isc_{xy} = 1/2 
\ensep\Longrightarrow\ensep
\ww^\ast_{xy} = \isc_{xy}.$

\vskip2pt
\noindent
In the complete case one has:

\iim{d}$\isc_{xy} < 1/2
\ensep\Longrightarrow\ensep
\uu^\ast_{xy} = \ww^\ast_{xy} = 0.$ 

\iim{e}$\isc_{xy} = 1/2 
\ensep\Longrightarrow\ensep
\uu^\ast_{xy} = 0.$

\end{lemma}

\begin{proof}\hskip.5em
Part~(a).\ensep It is simply a matter of noticing that $\uu_{xy} \le \ww_{xy} \le v_{xy}$ and checking that the inequality $p_{xy}\le q_{xy}$ for all $x,y$ implies $p^\ast_{xy}\le q^\ast_{xy}$ for all $x,y$. As an intermediate result towards this implication, one can see that $p_\gamma\le q_\gamma$ for all paths $\gamma$. In fact, if $\gamma=x_0x_1 \dots x_n$ and $i$~is such that $q_\gamma = q_{x_ix_{i+1}}$, the definition of $p_\gamma$ and the inequality between $p_{xy}$ and $q_{xy}$ give $p_\gamma \le p_{x_ix_{i+1}} \le q_{x_ix_{i+1}} = q_\gamma$.\ensep
The second step of that implication uses an analogous argument: if $\gamma$ is a path from $x$ to $y$ such that $p^\ast_{xy} = p_\gamma$, we can write $p^\ast_{xy} = p_\gamma \le q_\gamma \le q^\ast_{xy}$,
where we are using the intermediate result and the definition of $q^\ast_{xy}$.

\halfsmallskip
Part~(b).\ensep Let $\gamma=x_0x_1 \dots x_n$ be a path from $x$ to $y$ such that $\isc_{xy} = v_\gamma$.\linebreak 
Since $\isc_{xy} > 1/2$, every link of that path satisfies $v_{x_ix_{i+1}}>1/2$, which implies that $v_{x_ix_{i+1}}>v_{x_{i+1}x_i}$ (because $v_{x_ix_{i+1}}+v_{x_{i+1}x_i}\le 1$). Now, that inequality entails that $\uu_{x_ix_{i+1}} = v_{x_ix_{i+1}}$, from which it follows that $\uu_\gamma = v_\gamma$. Finally, it suffices to combine these facts with the inequality $\uu_\gamma\le\uu^\ast_{xy}$ and the inequalities of part~(a):
$$v_{xy}^*=v_\gamma =\uu_\gamma\le\uu_{xy}^*\le\ww_{xy}^* \le v_{xy}^*.$$

\halfsmallskip
Part~(c).\ensep The proof is similar to that of part~(b). Here we deal with the non-strict inequality $v_{x_ix_{i+1}}\ge 1/2$, which entails $v_{x_ix_{i+1}}\ge v_{x_{i+1}x_i}$ and $\ww_{x_ix_{i+1}} = v_{x_ix_{i+1}}$. These facts allow to conclude  that
$$v_{xy}^*=v_\gamma =\ww_\gamma\le\ww_{xy}^* \le v_{xy}^*.$$

\halfsmallskip
Part~(d).\ensep The hypothesis that $\isc_{xy} < 1/2$ means that for every path $\gamma=x_0x_1 \dots x_n$ from $x$ to $y$ there exists at least one $i$ such that $v_{x_ix_{i+1}}< 1/2$. By the assumption of completeness, this implies that $v_{x_ix_{i+1}}<v_{x_{i+1}x_i}$, so that $\uu_{x_ix_{i+1}}=\ww_{x_ix_{i+1}}=0$. This implies that $\uu_\gamma=\ww_\gamma=0$. Since $\gamma$ is arbitrary, it follows that 
$\uu^\ast_{xy}=\ww^\ast_{xy}=0$.

\halfsmallskip
Part~(e).\ensep The proof is similar to that of part~(d). Here we deal with the non-strict inequality $v_{x_ix_{i+1}}\le 1/2$, which implies $v_{x_ix_{i+1}}\le v_{x_{i+1}x_i}$ and $\uu_{x_ix_{i+1}} = 0$. This holds for at least one link of every path $\gamma$ from $x$ to $y$. So, $\uu^\ast_{xy}=0$.
\end{proof}

\smallskip
\begin{theorem}\hskip.5em
\label{st:laia}
In the complete case one has \,$\mu(\uu^\ast) = \mu(\ww^\ast) = \mu(\isc)$.
\end{theorem}

\begin{proof}\hskip.5em
It suffices to prove the three following statements:
\begin{gather}
\label{eq:laia1}
\isc_{xy} > \isc_{yx}
\ensep\Longrightarrow\ensep
\uu^\ast_{xy} > \uu^\ast_{yx}
\quad\hbox{and}\quad
\ww^\ast_
{xy} > \ww^\ast_
{yx}
\\[2.5pt]
\label{eq:laia2}
\ww^\ast_{xy} > \ww^\ast_{yx}
\ensep\Longrightarrow\ensep
\isc_{xy} > \isc_{yx},
\\[2.5pt]
\label{eq:laia3}
\uu^\ast_{xy} > \uu^\ast_{yx}
\ensep\Longrightarrow\ensep
\isc_{xy} > \isc_{yx},
\end{gather}

\halfsmallskip
Proof of~(\ref{eq:laia1}).\ensep
By combining the completeness assumption with the hypothesis of (\ref{eq:laia1}) we can write $1 = v_{xy}+v_{yx}\le \isc_{xy}+\isc_{yx}< 2 \isc_{xy}$, so that $\isc_{xy} > 1/2$. According to part~(b) of Lemma~\ref{st:laiarosa}, this inequality implies that $\uu^\ast_{xy} = \ww^\ast_{xy} = \isc_{xy}$. On the other hand, part~(a) of the same lemma ensures that $\uu^\ast_{yx} \le \ww^\ast_{yx} \le \isc_{yx}$. By combining these facts with the hypothesis of (\ref{eq:laia1}) we obtain the right-hand side of it.

\halfsmallskip
Proof of~(\ref{eq:laia2}).\ensep
Here we begin by noticing that the left-hand side implies $\ww^\ast_{xy} > 0$, which by part~(d) of Lemma~\ref{st:laiarosa} entails $\isc_{xy} \ge 1/2$. If $\isc_{yx}<1/2$, we are finished. If, on the contrary, $\isc_{yx}\ge 1/2$, then parts~(a), (b) and (c) of Lemma~\ref{st:laiarosa} allow to conclude that $\isc_{yx}=\ww^\ast_{yx}<\ww^\ast_{xy}\le \isc_{xy}.$

\halfsmallskip
Proof of~(\ref{eq:laia3}).\ensep
Similarly to above, the left-hand side implies $\uu^\ast_{xy} > 0$, which by parts~(d) and (e) of Lemma~\ref{st:laiarosa} entails $\isc_{xy} > 1/2$. If $\isc_{yx}\le 1/2$, we are finished. If, on the contrary, $\isc_{yx}> 1/2$, then parts~(a) and~(b) of Lemma~\ref{st:laiarosa} allow to conclude that $\isc_{yx}=\uu^\ast_{yx}<\uu^\ast_{xy}\le \isc_{xy}.$
\end{proof}

\section{Admissible orders}

Let us recall that an admissible order is a total order $\xi$ such that $\nu \sbseteq \xi \sbseteq \hat\nu$.
Here $\nu$ is the indirect comparison relation $\nu = \mu(\isc)$. So $xy\in\nu$ \ifoi $\img_{xy} = \isc_{xy} - \isc_{yx} > 0$,\ensep and $xy\in\hat\nu$ \ifoi $\img_{xy} \ge 0$.

\bigskip
\begin{lemma}\hskip.5em
\label{st:ore}
Assume that $\rho$ is an antisymmetric and transitive relation. If~$\rho $ contains neither~$xy$ nor~$yx$, then $\tcl{(\rho\cup\{xy\})}$ is also antisymmetric.
\end{lemma}
\begin{proof}\hskip.5em
We will proceed by contradiction.\ensep According to Lemma~\ref{st:trancycle}, if\linebreak $\tcl{(\rho\cup\{xy\})}$ were not antisymmetric, $\rho\cup\{xy\}$ would contain a cycle~$\gamma$. On~the~other hand, the hypotheses on $\rho$ ensure, by the same lemma, that $\rho$~contains no cycles. Therefore, $\gamma$ must involve the pair $xy$. By following this cycle from one ocurrence of the pair $xy$ until the next ocurrence of $x$, one obtains a path from $y$ to $x$ which is contained in $\rho$. But, since $\rho$ is transitive, this entails that $yx\in\rho$, which contradicts one of the hypotheses.
\end{proof}

\smallskip
\begin{theorem}\hskip.5em
\label{st:existenceXiThm}
Given a transitive antisymmetric relation~$\rho$ on a finite set~$\ist$, one can always find a total order~$\xi$ such that $\rho \sbseteq \xi \sbseteq \hat\rho$.
If $\rho$ contains neither~$xy$ nor~$yx$, one can constrain $\xi$ to include the pair $xy$.
\end{theorem}
\begin{proof}\hskip.5em
If $\rho$ is total, it suffices to take $\xi=\rho$ (notice that $\hat\rho=\rho$ because of statements~(c) and (d) of Lemma~\ref{st:adjoint}).
\ensep
Otherwise, let us consider the relation $\rho_1 = \tcl{(\rho\cup\{xy\})}$, where $xy$ is any pair such that $\rho$ contains neither $xy$\linebreak nor $yx$. According to Lemma~\ref{st:ore}, $\rho_1$ is antisymmetric. Furthermore, it is obvious that $\rho\sbset\rho_1$. Therefore, the statements~(b) and (c) of Lemma~\ref{st:adjoint} ensure that $\rho \sbset \rho_1 \sbseteq \hat\rho_1 \sbset \hat\rho$.
\ensep
From here, one can repeat the same process with $\rho_1$ substituted for~$\rho$: if $\rho_1$ is total we take $\xi=\rho_1$; otherwise we consider $\rho_2 = \tcl{(\rho_1\cup\{x_1y_1\})}$, where $x_1y_1$ is any pair such that $\rho_1$ contains neither $x_1y_1$ nor $y_1x_1$, and so on. This iteration will conclude in a finite number of steps since $\ist$ is finite.
\end{proof}

\smallskip
\begin{corollary}\hskip.5em
\label{st:existenceXiCor}
One can always find an admissible order~$\xi$.
\end{corollary}

\begin{proof}\hskip.5em
It follows from Theorem~\ref{st:existenceXiThm} 
because $\nu=\mu(\isc)$ is certainly antisymmetric and
Theorem~\ref{st:transSchulze}
ensures that it is transitive.
\end{proof}

\medskip
Later on we will make use of the following fact:

\smallskip
\begin{theorem}\hskip.5em
\label{st:existenceXiCluster}
Given a transitive antisymmetric relation~$\rho$ on a finite set~$\ist$ and a set $\cst$ which is a cluster for~$\rho$, one can always find a total order~$\xi$ such that $\rho \sbseteq \xi \sbseteq \hat\rho$ and such that $\cst$ is a cluster for~$\xi$.
\end{theorem}
\begin{proof}\hskip.5em
As in the proof of Theorem~\ref{st:existenceXiThm}, we will progressively extend $\rho$ until we get a total order. Here, we will take care that besides being transitive and antisymmetric, the successive extensions $\rho_i$ keep the property that $\cst$ be a cluster for~$\rho_i$. To this effect, the successive additions to $\rho$ will follow a certain specific order, and we will make an extensive use of the necessary and sufficient condition given by Lemma~\ref{st:clusterLemma}.

\newcommand\cc{c}
\newcommand\dd{d}

\halfsmallskip
In a first phase we will deal with pairs of the form $\cc\dd$ with $\cc,\dd\in\cst$.\linebreak 
Let us assume that neither $\cc\dd$ nor $\dd\cc$ is contained in $\rho$, and let us consider $\rho_1=\left(\rho\cup\{\cc\dd\}\right)^*$. Besides the properties mentioned in the proof of Theorem~\ref{st:existenceXiThm}, we~claim that this relation has the property that $\cst$ is a cluster for~$\rho_1$. According to Lemma~\ref{st:clusterLemma}, it suffices to  check that the simultaneous occurrence of $ax\in\rho_1$ and $xb\in\hat\rho_1$ with $a,b\in\cst$ implies $x\in\cst$, and\linebreak[3] similarly, that the simultaneous occurrence of $ax\in\hat\rho_1$ and $xb\in\rho_1$ with $a,b\in\cst$ implies also $x\in\cst$. So, let us assume first that $ax\in\rho_1$ and $xb\in\hat\rho_1$ with $a,b\in\cst$. Since $\rho\sbset\rho_1$, we have $xb\in\hat\rho$ (because $\hat\rho_1\sbset\hat\rho$). If~$ax\in\rho$, we immediately get $x\in\cst$ since $\cst$ is known to be a cluster for~$\rho$ (Lemma~\ref{st:clusterLemma}). Otherwise, \ie if $ax\in\rho_1\setminus\rho$, we see that $\rho_1$ contains a path of the form $\gamma=a\dots \cc\dd\dots x$. But this entails the existence of a path from $\dd$ to $x$ contained in $\rho$. So, by transitivity, $\dd x\in\rho$. Again, this fact together with $xb\in\hat\rho$ ensures that $x\in\cst$\, since $\cst$ is known to be a cluster for~$\rho$. A~similar argument takes care of the case where $ax\in\hat\rho_1$ and $xb\in\rho_1$ with $a,b\in\cst$.

By repeating the same process we will eventually get an extension
of $\rho$ with the same properties plus the following one: it includes either $\cc\dd$ or $\dd\cc$ for any $\cc,\dd\in\cst$. In other words, its restriction to $\cst$ is a total order.
In the following, this relation will 
be denoted by $\rhoi$.

\halfsmallskip
Now we will deal with pairs of the form $\cc\qq$ or $\qq\cc$ with $\cc\in\cst$ and $\qq\not\in\cst$. Let us assume that neither $\cc\qq$ nor $\qq\cc$ belong to $\rhoi$. In this case we will
proceed by taking $\rhoi_1=\left(\rhoi\cup\{\last\qq\}\right)^*$, where $\last$ denotes the last element of~$\cst$ according to the total order determined by~$\rhoi$ (alternatively, one could take $\rhoi_1=\left(\rhoi\cup\{\qq f\}\right)^*$, where $f$ denotes the first element of~$\cst$ by $\rhoi$). By so doing, we make sure that $\rhoi_1$ contains all pairs of the form $z\qq$ with $z\in\cst$. As a consequence, $\cst$ will keep the property of being a cluster for~$\rhoi_1$. In~fact, let us assume, in the lines of Lemma~\ref{st:clusterLemma}, that $ax\in\rhoi_1$ and $xb\in\hat\rhoi_1$ with $a,b\in\cst$. The hypothesis that $ax\in\rhoi_1$ can be divided in two cases, namely either $ax\in\rhoi$ or $ax\in\rhoi_1\setminus\rhoi$. Let us consider first the case $ax\in\rhoi_1\setminus\rhoi$. By the definition of $\rhoi_1$, this means that $\rhoi_1$ contains a path of the form $\gamma=a\dots \last\qq\dots x$, whose final part shows that $\qq x\in\rhoi_1$. On the other hand, we know that $\rhoi_1$ contains $b\qq$ (since $b\in\cst$). By transitivity, this entails $bx\in\rhoi_1$ and therefore $xb\not\in\hat\rhoi_1$, in contradiction with the hypothesis that $xb\in\hat\rhoi_1$.
\ensep
So, the only possibility of having $ax\in\rhoi_1$ and $xb\in\hat\rhoi_1$ is $ax\in\rhoi$. Besides, $xb\in\hat\rhoi_1$ implies that $xb\in\hat\rhoi$. So $x\in\cst$ because $\cst$ is a cluster for~$\rhoi$ (Lemma~\ref{st:clusterLemma}).
\ensep
Let us assume now that $ax\in\hat\rhoi_1$ and $xb\in\rhoi_1$. 
Like before, the former implies $ax\in\hat\rhoi$.
Again, the hypothesis that $xb\in\rhoi_1$ can be divided in two cases, namely either $xb\in\rhoi$ or $xb\in\rhoi_1\setminus\rhoi$.
In the first case we have $ax\in\hat\rhoi$ and $xb\in\rhoi$. So $x\in\cst$ because $\cst$ is a cluster for~$\rhoi$. In~the second case we can still use the same argument since $\rhoi_1$ contains a~path of the form $\gamma=x\dots \last\qq\dots b$, which shows that $x\last\in\rhoi$.

By repeating the same process we will eventually get an extension
of $\rhoi$ with the same properties plus the following one: it includes either $\cc\qq$ or $\qq\cc$ for any $\cc\in\cst$ and $\qq\not\in\cst$.

\halfsmallskip
Finally, it rests to deal with any pairs of the form $\pp\qq$ with $\pp,\qq\not\in\cst$. However, these pairs do not cause any problems since they do not appear in the definition of $\cst$ being a cluster.
\end{proof}

\medskip
In practice, one can easily obtain admissible orders by suitably arranging the elements of~$\ist$ according to their number of victories, ties and defeats against the others 
according to the indirect comparison relation $\nu$.
More precisely, it suffices to arrange the elements of~$\ist$ by non-decreasing values of their rank~$\rank x$ in $\nu$ as defined in~(\ref{eq:rank}).
\ensep
According to the the particular nature of~$\nu$ and the definitions given in~\secpar{5}, the sets $\pre{x}$, $\suc{x}$ and $\col{x}$ which appear in~(\ref{eq:rank}) are given by
\begin{align}
\pre{x} &= \{\,y\mid y\neq x,\ \img_{xy} < 0\,\}, \\
\suc{x} &= \{\,y\mid y\neq x,\ \img_{xy} > 0\,\}, \\
\col{x} &= \{\,y\mid y\neq x,\ \img_{xy} = 0\,\}.
\end{align}
So, ranking by~$\rank x$ amounts to applying the Copeland rule to the tournament defined by the indirect comparison relation $\nu = \mu(\isc)$ (see for instance \cite[p.\,206--209]{t6}).

\smallskip
\begin{proposition}\hskip.5em
\label{st:Copeland}
Any total ordering of the elements of $\ist$ by non-decreasing values of $\ranka x\nu$ is an admissible order. This is true for any fixed value of $\vartheta$ in the interval $0\le \vartheta\le 1$.
\end{proposition}
\begin{proof}\hskip.5em
Let $\xi$ be a total order of~$\ist$ for which $x\mapsto\rank x$ does not decrease. This means that
$$
xy\in\xi \,\Longrightarrow\, \rank x\le \rank y,
$$
or equivalently,
$$
\rank y<\rank x \,\Longrightarrow\, xy\notin\xi.
$$
Furthermore, the total character of~$\xi$ allows to derive that
$$
\rank y<\rank x \,\Longrightarrow\, yx\in\xi.
$$
On the other hand, we know by Theorem~\ref{st:transSchulze} that $\nu=\mu(\isc)$ is transitive. As a consequence, by Lemma~\ref{st:rank}, $xy\in\nu$ implies $\rank x<\rank y$. By combining this with the preceding implication (with $x$ and $y$ interchanged with each other), we get that $\nu\sbseteq\xi$. In order to complete the proof that $\xi$ is admissible, we must check that $\xi\sbseteq\hat\nu$, or equivalently, that $xy\not\in\hat\nu$ implies $xy\not\in\xi$. This is true because of the following chain of implications:
$$
xy\not\in\hat\nu \,\Longleftrightarrow\,
yx\in\nu \,\Longrightarrow\,
\rank y < \rank x \,\Longrightarrow\,
xy\not\in\xi,
$$
where we used respectively the definition of $\hat\nu$, Lemma~\ref{st:rank}, and the hypothesis that $\rank x$ does not decrease along $\xi$. 
\end{proof}

\bigskip
In the following section we will make use of the following fact:
\begin{lemma}\hskip.5em
\label{st:rodrigues}
Given two admissible orders $\xi$ and $\xibis$, one can find a sequence of admissible orders $\xi_i\ (i=0\dots n)$ such that $\xi_0=\xi$, $\xi_n=\xibis$, and such that $\xi_{i+1}$ differs from $\xi_i$ only by the transposition of two consecutive elements.
\end{lemma}
\begin{proof}\hskip.5em
Given two total orders $\rho$ and $\sigma$, we will denote as $d(\rho,\sigma)$ the number of pairs $ab$ such that $ab\in\rho\setminus\sigma$. Obviously, $\rho=\sigma$ \ifoi $d(\rho,\sigma)=0$. Furthermore, we will say that $ab$ is a~\emph{consecutive} pair in  $\rho$ whenever $ab\in\rho$ and there is no $x\in\ist$ such that $ax,xb\in\rho$.
\ensep
If all pairs $ab$ which are consecutive in $\xi$ belong to $\xibis$, the~transitivity of $\xibis$ allows to derive that $\xi\sbseteq\xibis$; furthermore, the fact that all total orders on the finite set $\ist$ have the same number of pairs allows to conclude that $\xi=\xibis$.
\ensep
So, if~$\xibis\neq \xi$, there must be some pair $ab$ which is consecutive in $\xi$ but it does not belong to $\xibis$.
Since $ab$ belongs to the admissible order $\xi$ and $ba$ belongs to the admissible order $\xibis$, it follows that
$\img_{ab} = 0$.
\ensep
Let us take as $\xi_1$ the total order which differs from $\xi$ only by the transposition of the two consecutive elements $a$ and $b$; \ie $\xi_1 = (\xi\setminus\{ab\})\cup\{ba\}$. This order is admissible since $\xi$ is so and
$\img_{ab} = 0$.
Obviously, $d(\xi_1,\xibis) = d(\xi,\xibis) - 1$.
\ensep
From here, one can repeat the same process with $\xi_1$ substituted for~$\xi$: if $\xi_1$ still differs from $\xibis$ we take $\xi_2 = (\xi_1\setminus\{a_1b_1\})\cup\{b_1a_1\}$, where $a_1b_1$ is any pair which is consecutive in $\xi_1$ but it does not belong to $\xibis$, and so on. This iteration will conclude in a number of steps equal to $d(\xi,\xibis)$, since $d(\xi_i,\xibis)$ decreases by one unit in each step.
\end{proof}

\section{The projection}

Let us recall that our rating method is based upon certain projected scores~$\psc_{xy}$.
These quantities (or~equivalently, the projected margins $\pmg_{xy} = \psc_{xy} - \psc_{yx}$ and the projected turn\-overs~$\pto_{xy} = \psc_{xy} + \psc_{yx}$)
are worked out by means of the procedure {(\ref{eq:projection1}--\ref{eq:projection7})} of page~\pageref{eq:projection7}.
Its starting point are the indirect margins~$\img_{xy} = \isc_{xy} - \isc_{xy}$ and the original turnovers~$t_{xy} = v_{xy} + v_{yx}$.
From these quantities, equations~(\ref{eq:projection2}.1) and (\ref{eq:projection2}.2), used in this order, determine what we called the intermediate projected margins and turnovers, $\ppmg_{xy}$ and~$\ppto_{xy}$. After their construction, one becomes interested only in their superdiagonal elements~$\ppmg_{xx'}$ and~$\ppto_{xx'}$. In fact, these quantities are combined into certain intervals~$\gamma_{xx'}$ whose unions give rise to the whole set of projected scores.

\medskip
Let us recall in more detail the meaning of
the operator $\Psi$ which appears in step~(\ref{eq:projection2}.2).
This operator produces the intermediate projected turnovers $(\ppto_{xy})$ as a function of the original turnovers $(t_{xy})$
and the superdiagonal intermediate projected margins $(\ppmg_{pp'})$.
Here we are using parentheses to emphasize that we are dealing with the whole collection of turnovers and the whole collection of superdiagonal intermediate projected margins. Specifically, $(\ppto_{xy})$ is found by imposing certain conditions, namely (\ref{eq:bounded}--\ref{\incrementtwo}), and minimizing the function~(\ref{eq:phi}), which is nothing else than the euclidean distance to $(t_{xy})$. Equivalently, we can think in the following way (where the pair $xy$ is not restricted to belong to $\xi$): we consider a candidate $(\gto_{xy})$ which varies over the set~$T$ which is determined by the following conditions:
\begin{align}
&\hskip2em\gto_{yx} \,=\, \gto_{xy},
\label{eq:symmetry}
\\[2.5pt]
&\hskip.5em m{}^\sigma_{xx'} \,\le\, \gto_{xx'} \,\le\, 1;
\label{eq:taubounded}
\\[2.5pt]
&0 \,\le\, \gto_{xy} - \gto_{xy'} \,\le\, \ppmg_{yy'},
\label{eq:yincrement}
\end{align}
we associate each candidate $(\gto_{xy})$ with its euclidean distance from $(t_{xy})$;
finally, we define $(\ppto_{xy})$ as the only value of $(\gto_{xy})$ which minimizes such a distance. The minimizer exists and it is unique as a consequence of the fact that $T$ is a closed convex set
\cite[ch.\,I, \secpar{2}]{ks}.
In~this connection, one can say that $(\ppto_{xy})$ is the orthogonal projection of $(t_{xy})$ onto the convex set $T$.

\medskip
The procedure {(\ref{eq:projection1}--\ref{eq:projection7})} produces the projected scores as the end points of the intervals
\begin{equation}
\begin{repeatedNN}{eq:projection4}
\gamma_{xy} = \bigcup\,\,\{\gamma_{pp'}\mid x\rxieq p\rxi y\},
\end{repeatedNN}
\label{eq:projection4R}
\end{equation}
where 
\begin{equation}
\begin{repeatedNN}{eq:projection3}
\gamma_{xx'} = [\,(\ppto_{xx'}-\ppmg_{xx'})/2\,,\, (\ppto_{xx'}+\ppmg_{xx'})/2\,].
\end{repeatedNN}
\label{eq:projection3R}
\end{equation}
The desired properties of the projected scores and the associated margins and turnovers will be based upon the following properties of the intervals~$\gamma_{xy}$, where we recall that 
$|\gamma|$~means the length of an interval,
and $\centre\gamma$~means its barycentre, or centroid,
\ie the~number $(a+b)/2$ \,if $\gamma=[a,b]$.

\smallskip
\begin{lemma}\hskip.5em
\label{st:intervals}
The sets $\gamma_{xy}$ have the following properties
for~$x\rxi y\rxi z$:

\iim{a} $\gamma_{xy}$ is a closed interval.

\iim{b} $\gamma_{xy} \,\,\sbseteq\,\, [0,1]$.

\iim{c} $\gamma_{xz} \,\,=\,\, \gamma_{xy} \,\cup\, \gamma_{yz}$.

\iim{d} $\gamma_{xy} \,\cap\, \gamma_{yz} \,\,\neq \,\, \emptyset$.

\iim{e} $|\gamma_{xz}| \,\,\ge\,\, \max\,(\,|\gamma_{xy}|\,,\, |\gamma_{yz}|\,)$.

\iim{f} $\centre\gamma_{xy} \,\,\ge\,\, \centre\gamma_{xz} \,\,\ge \,\,\centre\gamma_{yz}$.

\iim{g} $|\gamma_{xz}|/\centre\gamma_{xz} \,\,\ge\,\, \max\,(\,|\gamma_{xy}|/\centre\gamma_{xy}\,,\, |\gamma_{yz}|/\centre\gamma_{yz}\,)$.
\end{lemma}
\begin{proof}\hskip.5em
Let us start by recalling that the superdiagonal intermediate turnovers and margins are ensured to satisfy the following inequalities:
\begin{gather}
\label{eq:bounded0}
0 \,\le\, \ppmg_{xx'} \,\le\, \ppto_{xx'} \,\le\, 1
\\[2.5pt]
\label{eq:overlapR}
\begin{repeated}{eq:overlap}
0 \,\le\, \ppto_{xx'} - \ppto_{x'x''} \,\le\, \ppmg_{xx'} + \ppmg_{x'x''}.
\end{repeated}
\end{gather}
From (\ref{eq:bounded0}) it follows that 
$0\le(\ppto_{xx'}-\ppmg_{xx'})/2\le(\ppto_{xx'}+\ppmg_{xx'})/2\le1$.
So, every $\gamma_{xx'}$ is an interval (possibly reduced to one point) and this interval is contained in $[0,1]$.
Also, the inequalities of~(\ref{eq:overlapR}) ensure on the one hand that $\centre\gamma_{xx'}\ge\centre\gamma_{x'x''}$, and on the other hand that the intervals $\gamma_{xx'}$ and $\gamma_{x'x''}$ overlap each other.
\ensep
In the following we will see that these facts about the elementary intervals $\gamma_{xx'}$ entail the stated properties of the sets $\gamma_{xy}$ defined by~(\ref{eq:overlapR}).

\newcommand\ppart[1]{\halfsmallskip\leavevmode
\hbox to4em{Part~(#1).\hss}}

\ppart{a} This is an obvious consequence of the fact that $\gamma_{pp'}$ and $\gamma_{p'p''}$ overlap each other.

\ppart{b} This follows from the fact that $\gamma_{pp'}\sbseteq[0,1]$.

\ppart{c} This is a consequence of the associative property enjoyed by the set-union operation.

\ppart{d} This is again an obvious consequence of the fact that $\gamma_{pp'}$ and $\gamma_{p'p''}$ overlap each other (take $p'=y$).

\ppart{e} This follows from~(c) because $\intg\sbseteq\inth$ implies $|\intg|\le|\inth|$.

\ppart{f} This follows from the fact that $\centre\gamma_{pp'}\ge\centre\gamma_{p'p''}$ because of the following general fact:
\ensep
If~$\intg$ and~$\inth$ are two intervals with $\centre\intg\ge\centre\inth$ then\linebreak 
$\centre\intg\ge \pcentre{\intg\cup\inth} \ge \centre\inth$. This is clear if $\intg$ and $\inth$ are disjoint and also if one of them is contained in the other. Otherwise, $\intg\setminus\inth$ and $\inth\setminus\intg$ are nonempty intervals and the preceding disjoint case
allows to proceed in the following way:
$$
\centre\intg
\,\ge\, \pcentre{\intg\cup (\inth\setminus\intg)}
\,=\, \pcentre{\intg\cup\inth}
\,=\, \pcentre{(\intg\setminus\inth)\cup\inth}
\,\ge\, \centre\inth.
$$ 

\ppart{g} This follows from (c) and (d) because of the following general fact: If~$\intg$ and~$\inth$ are two closed intervals with $\intg\sbseteq\inth\sbset[0,+\infty)$ then $|\intg|/\centre\intg \le |\inth|/\centre\inth$. In~fact, let $\intg=[a,b]$ and $\inth=[c,d]$. The hypothesis that $\intg\sbseteq\inth$ takes then the following form : $c\le a$ and $b\le d$. On the other hand, the claim that $|\intg|/\centre\intg \le |\inth|/\centre\inth$ takes the following form: $(b-a)/(b+a)\le (d-c)/(d+c)$. An~elementary computation shows that the latter is equivalent to $bc\le ad$, which is a consequence of the preceding inequalities.
\end{proof}

\medskip
The projection procedure makes use of a particular admissible order~$\xi$. In~fact, this order occurs in equations (\ref{eq:projection2}--\ref{eq:projection4}), as well as in conditions (\ref{eq:bounded}--\ref{\incrementtwo}).
In spite of this, the next theorem claims that
the final results are independent of~$\xi$.
The proof is not difficult, but it is rather long.

\smallskip
\begin{theorem} 
\label{st:independenceOfXi}
The projected scores do not depend on the admissible order~$\xi$ used for their calculation,
\ie the value of\, $\psc_{xy}$ is independent of\, $\xi$ \,for every $xy\in\tie$.
On the other hand, the matrix of the projected scores
in an admissible order~$\xi$
is also independent of\, $\xi$;
\ie if $x_i$ denotes the element of rank~$i$ in~$\xi$,
the value of $\psc_{x_ix_j}$ is independent of\, $\xi$ 
\,for every pair of indices~$i,j$.
\end{theorem}

\noindent
\textit{Remark}.\hskip.5em
The two statements say different things since the identity of $x_i$ and $x_j$ may depend on the admissible order~$\xi$.

\begin{proof}\hskip.5em
For the purposes of this proof it becomes necessary to change our set-up in a certain way. In fact, until now the intermediate objects $m^\sigma_{xy}$, $t^\sigma_{xy}$ and $\gamma_{xy}$ were considered only for~$x\rxi y$, \ie $xy\in\xi$. However, since we have to deal with changing the admissible order $\xi$, here we will allow their argument $xy$ to be any pair (of different elements), no matter whether it belongs to $\xi$ or not. In~this connection, we will certainly put $m^\sigma_{yx}=-m^\sigma_{xy}$ and $t^\sigma_{yx}= t^\sigma_{xy}$. On~the other hand, concerning $\gamma_{xy}$ and $\gamma_{yx}$, we will proceed in the following way: if $\gamma_{xy}=[a,b]$ then $\gamma_{yx}=[b,a]$. So,~generally speaking the $\gamma_{xy}$ are here ``oriented intervals'', \ie ordered pairs of real numbers. However, $\gamma_{xy}$ will always be ``positively oriented'' when $xy$ belongs to an admissible order (but it will be reduced to a point whenever there is another admissible order which includes~$yx$). In~particular, the $\gamma_{pp'}$ which are combined in (\ref{eq:projection4}) are always positively oriented intervals; so, the union operation performed in that equation can always be understood in the usual sense.
\ensep
In the following, $\gamma\rev$ denotes the oriented interval ``reverse'' to $\gamma$, \ie $\gamma\rev=[b,a]$ if $\gamma=[a,b]$.

So, let us consider the effect of replacing $\xi$ by another admissible order~$\xibis$. In the following, the tilde is systematically used to distinguish between hom\-olo\-gous objects which are associated respectively with $\xi$~and~$\xibis$; in particular, such a notation will be used in connection with the labels of the equations which are formulated in terms of
the assumed admissible order.

With this terminology, we will prove the two following equalities. First, 
\begin{equation}
\hskip.75em\gamma_{xy} \,=\, \gammabis_{xy},\hskip.75em\qquad
\hbox to60mm{for any pair\, $xy\,\ (x\neq y)$,\hfil}
\label{eq:ginvariance}
\end{equation}
where $\gamma_{xy}$ are the intervals produced by (\ref{eq:projection2}--\ref{eq:projection4}) together with the operation $\gamma_{yx}=\gpxy$, 
and $\gammabis_{xy}$ are those produced by (\tref{eq:projection2}--\tref{eq:projection4}) together with the operation $\gammabis_{yx}=\gbispxy$.
\ensep
Secondly, we will see also that
\begin{equation}
\gamma_{x_ix_j} \,=\, \gammabis_{\tilde x_i\tilde x_j},\qquad
\hbox to60mm{for any pair of indices\, $ij\,\ (i\neq j),$\hfil}
\label{eq:gijinvariance}
\end{equation}
where $x_i$ denotes the element of rank $i$ in $\xi$, and analogously for~$\tilde x_i$ in $\xibis$.
\ensep
These equalities contain the statements of the theorem since the projected scores are nothing else than the end points of the $\gamma$ intervals. 

\medskip
Now, by Lemma~\ref{st:rodrigues}, it suffices to deal with the case of two admissible orders $\xi$~and~$\xibis$ which differ from each other by one inversion only.
So, we will assume that there are two elements $a$ and $b$ such that the only difference between $\xi$~and $\xibis$ is that $\xi$ contains $ab$ whereas $\xibis$~contains~$ba$.
According to the definition of an admissible order, this implies that
$\img_{ab} = \img_{ba} = 0$.

\newcommand\Ant{P}
\newcommand\ant{p}
\newcommand\Pos{Q}
\newcommand\pos{q}

In order to control the effect of the differences between $\xi$ and $\xibis$, we will make use of the following notation:\ensep
$\Ant$ and $\ant$ will denote respectively the set of predecessors of~$a$ in $\xi$ \ensep and its lowest element, \ie  the immediate predecessor of~$a$ in $\xi$; in~this connection, any statement about $\ant$ will be understood to imply the assumption that $\Ant$ is not empty. Similarly, $\Pos$~and~$\pos$ will denote respectively the set of successors of~$b$ in $\xi$\ensep and its top element, \ie  the immediate successor of~$b$ in $\xi$; here too, any statement about $\pos$ will be understood to imply the assumption that $\Pos$ is not empty. So, $\xi$ and $\xibis$ contain respectively the paths $\ant ab\pos$ and $\ant ba\pos$.

\medskip
Let us look first at the superdiagonal intermediate projected margins~$\ppmg_{hh'}$.
According to (\ref{eq:projection2}.1), $\ppmg_{hh'}$ is the minimum of  a certain set
of values of $\img_{xy}$. In a table where $x$ and $y$ are ordered according to $\xi$, this set is an upper-right rectangle with lower-left vertex at $hh'$. Using \smash{$\xibis$} instead of $\xi$ amounts to interchanging two consecutive columns and the corresponding rows of that table, namely those labeled by $a$ and $b$. In spite of such a rearrangement, in all cases but one the underlying set from which the minimum is taken is exactly the same, so the mininum is the same. The only case where the underlying set is not the same occurs for~$h=a$ in the order $\xi$, or $h=b$ in the order $\xibis$; but then the minimum is still the same because the underlying set includes $\img_{ab} = \img_{ba} = 0$. So,
\begin{equation}
\ppmg_{x_ix_{i+1}} \,=\, \ppmgbis_{\tilde x_i\tilde x_{i+1}},\qquad \hbox{for any $i=1,2,\dots N\!-\!1$.}
\label{eq:mequivariance}
\end{equation}
In more specific terms, we have 
\begin{align}
\ppmg_{xx'}\! &\,=\, \ppmgbis_{xx'},\qquad \hbox{whenever \,$x\neq \ant,a,b$},
\label{eq:mxxp}
\\[2.5pt]
\ppmg_{\ant a} &\,=\, \ppmgbis_{\ant b},
\label{eq:mcacb}
\\[2.5pt]
\ppmg_{ab} &\,=\, \ppmgbis_{ba} \,=\, 0,
\label{eq:mabba}
\\[2.5pt]
\ppmg_{b\pos} &\,=\, \ppmgbis_{a\pos}.
\label{eq:mbdad}
\end{align}
In connection with equation~(\ref{eq:mxxp}) it should be clear that \emph{for~$x\neq \ant,a,b$ the immediate successor $x'$ is the same in both orders $\xi$ and $\xibis$}.

\medskip
Next we will see that the intermediate projected turnovers $\ppto_{xy}$ are invariant with respect to $\xi$:
\begin{equation}
\ppto_{xy} \,=\, \pptobis_{xy},\qquad\hbox{for any pair\, } xy\,\ (x\neq y),
\label{eq:tinvariance}
\end{equation}
where $\ppto_{xy}$ are the numbers produced by (\ref{eq:projection2}.2) together with the symmetry $\ppto_{yx}=\ppto_{xy}$, 
and $\pptobis_{xy}$ are those produced by (\tref{eq:projection2}.2) together with the symmetry $\pptobis_{yx}=\pptobis_{xy}$.

We will prove (\ref{eq:tinvariance}) by seeing that the set $T$ determined by conditions (\ref{eq:symmetry},\ref{eq:taubounded},\ref{eq:yincrement}) coincides exactly with the set $\widetilde T$ determined by (\ref{eq:symmetry},\tref{eq:taubounded},\tref{eq:yincrement}). In~other words, conditions (\ref{eq:taubounded}--\ref{eq:yincrement}) are exactly equivalent to (\tref{eq:taubounded}--\tref{eq:yincrement}) under condition (\ref{eq:symmetry}), which does not depend on~$\xi$.

In order to prove this equivalence we begin by noticing that condition~(\ref{eq:taubounded}) coincides exactly with (\tref{eq:taubounded}) when $x\neq \ant,a,b$. This is true because, on the one hand, $x'$ is then the same in both orders $\xi$ and $\xibis$, and, on the other hand, (\ref{eq:mxxp}) ensures that the right-hand sides have the same value. Similarly happens with conditions~(\ref{eq:yincrement}) and (\tref{eq:yincrement}) when $y\neq \ant,a,b$.
So, it remains to deal with conditions (\ref{eq:taubounded}) and (\tref{eq:taubounded}) for~$x=\ant,a,b$, and with conditions (\ref{eq:yincrement}) and (\tref{eq:yincrement}) for~$y=\ant,a,b$.\ensep
Now, on account of the symmetry~(\ref{eq:symmetry}), one easily sees that  condition (\ref{eq:taubounded}) with $x=a$ is equivalent to (\tref{eq:taubounded}) with $x=b$. In~fact, both of them reduce to $0 \le \gto_{ab} \le 1$ since $\ppmg_{ab} = \ppmgbis_{ba} = 0$, as it was   obtained in (\ref{eq:mabba}). This last equality ensures also the equivalence between condition~(\ref{eq:yincrement}) with~$y=a$ and condition~(\tref{eq:yincrement}) with~$y=b$. In this case both of them reduce to
\begin{equation}
\gto_{xa} \,=\, \gto_{xb}.
\label{eq:xaxb}
\end{equation}
This common equality plays a central role in
the equivalence between the remaining conditions.
\ensep
Thus, its combination with (\ref{eq:mbdad}) ensures the equivalence between  (\ref{eq:taubounded}) with $x=b$ and (\tref{eq:taubounded}) with $x=a$, as well as the equivalence between (\ref{eq:yincrement}) with $y=b$ and (\tref{eq:yincrement}) with $y=a$ when $x\neq a,b$. 
\ensep
On the other hand, its combination with (\ref{eq:mcacb}) ensures the equivalence between (\ref{eq:taubounded}) and (\tref{eq:taubounded}) when $x=\ant$, as well as the equivalence between (\ref{eq:yincrement}) and (\tref{eq:yincrement}) when $y=\ant$ and $x\neq a,b$. 
\ensep
Finally, we have the two following equivalences: (\ref{eq:yincrement})~with $y=\ant$ and $x=b$ is equivalent to (\tref{eq:yincrement}) with $y=\ant$ and $x=a$ because of the same equality  (\ref{eq:xaxb}) together with~(\ref{eq:mcacb}) and the symmetry~(\ref{eq:symmetry}); 
and similarly, (\ref{eq:yincrement}) with $y=b$ and $x=a$ is equivalent to (\tref{eq:yincrement}) with $y=a$ and $x=b$ because of (\ref{eq:xaxb}) together with~(\ref{eq:mbdad}) and~(\ref{eq:symmetry}).
\ensep
This completes the proof of~(\ref{eq:tinvariance}).

\medskip
Having seen that condition~(\ref{eq:xaxb}) is included in both~(\ref{eq:yincrement}) and~(\tref{eq:yincrement}), it follows that the intermediate projected turnovers satisfy
\begin{equation}
\ppto_{xa} \,=\, \ppto_{xb},\qquad \pptobis_{xa} \,=\, \pptobis_{xb}.
\end{equation}
By taking $x=\ant,\pos$ and using also~(\ref{eq:tinvariance}), it follows that
\begin{align}
\ppto_{xx'}\! &\,=\, \pptobis_{xx'},\qquad \hbox{whenever \,$x\neq \ant,a,b$},
\label{eq:txxp}
\\[2.5pt]
\ppto_{\ant a} &\,=\, \pptobis_{\ant b},
\label{eq:tcacb}
\\[2.5pt]
\ppto_{ab} &\,=\, \pptobis_{ba},
\label{eq:tabba}
\\[2.5pt]
\ppto_{b\pos} &\,=\, \pptobis_{a\pos}.
\label{eq:tbdad}
\end{align}
In other words, the superdiagonal intermediate turnovers satisfy
\begin{equation}
\ppto_{x_ix_{i+1}} \,=\, \pptobis_{\tilde x_i\tilde x_{i+1}},\qquad \hbox{for any $i=1,2,\dots N\!-\!1$.}
\label{eq:tequivariance}
\end{equation}
On account of the definition of $\gamma_{x_ix_{i+1}}$ and $\gammabis_{\tilde x_i\tilde x_{i+1}}$, 
the combination of (\ref{eq:mequivariance}) and (\ref{eq:tequivariance}) results in
\begin{equation}
\gamma_{x_ix_{i+1}} \,=\, \gammabis_{\tilde x_i\tilde x_{i+1}},\qquad \hbox{for any $i=1,2,\dots N\!-\!1$,}
\label{eq:gequivariance}
\end{equation}
from which the union operation~(\ref{eq:projection4}) produces~(\ref{eq:gijinvariance}).

\medskip
Finally, let us see that (\ref{eq:ginvariance}) holds too. To this effect, we begin by noticing that (\ref{eq:mabba}) together with (\ref{eq:tabba}) are saying not only that $\gamma_{ab} = \gammabis_{ba}$ but also that this interval reduces to a point. As a consequence, we have
\begin{equation}
\gamma_{ba} \,=\, \gamma_{ab} \,=\, \gammabis_{ba} \,=\, \gammabis_{ab}.
\label{eq:4gabba}
\end{equation}
Let us consider now the equation $\gamma_{\ant a} = \gammabis_{\ant b}$, which is contained in (\ref{eq:gequivariance}).
Since $\gamma_{ab}$ reduces to a point, the overlapping property $\gamma_{\ant a}\cap\gamma_{ab} \neq \emptyset$ (part~(d) of Lemma~\ref{st:intervals}) reduces to $\gamma_{ab}\sbseteq\gamma_{\ant a}$. Therefore, $\gamma_{\ant b} = \gamma_{\ant a}\cup\gamma_{ab} = \gamma_{\ant a}$ (where we used part~(c) of Lemma~\ref{st:intervals}). 
Analogously, $\gammabis_{pa} = \gammabis_{\ant b}\cup\gammabis_{ba} = \gammabis_{\ant b}$. Altogether, this gives
\begin{equation}
\gamma_{\ant b} \,=\, \gamma_{\ant a} \,=\, \gammabis_{\ant b} \,=\, \gammabis_{\ant a}.
\label{eq:4gcacb}
\end{equation}
By means of an analogous argument, one obtains also that
\begin{equation}
\gamma_{a\pos} \,=\, \gamma_{b\pos} \,=\, \gammabis_{a\pos} \,=\, \gammabis_{b\pos}.
\label{eq:4gbdad}
\end{equation}
On the other hand, (\ref{eq:gequivariance}) ensures that
\begin{equation}
\gamma_{xx'} \,=\, \gammabis_{xx'},\qquad \hbox{whenever \,$x\neq \ant,a,b$}.
\label{eq:4gxxp}
\end{equation}
Finally, part~(c) of Lemma~\ref{st:intervals} allows to go from (\ref{eq:4gabba}--\ref{eq:4gxxp}) to the desired general equality~(\ref{eq:ginvariance}).
\end{proof}

\smallskip
\begin{theorem}\hskip.5em
\label{st:propertiesOfProjection}
The projected scores and their asssociated margins and turn\-overs satisfy the following properties with respect to any admissible order~$\xi$:

\halfsmallskip\noindent
\textup{(a)}~The following inequalities hold whenever $x\rxi y$:
\begin{alignat}{2}
\label{eq:vxyinequality}
\psc_{xy} \,&\ge\, \psc_{yx}
&\qquad
\pmg_{xy} \,&\ge\, 0,
\\[2.5pt]
\label{eq:vinequalities}
\psc_{xz} \,&\ge\, \psc_{yz},
&\qquad
\psc_{zx} \,&\le\, \psc_{zy},
\\[2.5pt]
\label{eq:minequalities}
\pmg_{xz} \,&\ge\, \pmg_{yz},
&\qquad
\pmg_{zx} \,&\le\, \pmg_{zy},
\\[2.5pt]
\label{eq:tinequalities}
\pto_{xz} \,&\ge\, \pto_{yz},
&\qquad
\pto_{zx} \,&\ge\, \pto_{zy},
\\[2.5pt]
\label{eq:mtinequalities}
\pmg_{xz}/\pto_{xz} \,&\ge\, \pmg_{yz}/\pto_{yz},
&\qquad
\pmg_{zx}/\pto_{zx} \,&\le\, \pmg_{zy}/\pto_{zy}.
\end{alignat}

\halfsmallskip\noindent
\textup{(b)} If~$\psc_{xy}=\psc_{yx}$, or equivalently~$\pmg_{xy}=0$,  then \textup{(\ref{eq:vinequalities}--\ref{eq:mtinequalities})} are satisfied all of them with an equality sign.

\halfsmallskip\noindent
\textup{(c)}~In the complete case, the projected margins satisfy the following property:
\begin{equation}
\label{eq:equaltomaxBis}
\pmg_{xz} \,=\, \max\,(\pmg_{xy},\pmg_{yz}),\qquad \hbox{whenever $x\rxi y\rxi z$}.
\end{equation}
\end{theorem}

\begin{proof}\hskip.5em
We will see that these properties derive from those satisfied by the $\gamma$ intervals, which are collected in Lemma~\ref{st:intervals}.
For the derivation one has to bear in mind that $\psc_{xy}$ and $\psc_{yx}$ are respectively the right and left end points of the interval $\gamma_{xy}$, and that $\pmg_{xy} = -\pmg_{yx}$ and $\pto_{xy} = \pto_{yx}$ are respectively the width and twice the barycentre of $\gamma_{xy}$.

\smallskip
Part~(a).\ensep
Let us begin by noticing that (\ref{eq:minequalities}) will be an immediate consequence of (\ref{eq:vinequalities}), since $\pmg_{xz} = \psc_{xz} - \psc_{zx}$ and $\pmg_{yz} = \psc_{yz} - \psc_{zy}$. 
On the other hand, (\ref{eq:tinequalities}.2) is equivalent to (\ref{eq:tinequalities}.1) \,and\, (\ref{eq:mtinequalities}.2) is equivalent to (\ref{eq:mtinequalities}.1). These equivalences hold because the turnovers and margins are respectively symmetric and antisymmetric. Now, (\ref{eq:vxyinequality}) holds as soon as $\gamma_{xy}$ is an interval, as it is ensured by part~(a) of Lemma~\ref{st:intervals}. So, it remains to prove the inequalities (\ref{eq:vinequalities}), (\ref{eq:tinequalities}.1) and (\ref{eq:mtinequalities}.1). In order to prove them we will distinguish three cases, namely:\ensep
(i)~$x\rxi y\rxi z$;\ensep
(ii)~$z\rxi x\rxi y$;\ensep
(iii)~$x\rxi z\rxi y$.\ensep

Case~(i)\,: By part~(c) of Lemma~\ref{st:intervals}, in this case we have $\gamma_{xz}\spseteq\gamma_{yz}$. This immediately implies (\ref{eq:vinequalities}) because $[a,b]\spseteq[c,d\,]$ is equivalent to saying that $b\ge d$ and $a\le c$. On the other hand, the inequalities~(\ref{eq:tinequalities}.1) and (\ref{eq:mtinequalities}.1) are contained in parts (f) and (g) of Lemma~\ref{st:intervals}.
Case~(ii)\, is analogous to case~(i).

Case~(iii)\,: In this case, (\ref{eq:vinequalities}) follows from part~(d) of Lemma~\ref{st:intervals} since $[a,b]\cap[c,d\hskip2pt]\neq \emptyset$ is equivalent to saying that $b\ge c$ and $a\le d$. On the other hand,  (\ref{eq:tinequalities}.1) is still contained in part~(f) of Lemma~\ref{st:intervals} (because of the symmetric character of the turnovers), and (\ref{eq:mtinequalities}.1) holds since $\pmg_{xz} \ge 0 \ge \pmg_{yz}$.

\smallskip
Part~(b).\ensep
The hypothesis that $\psc_{xy}=\psc_{yx}$ is equivalent to saying that $\gamma_{xy}$ reduces to a point, \ie $\gamma_{xy}=[v,v]$ for some $v$. The claimed equalities will be obtained by showing that in these circumstances one has $\gamma_{xz} = \gamma_{yz}$. 
We will distinguish the same three cases as in part~(a).

Case~(i)\,: On account of the overlapping property $\gamma_{xy}\cap\gamma_{yz}\neq \emptyset$ (part~(d) of Lemma~\ref{st:intervals}), the one-point interval $\gamma_{xy}=[v,v]$ must be contained in $\gamma_{yz}$. So,~$\gamma_{xz} = \gamma_{xy}\cup\gamma_{yz} = \gamma_{yz}$ (where we used part~(c) of Lemma~\ref{st:intervals}).
Case~(ii)\, is again analogous to case~(i).

Case~(iii)\,: By part~(c) of Lemma~\ref{st:intervals} (with $y$ and $z$ interchanged with each other), the fact that $\gamma_{xy}$ reduces to the one-point interval $[v,v]$ implies that  both $\gamma_{xz}$ and $\gamma_{zy}$ reduce also to this one-point interval

\smallskip
Part~(c).\ensep
In the complete case the intermediate projected turnovers are all of them equal to $1$, so the intervals $\gamma_{pp'}$ and $\gamma_{xy}$ are all of the centred at~$1/2$. In these circumstances, (\ref{eq:equaltomaxBis}) is exactly equivalent to part~(c) of Lemma~\ref{st:intervals}.
\end{proof}

\medskip
The following propositions identify certain situations where the preceding projection reduces to the identity.

\smallskip
\begin{proposition}\hskip.5em
\label{st:plumpps}
In the case of plumping votes the projected scores coincide with the original ones.
\end{proposition}

\begin{proof}\hskip.5em
Let us begin by recalling that in the case of plumping votes the binary scores have the form $v_{xy} = \plumpf_x$ for~every $y\neq  x$,
where $\plumpf_x$ is the fraction of voters who choose~$x$.
\ensep
This implies that $\isc_{xy} = v_{xy} = \plumpf_x$. In fact, any path $\gamma$ from $x$ to $y$ starts with a link of the form $xp$, whose associated score is $v_{xp} = \plumpf_x$. So $v_\gamma \le \plumpf_x$ and therefore $\isc_{xy} \le \plumpf_x$. But on the other hand $\plumpf_x = v_{xy} \le \isc_{xy}$.
\ensep
Consequently, we get $\img_{xy} = \isc_{xy} - \isc_{yx} = v_{xy} - v_{yx} = \plumpf_x - \plumpf_y$, and the admissible orders are those for which the $\plumpf_x$ are non-increasing.
\ensep
Owing to this non-increasing character, the intermediate projected margins are~$\ppmg_{xx'} = m_{xx'} = \plumpf_x - \plumpf_{x'}$.
\ensep
On the other hand, the intermediate projected turnovers are~$\ppto_{xy} = t_{xy} = \plumpf_x + \plumpf_{y}$. In fact these numbers are easily seen to satisfy conditions (\ref{eq:bounded}--\ref{eq:qincrement}) and they obviously minimize~(\ref{eq:phi}).
\ensep
As a consequence, $\gamma_{xx'} = [\plumpf_{x'},\plumpf_x]$. In particular, the intervals $\gamma_{xx'}$ and $\gamma_{x'x''}$ are adjacent to each other (the right end of the latter coincides with the left end of the former). This fact entails that $\gamma_{xy} = [\plumpf_y,\plumpf_x]$ whenever $x\rxi y$.
\ensep
Finally, the projected scores are the end points of these intervals, namely $\psc_{xy} = \plumpf_x = v_{xy}$ and $\psc_{yx} = \plumpf_y = v_{yx}$.
\end{proof}

\smallskip
\begin{proposition}\hskip.5em
\label{st:condrg}
Assume that the votes are complete.
Assume also that there exists a total order $\xi$ such that $\mu(v) \sbseteq \xi \sbseteq \hat\mu(v)$ and such that the original margins satisfy
\begin{equation}
m_{xz} \,=\, \max\,(m_{xy},m_{yz}),\qquad \hbox{whenever $x\rxi y\rxi z$ in $\xi$}.
\label{eq:mequaltomax}
\end{equation}
In that case, the projected scores coincide with the original ones. Besides, condition $(\ref{eq:mequaltomax})$ holds also for any other total order $\xibis$ which satisfies $\mu(v) \sbseteq \xibis \sbseteq \hat\mu(v)$.
\end{proposition}

\noindent
\textit{Remark}.\hskip.5em
The hypothesis that $\mu(v)\sbseteq\xi\sbseteq\hat\mu(v)$ is not the one which defines an admissible order, namely $\mu(\isc)\sbseteq\xi\sbseteq\hat\mu(\isc)$. However, in the course of the proof we will see that $\isc_{xy} = v_{xy}$. So, $\xi$ will be after all an admissible order.

\begin{proof}\hskip.5em
Since we are in the complete case, the scores $v_{xy}$ and the margins $m_{xy}$ are related to each other by the monotone increasing transformation $v_{xy}=(1+m_{xy})/2.$ Therefore, condition~(\ref{eq:mequaltomax}) on the margins is equivalent to the following one on the scores:
\begin{equation}
v_{xz} \,=\, \max\,(v_{xy},v_{yz}),\qquad \hbox{whenever $x\rxi y\rxi z$ in $\xi$}.
\label{eq:vequaltomax}
\end{equation}
On the other hand, since $v_{xy} + v_{yx} = 1$, the preceding condition is also equivalent to the following one:
\begin{equation}
v_{zx} \,=\, \min\,(v_{zy},v_{yx}),\qquad \hbox{whenever $x\rxi y\rxi z$ in $\xi$}.
\label{eq:vequaltomin}
\end{equation}
In fact, $v_{zx}=1-v_{xz}=1-\max\,(v_{xy},v_{yz})=\min\,(v_{zy},v_{yx})$.

\halfsmallskip
Now, we claim that these properties imply the following one:
\begin{equation}
v_{xz} \ge \min(v_{xy},v_{yz}),\qquad
\hbox{for any $x,y,z$.}
\label{eq:vgreaterthanmin}
\end{equation}
In order to prove~(\ref{eq:vgreaterthanmin}) we will distinguish four cases depending on whether or not do $xy$ and $yz$ belong to $\xi$:
\ensep
(a)~If $xy,yz\in\xi$, then (\ref{eq:vgreaterthanmin}) is an immediate consequence of (\ref{eq:vequaltomax}).
\ensep
(b)~Similarly, if $xy,yz\not\in\xi$, then (\ref{eq:vgreaterthanmin}) is an immediate consequence of (\ref{eq:vequaltomin}) with $x$ and $z$ interchanged with each other.
\ensep
(c)~Consider now the case where $xy\not\in\xi$ and $yz\in\xi$. In~this case we have $v_{xy}\le 1/2\le v_{yz}$, so $\min(v_{xy},v_{yz})=v_{xy}$. Now we must distinguish two subcases: If $xz\in\xi$, then $v_{xy}\le 1/2\le v_{xz}$, so we get (\ref{eq:vgreaterthanmin}). If, on the contrary, $zx\in\xi$, then (\ref{eq:vequaltomin}) applied to $y\rxi z\rxi x$ gives $v_{xy} = \min\,(v_{xz},v_{zy}) \le v_{xz}$ as claimed.
\ensep
(d)~Finally, the case where $xy\in\xi$ and $yz\not\in\xi$ is analogous to the preceding one.

\halfsmallskip
Now we invoke Lemma~\ref{st:indirectEqualToDirect}, according to which (\ref{eq:vgreaterthanmin}) implies that $\isc_{xy}=v_{xy}$.
\ensep
In particular, $\xi$ is ensured to be an admissible order.
Let us consider any pair $xy$ contained in $\xi$. By applying condition (\ref{eq:mequaltomax}) we see that $\ppmg_{xy} = \img_{xy} = m_{xy}$. On the other hand, since the votes are complete we have $\ppto_{xy}=1$. So,~the intervals $\gamma_{pp'}$ and their unions are all of them centred at $1/2$. In this case, the union operation of (\ref{eq:projection4}) is equivalent to a maximum operation performed upon the margins. On account of (\ref{eq:mequaltomax}), this implies that $\pmg_{xy}=m_{xy}$. Since we also have $\pto_{xy} = 1 = t_{xy}$, it follows that $\psc_{xy}=v_{xy}$ and $\psc_{yx}=v_{yx}$.

\halfsmallskip
Having proved that $v_{xy}=\isc_{xy}=\psc_{xy}$, and taking into account that this entails $m_{xy}=\pmg_{xy}$, one easily sees that condition~(\ref{eq:mequaltomax}) holds also for any other total order $\xibis$ such that $\mu(v)\subseteq\xibis\subseteq\hat\mu(v)$. In fact, such an order is an admissible one, since $\mu(\isc) = \mu(v)$, and that condition is guaranteed by part~(c) of Theorem~\ref{st:propertiesOfProjection}. 
\end{proof}

\section{The rank-like rates}

Let us recall that the rank-like rates $r_x$ are given by the formula
\begin{equation}
\label{eq:rratesR}
\begin{repeated}{eq:rrates}
r_x \,=\, N - \sum_{y\neq x}\,\psc_{xy}.
\end{repeated}
\end{equation}
where $\psc_{xy}$ are the projected scores. In the special case of complete votes, where $\psc_{xy}+\psc_{yx}=1$, the preceding formula is equivalent to the following one:\linebreak
\begin{equation}
\label{eq:rratesfrommarginsR}
\begin{repeated}{eq:rratesfrommargins}
r_x \,=\, (N+1 - \sum_{y\neq x}\,\pmg_{xy}\,) \,/\, 2.
\end{repeated}
\end{equation}
Let us remark also that in this special case the rank-like rates have the property that
\begin{equation}
\label{eq:sumrrates}
\sum_{x\in\ist}r_x \,=\, N(N+1)/2.
\end{equation}

\pagebreak 

\bigskip
In view of formula~(\ref{eq:rratesR}),
the properties of the projected scores obtained in Theorem~\ref{st:propertiesOfProjection} imply the following facts:

\smallskip
\begin{lemma}\hskip.5em
\label{st:obsRosa}

\iim{a} If $x\rxi y$ in an admissible order $\xi$,
then $r_x\le r_y$.

\iim{b} $r_x=r_y$ \ifoi $\psc_{xy}=\psc_{yx}$, \ie $\pmg_{xy}=0$.

\iim{c} The inequalities $(\ref{eq:vxyinequality}\text{--\,}\ref{eq:mtinequalities})$ are satisfied whenever $r_x\le r_y$.\hfil\break
In particular, $\psc_{xy}>\psc_{yx}$ implies $r_x<r_y$.
\end{lemma}
\begin{proof}\hskip.5em
Part~(a).\ensep
It is an immediate consequence of the preceding formula together with the inequalities (\ref{eq:vxyinequality}) and (\ref{eq:vinequalities}.1) ensured by Theorem~\ref{st:propertiesOfProjection}.

\halfsmallskip
Part~(b).\ensep
According to the formula above,
\begin{equation}
\label{eq:riguals}
r_y - r_x \,=\, (\psc_{xy} - \psc_{yx}) \,+\,
\sum_{\latop{\scriptstyle z\neq x}{\scriptstyle z\neq y}}%
\,(\psc_{xz}-\psc_{yz}).
\end{equation}
Let $\xi$ be an admissible order. By symmetry we can assume $xy\in\xi$. As~a consequence, Theorem~\ref{st:propertiesOfProjection} ensures that the terms of (\ref{eq:riguals}) which appear in parentheses are all of them greater than or equal to zero. So the only possibility for their sum to vanish is that each of them vanishes separately, \ie $\psc_{xy}=\psc_{yx}$ and $\psc_{xz}-\psc_{yz}$ for any $z\not\in\{x,y\}$. Finally, part~(b) of Theorem~\ref{st:propertiesOfProjection}  ensures that all of these equalities hold as soon as the first one is satisfied.

\halfsmallskip
Part~(c).\ensep
It suffices to use the contrapositive of~(a) in the case of a strict inequality and (b)~together with Theorem~\ref{st:propertiesOfProjection}.(b) in the case of an equality.
\end{proof}

\smallskip
\begin{theorem}\hskip.5em
\label{st:RvsNu}
The rank-like rating given by~\textup{(\ref{eq:rratesR})}
is related to the indirect comparison relation~$\nu = \mu(\isc)$ in the following way:

\iim{a}$xy\in\hat\nu \,\Rightarrow\, r_x \le r_y$.

\iim{b}$r_x < r_y \,\Rightarrow\, xy\in\nu$.

\iim{c}If $\nu$ contains a set of the form $\xst\times\yst$ with $\xst\cup\yst=\ist$,\hfil\break 
then $r_x < r_y$ for any $x\in\xst$ and $y\in\yst$.

\iim{d}If $\nu$ is total, \ie $\hat\nu=\nu$, then $xy\in\nu \,\Leftrightarrow\, r_x < r_y$.
\end{theorem}

\begin{comment}
More generally than (d): Assume that $\nu$ has the following properties: (i)~$xy\in\nu$ plus $yz\in\hat\nu$ imply $xz\in\nu$; (ii)~$xy\in\hat\nu$ plus $yz\in\hat\nu$ plus $xz\in\nu$ imply $xy\in\nu$ or $yz\in\nu$. The $xy\in\nu \,\Leftrightarrow\, r_x < r_y$.

A particular case of the latter: There exists $\lambda:\ist\rightarrow\bbr$ such that $\nu=\{xy\mid \lambda_x>\lambda_y\}$

Proof: Generalize that of Theorem~\ref{st:av4}.
\end{comment}

\begin{proof}\hskip.5em
Part~(a).\ensep 
Let us begin by noticing that $xy\in\nu$ implies $r_x \le r_y$.
This follows from part~(a) of Lemma~\ref{st:obsRosa} since $\nu$ is included in any admissible ordering $\xi$.
\ensep
Consider now the case $xy\in\hat\nu\setminus\nu$. This is equivalent to saying that $\nu$ contains neither $xy$ nor $yx$. Now, in this case Theorem~\ref{st:existenceXiThm} ensures the existence of an admissible order which contains such a pair $xy$. So, using again the preceding proposition, we are still ensured that $r_x\le r_y$.

\halfsmallskip
Part~(b).\ensep 
It reduces to the the contrapositive of~(a).

\halfsmallskip
Part~(c).\ensep 
Let $x\in\xst$ and $y\in\yst$. Since $\xst\times\yst \sbset \nu$, part~(a) ensures that $r_x \le r_y$. So, it suffices to exclude the possibility that $r_x = r_y$. This will be done by showing that this equality leads to a contradiction. By part~(b) of Lemma~\ref{st:obsRosa}, that equality implies $\psc_{xy} = \psc_{yx}$, or equivalently, $\pmg_{xy} = 0$. But according to (\ref{eq:projection3}--\ref{eq:projection5}), this means that $\ppmg_{hh'}=0$ for all $h$ such that $x\rxieq h\rxi y$.
Here we are making use of an admissible order~$\xi$.
In particular we have $\ppmg_{\last\last'}=0$,
where $\last$ denotes the lowest element of $\xst$ according to~$\xi$, 
and $\last'$ is the top element of $\yst$. But this contradicts the fact that $\last\last'\in\xst\times\yst\sbset\nu$.

\halfsmallskip
Part~(d).\ensep
It suffices to show that $r_x < r_{x'}$, where $x'$ denotes the item that immediately follows $x$ in the total order $\nu$. This follows from part~(c) by taking $\xst = \{p\mid p\rxieq x\}$ and $\yst = \{q\mid x'\rxieq q\}$ and using the transitivity of~$\nu$.
\end{proof}

\medskip
By construction, 
the rank-like rates are related to the projected scores
in the same way as
the average ranks are related to the original scores
when the votes are complete rankings~(\secpar{2.5}).
Therefore, if we are in the case of complete ranking votes
and the projected scores coincide with the original ones,
then the rank-like rates coincide with the average ranks:

\smallskip
\begin{proposition}\hskip.5em
\label{st:avranks}
Assume that the votes are complete rankings. Assume also that the Llull matrix satisfies the hypothesis of Proposition~\ref{st:condrg}.
In that case, the rank-like rates~$r_x$ coincide exactly with the average ranks~$\bar r_x$.
\end{proposition}

\begin{proof}\hskip.5em
This is an immediate consequence of Proposition~\ref{st:condrg}.
\end{proof}

\section{Zermelo's method}

The Llull matrix of a vote can be viewed as corresponding to a tournament between the members of~$\ist$ where $x$ and $y$ have played $T_{xy}$ matches (the number of voters who made a comparison between $x$ and $y$, even if this comparison resulted in a tie) and $V_{xy}$ of these matches were won by $x$, whereas the other $V_{yx}$ were won by~$y$ (one tied match will be counted as half a match in favour of $x$ plus half a match in favour of $y$).
For such a scenario, Ernst Zermelo~\cite{ze} devised in 1929 a rating method which turns out to be quite suitable to convert our rank-like rates into fraction-like ones. This method was rediscovered later on by other autors~\cite{bt,fo}.

\medskip
Zermelo's method is based upon a probabilistic model for the outcome of a match between two items $x$ and $y$. This model assumes that such a match is won by $x$ with probability $\flr_x/(\flr_x+\flr_y)$ whereas it is won by $y$ with probability $\flr_y/(\flr_x+\flr_y)$, where $\flr_x$ is a non-negative parameter associated with each player~$x$, usually referred to as its {\df strength}. If~all matches are independent events, the probability of obtaining a particular system of values for the scores $(V_{xy})$ 
is~given by
\begin{equation}
P \,=\, \prod_{\{x,y\}}\,\left(\vbox{\halign{\hfil#\hfil\cr$T_{xy}$\cr$V_{xy}$\cr\noalign{\vskip-11pt}}}\right)
\left(\frac{\flr_x}{\flr_x+\flr_y}\right)^{\hskip-3ptV_{xy}}
\left(\frac{\flr_y}{\flr_x+\flr_y}\right)^{\hskip-3ptV_{yx}},
\label{eq:probabilitat}
\end{equation}
where the product runs through all unordered pairs $\{x,y\}\sbseteq\ist$ with $x\neq y$. Notice that $P$~depends only on the strength ratios; in other words, multiplying all the strengths by the same value has no effect on the result. On account of this, we will normalize the strengths by requiring their sum to take a fixed positive value~$f$.
\ensep
In order to include certain extreme cases, one must allow for some of the strengths to vanish. However, this may conflict with $P$ being well defined, since it could lead to indeterminacies of the type $0/0$ or $0^0$. So,~one should be careful in connection with vanishing strengths.
\ensep
With all this in mind, for the moment we will let the strengths vary in the following set:
\begin{equation}
\phiset = \{\,\flr\in\bbr^\ist\mid\flr_x>0\text{ \,for all }x\in\ist,\ \sum_{x\in\ist}\flr_x = \nev\,\}.
\label{eq:conjuntphis}
\end{equation}
Together with this set, in the following we will consider also its closure~$\phisetb$, which includes vanishing strengths, and its boundary $\partial\phiset = \phisetb\setminus\phiset$.

\medskip
In connection with our interests,
it is worth noticing that Zermelo's model can be viewed as a special case of a more general one, proposed in 1959 by~Robert Duncan Luce, which considers the outcome of making a choice out~of~multiple options \cite{luce}.
According to Luce's ‘choice axiom’, the~probabilities of two different choices $x$ and $y$ \,are in a ratio which does not depend on~which other options are present.
As a consequence, it follows that 
every option $x$ can be associated a number~$\flr_x$ so that the probability of~choosing $x$ out of a set~$\xst\ni x$ is given by $\flr_x/(\sum_{y\in\xst}\flr_y)$. Obviously, Zermelo's model corresponds to considering binary choices only.
\ensep
It~is interesting to notice that 
Luce's model allows to associate every ranking with a certain probability. In~fact, a ranking can be viewed as the result of first choosing the winner out of the whole set $\ist$, then choosing the best of the remainder, and so on. If~these successive choices are assumed to be independent events, then one can easily figure out the corresponding probability.
\ensep
\begin{comment}
The expected rank of $x$ is~$N - \sum_{y\ne x}\flr_x/(\flr_x+\flr_y)$.
\end{comment}
\noindent
Anyway, when the normalization condition $\sum_{x\in\ist}\flr_x = \nev$ ($\le1$) is adopted,
Luce's theory of choice allows
to view $\flr_x$ as the first-choice probability of~$x$,
and to view $1-\nev$ as the probability of abstaining from making a choice out~of~$\ist$.

Let us mention here also that the hypothesis of independence which lies behind formula~(\ref{eq:probabilitat}) is certainly not satisfied by the binary comparisons which arise out of preferential voting.
In order to satisfy that hypothesis, the individual votes should be  based upon independent binary comparisons, in which case they could take the form of an arbitrary binary relation, as we considered in~\secpar{3.3}. However, even if the independence hypothesis is not satisfied, we will see that Zermelo's
method, which we are about to discuss, 
has good properties for 
transforming our projected scores into fraction-like rates.

\medskip
Zermelo's method corresponds to a maximum likelihood estimate of the parameters~$\flr_{x}$ from a given set of actual values of~$V_{xy}$ (and of $T_{xy}=V_{xy}+V_{yx}$). In other words, given the values of~$V_{xy}$, one looks for the values of $\flr_x$ which maximize the probability~$P$. Since $V_{xy}$ and $T_{xy}$ are now fixed, this is equivalent to maximizing the following function of the~$\flr_x$:
\begin{equation}
F(\flr) \,=\, \prod_{\{x,y\}}\,
\frac{{\flr_x}^{v_{xy}}\,{\flr_y}^{v_{yx}}}{(\flr_x+\flr_y)^{t_{xy}}},
\label{eq:funcioF}
\end{equation}
(recall that $v_{xy} = V_{xy}/V$ and $t_{xy} = T_{xy}/V$ where $V$~is a positive constant greater than or equal to any of the turnovers $T_{xy}$; going from (\ref{eq:probabilitat}) to (\ref{eq:funcioF}) involves taking the power of exponent $1/V$ and disregarding a fixed multiplicative constant).
\ensep
The function $F$ is certainly smooth on $\phiset$. Besides, it~is clearly bounded from above, since the probability is always less than or equal to~$1$. However, generally speaking $F$~needs not to achieve a maximum in $\phiset$, because this set is not compact. In~the present situation, the~only general fact that one can guarantee in this connection is the existence of maximizing sequences, \ie sequences $\flrn$ in~$\phiset$ with the property that $F(\flrn)$ converges to the lowest upper bound $\lub=\sup\,\{F(\psi)\,|\,\psi\in\phiset\}$.

\medskip
In connection with maximizing the function $F$ defined by (\ref{eq:funcioF}) it makes a difference whether two particular items $x$ and $y$ satisfy or not the inequality $v_{xy}>0$, or more generally ---as we will see--- whether they satisfy\linebreak[4] $\isc_{xy}>0$. By the definition of $\isc_{xy}$, the~last inequality defines a transitive relation \hbox{---namely} the transitive closure of the one defined by the former inequality---. In the following we will denote this transitive relation by the symbol \,$\chrel$. Thus, 
\begin{equation}
x\chrel y \ \Longleftrightarrow\  \isc_{xy}>0.
\end{equation}
Associated with it, it is interesting to consider also the following derived relations,
which keep the property of transitivity and are respectively symmetric and antisymmetric:
\begin{align}
x\chrels y \ &\Longleftrightarrow\  \isc_{xy}>0 \text{ \,and\, } \isc_{yx}>0,\\[2.5pt]
x\chrela y \ &\Longleftrightarrow\  \isc_{xy}>0 \text{ \,and\, } \isc_{yx}=0.
\end{align}
Therefore, $\chrels$ is an equivalence relation
and $\chrela$ is a partial order.
In the following, the situation where $x\chrela y$ will be expressed by saying that $x$ \dfc{dominates} $y$.
The equivalence classes of~$\ist$ by $\chrels$ are called the \dfc{irreducible components} of~$\ist$ (for~$\vaa$). If there is only one of them, namely $\ist$ itself, then one says that the matrix $\vaa$ is irreducible. So, $\vaa$ is irreducible \ifoi $\isc_{xy}>0$ for any $x,y\in\ist$. It is 
not difficult to see that this property is equivalent to the following one formulated in terms of the direct scores only: there is no splitting of $\ist$
into two classes $\xst$ and~$\yst$ so that $v_{yx}=0$ for any $x\in\xst$ and $y\in\yst$; in other words, there is no ordering of $\ist$ for which the matrix~$\vaa$ takes the form
\begin{equation}\label{eq:blocs}
 \begin{pmatrix} \vxx & \vxy \\ \zeromatrix & \vyy \end{pmatrix},
\end{equation}
where $\vxx$ and $\vyy$ are square matrices and $\zeromatrix$ is a  zero matrix.
Besides, a subset $\xst\sbseteq\ist$ is an irreducible component \ifoi  $\xst$ is maximal, in~the sense of set inclusion, for the property of $\vxx$ being irreducible.
\ensep
On~the other hand, 
it~also happens that the relations \,$\chrel$\, and \,$\chrela$\, are compatible with the equivalence relation \,$\chrels$, \ie if $x\chrels\bar x$ and $y\chrels\bar y$ then \,$x\chrel y$ implies $\bar x\chrel\bar y$,\, and analogously $x\chrela y$ implies $\bar x\chrela\bar y$.
As a consequence, the relations \,$\chrel$\, and \,$\chrela$\, can be 
applied also to the irreducible components of $\ist$ for $\vaa$.
In the following we will be interested in the case where $\vaa$ is irreducible, or more generally, when there is a \dfc{top dominant irreducible component}, \ie an irreducible component which dominates any other.
\ensep
From now on we systematically use the notation $\vrs$ to mean the restriction of $(v_{xy})$ to~$x\in R$ and $y\in S$, where $R$ and $S$ are arbitrary non-empty subsets of $\ist$.
Similarly, $\flrr$~will denote the restriction of $(\flr_x)$ to~$x\in R$.


\medskip
The next theorems collect the basic results that we need about Zermelo's method.

\begin{theorem}[Zermelo, 1929 \cite{ze}; see also \cite{fo,ke}]\hskip.5em
\label{st:zermeloirre}
If $\vaa$ is irreducible,\, then:

\iim{a}There is a unique $\flr\in\phiset$ which maximizes $F$ on $\phiset$.

\iim{b}$\flr$ is the solution of the following system of equations:
\begin{align}
\sum_{y\neq x}\,t_{xy}\,\frac{\flr_x}{\flr_x+\flr_y} \,&=\, \sum_{y\neq x}\,v_{xy},
\label{eq:fratesZ}
\\[2.5pt]
\sum_x\,\flr_x \,&=\, \nev,
\label{eq:fratesaZ}
\end{align}
where $(\ref{eq:fratesZ})$~contains one equation for every~$x$.

\iim{c}$\flr$ is an infinitely differentiable function of the scores $v_{xy}$ as long as they keep satisfying the hypothesis of irreducibility.
\end{theorem}

\setlength\repskip{1.05em} 

\begin{proof}\hskip.5em
Let us begin by noticing that the hypothesis of irreducibility entails that $F$ can be extended to a continuous function on $\phisetb$ by putting $F(\psi)=0$ for $\psi\in\partial\phiset$. In order to prove this claim we must show that $F(\psi^n)\rightarrow 0$ whenever $\psi^n$ converges to a point $\psi\in\partial\phiset$. Let us consider the following sets associated with $\psi$: $\xst=\{x\,|\,\psi_x > 0\}$ 
and $\yst=\{y\,|\,\psi_y = 0\}$. The second one is not empty since we are assuming $\psi\in\partial\phiset$, whereas the first one is not empty because the strengths add up to the positive value~$f$. Now, for any $x\in\xst$ and $y\in\yst$, $F(\psi^n)$
contains a factor of the form $(\psi^n_y)^{v_{yx}}$, which tends to zero as soon as~$v_{yx}>0$. So, the only way for $F(\psi^n)$ not to approach zero would be $\vyx=\zeromatrix$, in contradiction with the irreducibility of $\vaa$.

After such an extension, $F$ is a continuous function on the compact set $\phisetb$.\linebreak 
So, there exists $\flr$ which maximizes $F$ on $\phisetb$. However, since $F(\psi)$ vanishes on $\partial\phiset$ whereas it is strictly positive for $\psi\in\phiset$, the maximizer $\flr$ must belong to $\phiset$.
This establishes the existence part of~(a).

Maximizing $F$ is certainly equivalent to maximizing $\log F$.
According to Lagrange, any $\flr\in\phiset$ which maximizes $\log F$  under the condition of a fixed sum is bound to satisfy
\begin{equation}
\label{eq:lagrange}
\frac{\partial\log F(\flr)}{\partial\flr_x} \,=\, \lambda,
\end{equation}
for some scalar $\lambda$ and every $x\in\ist$.
Now, a~straightforward computation gives
\begin{equation}
\label{eq:gradient}
\frac{\partial\log F(\flr)}{\partial\flr_x} \,=\,
\sum_{y\neq x}\left(\frac{v_{xy}}{\flr_x} - \frac{t_{xy}}{\flr_x+\flr_y}\right).
\end{equation}
On the other hand, using the fact that $v_{xy}+v_{yx} = t_{xy}$,
the preceding expression is easily seen to imply that
\begin{equation}
\label{eq:ortogonalitat}
\sum_x\,\frac{\partial\log F(\flr)}{\partial\flr_x}\,\flr_x \,=\, 0.
\end{equation}
In other words, the gradient of $\log F$ at $\flr$ is orthogonal to $\flr$, which was foreseeable since $F(\flr)$ remains constant when $\flr$ is multiplied by an arbitrary positive number. Notice that this is true for any $\flr$.
In particular, (\ref{eq:ortogonalitat}) entails that the above Lagrange multiplier $\lambda$ is equal to zero; in fact, it suffices to plug~(\ref{eq:lagrange}) in (\ref{eq:ortogonalitat}) and to use the fact that $\sum_x\flr_x = f$ is positive.
\ensep
So, the conditions (\ref{eq:lagrange}) reduce finally to
\begin{equation}
\label{eq:zerogradient}
\frac{\partial\log F(\flr)}{\partial\flr_x} \,=\, 0,
\end{equation}
for every $x\in\ist$, which is equivalent to (\ref{eq:fratesZ}) on account of (\ref{eq:gradient}) and the fact that $\flr_x>0$. So,~any maximizer must satisfy the conditions stated in~(b).

Let us see now that the maximizer is unique. 
Instead of following the interesting proof given by Zermelo,
here we will prefer to follow~\cite{ke}, which will have the advantage of preparing matters for part~(c). More specifically, the uniqueness will be obtained by seeing that any critical point of \,$\log F$ as a function on $\phiset$, \ie any solution of (\ref{eq:fratesZ}--\ref{eq:fratesaZ}), is a strict local maximum; this~implies that there is only one critical point, because otherwise
one should have other kinds of critical points~\cite[\S VI.6]{cou}
(we are invoking the so-called mountain pass theorem;
here we are using the fact that $\log F$ becomes $-\infty$ at $\partial\phiset$).
In order to study the character of a critical point we will look at the second derivatives of \,$\log F$ with respect to $\flr$. 
By~differentiating~(\ref{eq:gradient}), one obtains that
\begin{align}
\label{eq:hessianxx}
\frac{\partial^2\log F(\flr)}{\partial\flr_x{}^2} \,&=\,
-\, \sum_{y\neq x} \left(\frac{v_{xy}}{\flr^2_x} \,-\,
\frac{t_{xy}}{(\flr_x+\flr_y)^2}\right),
\\[2.5pt]
\label{eq:hessianxy}
\frac{\partial^2\log F(\flr)}{\partial\flr_x\,\partial\flr_y} \,&=\,
\frac{t_{xy}}{(\flr_x+\flr_y)^2},\qquad\text{ for }x\neq y.
\end{align}
On the other hand, when $\flr$ is a critical point, equation~(\ref{eq:fratesZ}) transforms (\ref{eq:hessianxx}) into the following expression:
\begin{equation}
\label{eq:hessianxxBis}
\frac{\partial^2\log F(\flr)}{\partial\flr_x{}^2} \,=\,
-\, \sum_{y\neq x}\,\frac{t_{xy}}{(\flr_x+\flr_y)^2}\,\frac{\flr_{y}}{\flr_{x}}.
\end{equation}
So, the Hessian bilinear form is as follows:
\begin{equation}
\label{eq:hessian}
\begin{split}
\sum_{x,y}\left(
\frac{\partial^2\log F(\flr)}{\partial\flr_x\,\partial\flr_y}
\right)\,\psi_x\,\psi_y
\,&=\, - \sum_{x,y\neq x}\,\frac{t_{xy}}{(\flr_x+\flr_y)^2}\,
\left(\frac{\flr_{y}}{\flr_x} \psi_x^2 - \psi_x\psi_y\right)\\
\,=\, - \sum_{x,y\neq x}\,&\frac{t_{xy}}{(\flr_x+\flr_y)^2\,\flr_x\flr_y}\,
\left(\flr_y^2\psi_x^2 - \flr_x\flr_y\psi_x\psi_y\right)\\
\,=\, -\, \sum_{\{x,y\}}\,&\frac{t_{xy}}{(\flr_x+\flr_y)^2\,\flr_x\flr_y}\,
\left(\flr_y\psi_x - \flr_x\psi_y\right)^2,
\end{split}
\end{equation}
where the last sum runs through all unordered pairs $\{x,y\}\sbseteq\ist$ with $x\neq y$.
The last expression is non-positive and it vanishes \ifoi
$\psi_x/\flr_x = \psi_y/\flr_y$ for any $x,y\in\ist$ 
(the ‘only~if’ part is immediate when $t_{xy}>0$; for~arbitrary  $x$ and $y$ the hypothesis of irreducibility allows to connect them through a path $x_0x_1\dots x_n$ ($x_0=x$, $x_n=y$) with the property that\linebreak 
$t_{x_ix_{i+1}}\ge v_{x_ix_{i+1}}>0$ for~any~$i$, so that one gets $\psi_x/\flr_x = \psi_{x_1}/\flr_{x_1} = \dots = \psi_y/\flr_y$).
So,~the vanishing of (\ref{eq:hessian}) happens \ifoi 
$\psi=\lambda\flr$ for some scalar $\lambda$. However, when $\psi$ is restricted to 
variations within $\phiset$, \ie to vectors in 
$T\phiset=\{\psi\in\bbr^\ist\,|\,\sum_x\psi_x=0\}$, the case $\psi=\lambda\flr$ reduces to $\psi=0$ (since $\sum_x\flr_x = f$ is~positive). So, the Hessian is negative definite on $T\phiset$. 
This ensures that $\flr$ is a strict local maximum of \,$\log F$ as a function on~$\phiset$.
In fact, one easily arrives at such a conclusion when Taylor's formula is used to analyse the behaviour of $\log F(\flr+\psi)$ for small $\psi$ in $T\phiset$.


Finally, let us consider the dependence of $\flr\in\phiset$ on the matrix~$\vaa$.\linebreak 
To~begin with, we notice that the set $\irreopen$ of irreducible matrices is open since it is a finite intersection of open sets, namely one open set for each splitting of $\ist$ into two sets $\xst$ and $\yst$.
The dependence of $\flr\in\phiset$ on $\vaa$ is due to the presence of $v_{xy}$ and $t_{xy}=v_{xy}+v_{yx}$ in the equations~(\ref{eq:fratesZ}--\ref{eq:fratesaZ}) which determine~$\flr$. However, we are not in the standard setting of the implicit function theorem since
we are dealing with a system of $N\cd+1$~equations whilst $\flr$ varies in a~space of dimension~$N\cd-1$.
In order to place oneself in a standard setting, it~is~convenient here to replace the condition of normalization $\sum_x\flr_x=\nev$ \,by~the alternative one\, $\flr_a=1$, where $a$ is a fixed element of $\ist$.
This change of normalization corresponds to mapping $\phiset$ \,to\, $U = \{\flr\in\bbr^\ist\,|\,\flr_x>0\text{ for all }x\in\ist,\, \flr_a=1\}$ by means of the diffeo\-morphism
$g:\flr\mapsto\flr/\flr_a$, which has the property that~$F(g(\flr)) = F(\flr)$. By taking as coordinates the $\flr_x$ with $x\in\ist\setminus\{a\}=:\ist'$, one easily checks that the function $F$ restricted to $U$ ---\ie restricted to $\flr_a=1$--- \,has the property that the matrix $(\,{\partial^2\log F(\flr)}/{\partial\flr_x\partial\flr_y}\mid x,y\cd\in\ist')$ is negative definite and therefore invertible, which entails that the system of equations $(\,{\partial\log F(\flr,\vaa)}/{\partial\flr_x}=0\mid x\cd\in\ist')$ determines $\flr\in U$ as a~smooth function of $\vaa\in\irreopen$.
\end{proof}

\medskip
Let us recall that a maximizing sequence 
means a sequence $\flrn\in\phiset$ such that $F(\flrn)$ approaches the lowest upper bound of $F$ on $\phiset$.

\smallskip
\begin{theorem}[Statements (a) and (b) are proved in \cite{ze}; results related to~(c) are contained in \cite{cg}]
\label{st:zermelore}
Assume that there exists
a top dominant irreducible component $\xst$.
In this case: 

\iim{a} There is a unique $\flr\in\phisetb$ such that 
any maximizing sequence 
converges to $\flr$.

\iim{b} $\flrx$ is the solution of a system analogous to \textup{(\ref{eq:fratesZ}--\ref{eq:fratesaZ})} where $x$ and $y$ vary only within $\xst$. $\flrax = 0$.

\iim{c} $\flr$ is a continuous function of the scores $v_{xy}$ as long as they keep satisfying the hypotheses of the present theorem.
\end{theorem}

\begin{proof}\hskip.5em
The definition of the lowest upper bound immediately implies the existence of maximizing sequences. On the other hand, the compactness of~$\phisetb$ guarantees that any maximizing sequence has a subsequence which converges in~$\phisetb$. 
Let $\flrn$ and $\flr$ denote respectively one of such convergent maximizing sequences and its limit.
In the following we will see that $\flr$ must be the unique point specified in statement~(b). This~entails that any maximizing sequence converges itself to $\flr$ (without extracting a subsequence).

So, our aim is now statement~(b).
From now on we will use the following notations:
a general element of $\phiset$ will be denoted by $\gflr$;
we will write\linebreak 
$\yst=\ist\setminus\xst$.
For convenience, in this part of the proof 
we will replace the condition $\sum_x\gflr_x = f$ by $\sum_x\gflr_x \le f$ (and similarly for $\flrn$ and $\flr$); since $F(\lambda\gflr) = F(\gflr)$ for any $\lambda>0$, 
the properties that we will obtain will be easily translated from 
$\widehat\phiset = \{\,\gflr\in\bbr^\ist\mid\gflr_x>0 
\text{ for all }x\in\ist,\ \sum_{x\in\ist}\gflr_x \le f\,\}$ \,to\, $\phiset$.
On the other hand, it will also be convenient to consider first the case where $\yst$ is also an irreducible component.
In such a case, it is interesting to  rewrite $F(\gflr)$ as a product of three factors:
\begin{equation}
F(\gflr) \,=\, \fxx(\gflrx)\,\fyy(\gflry)\,\fxy(\gflrx, \gflry),
\end{equation}
namely:
\begin{align}
\fxx(\gflrx) \,\,&=\, \prod_{\{x,\bar x\}\sbset\xst}\,
\frac{{\gflr_x}^{v_{x\bar x}}\,{\gflr_{\bar x}}^{v_{\bar xx}}}{(\gflr_x+\gflr_{\bar x})^{t_{x\bar x}}},
\label{eq:funcioFxx}
\\[2.5pt]
\fyy(\gflry) \,\,&=\, \prod_{\{y,\bar y\}\sbset\yst}\,
\frac{{\gflr_y}^{v_{y\bar y}}\,{\gflr_{\bar y}}^{v_{\bar yy}}}{(\gflr_y+\gflr_{\bar y})^{t_{y\bar y}}},
\label{eq:funcioFyy}
\\[2.5pt]
\fxy(\gflrx,\gflry) \,\,&=\, 
\kern6pt\prod_{{\scriptstyle x\in\xst \atop\scriptstyle y\in\yst }}\,\kern6pt
\left(\frac{\gflr_x}{\gflr_x+\gflr_y}\right)^{\hskip-3ptv_{xy}},
\label{eq:funcioFxy}
\end{align}
where we used that $v_{yx}=0$ and $t_{xy}=v_{xy}$. Now, let us look at the effect of replacing $\gflry$ by $\lambda\gflry$ without varying $\gflrx$. The values of $\fxx$ and $\fyy$ remain unchanged, but that of $\fxy$ varies in the following way:
\begin{equation}
\frac{\fxy(\gflrx,\lambda\gflry)}{\fxy(\gflrx,\gflry)} \,\,=\, 
\kern6pt\prod_{{\scriptstyle x\in\xst \atop\scriptstyle y\in\yst }}\,\kern6pt
\left(\frac{\gflr_x+\gflr_y}{\gflr_x+\lambda\gflr_y}\right)^{\hskip-3ptv_{xy}}.
\label{eq:factorFxy}
\end{equation}
In particular, for $0<\lambda<1$ each of the factors of the right-hand side of~(\ref{eq:factorFxy}) is greater than $1$. This remark leads to the following argument.
\ensep
First, we~can see that $\flrn_y/\flrn_x\rightarrow 0$ for any $x\in\xst$ and $y\in\yst$ such that $v_{xy}>0$ (such pairs $xy$ exist because of the hypothesis that $\xst$ dominates~$\yst$). Otherwise, the~preceding remark entails that the sequence $\flrbis^n = (\flrx^n,\lambda\flry^n)$ with $0<\lambda<1$ would satisfy \,$F(\flrbis^n) > K F(\flrn)$ for some $K>1$ and infinitely many $n$, in contradiction with the hypothesis that $\flrn$ was a maximizing sequence. On the other hand, we see also that $\fxy(\flrn)$ approaches its lowest upper bound, namely~$1$. 
\ensep
Having achieved such a property, the problem of maximizing $F$ reduces to separately maximizing $\fxx$ and $\fyy$, which is solved by Theorem~\ref{st:zermeloirre}. For~the moment we are dealing with relative strengths only, \ie without any normalizing condition like~(\ref{eq:fratesaZ}). So, we see that $\fyy$ gets optimized when each of the ratios $\flrn_y/\flrn_{\bar y}\ (y,\bar y\in\yst)$ approaches the homologous one for the unique maximizer of~$\fyy$, and analogously with~$\fxx$. Since these ratios are finite positive quantities, the statement that $\flrn_y/\flrn_x\rightarrow 0$ becomes extended to any $x\in\xst$ and $y\in\yst$ whatsoever (since one can write $\flrn_y/\flrn_x = (\flrn_y/\flrn_{\bar y})\times(\flrn_{\bar y}/\flrn_{\bar x})\times(\flrn_{\bar x}/\flrn_x)$ with $v_{\bar x\bar y}> 0$). Let us recover now the condition $\sum_{x\in\ist}\flrn_x=\nev$. The preceding facts imply that $\flry^n\rightarrow 0$, whereas $\flrx^n$ converges to the unique maximizer of $\fxx$. This establishes~(b) as well as the uniqueness part of~(a).

The general case where $\yst$ decomposes into several irreducible components, all of them dominated by $\xst$, can be taken care of by
induction over the different irreducible components of $\ist$.
At each step,
one deals with an irreducible component $\zst$ with the property of being minimal, in the sense of the dominance relation~$\chrela$, among those which are still pending. By~means of an argument analogous to that of the preceding paragraph, one sees that:\ensep
(i)~$\flrn_z/\flrn_x\rightarrow 0$ for any $z\in\zst$~and~$x$ such that $x\chrela z$ with $v_{xz}>0$;\ensep
(ii)~the ratios $\flrn_z/\flrn_{\bar z}\ (z,\bar z\in\zst)$ approach the homologous ones for the unique maximizer of~$\fzz$;\ensep
and (iii)~$\flrr^n$ is a maximizing sequence for $\frr$, where $R$~denotes the union of the pending components, $\zst$ excluded.
Once this induction process has been completed, one can combine its partial results to show that $\flrn_z/\flrn_x\rightarrow 0$ for any $x\in\xst$ and $z\not\in\xst$ (it suffices to consider a path $x_0x_1\dots x_n$ from $x_0\in\xst$ to $x_n=z$ with the property that $v_{x_ix_{i+1}}>0$ for any $i$ and to notice that each of the factors $\flrn_{x_{i+1}}/\flrn_{x_i}$ remains bounded while at least one of them tends to zero). As above, one concludes that $\flrax^n\rightarrow 0$, whereas $\flrx^n$ converges to the unique maximizer of $\fxx$.

\smallskip
The two following remarks will be useful in the proof of part~(c):
\linebreak 
(1)~\label{remarkOne}According to the proof above, $\flrx$~is determined (up to a multiplicative constant) by equations~(\ref{eq:fratesZ})
with $x$ and $y$ varying only within~$\xst$:
\smallskip\null\pagebreak 
\begin{equation}
\label{eq:yinx}
\fun_x(\flrx,\vaa) \,:=\, 
\sum_{\latop{\scriptstyle y\in\xst}{\scriptstyle y\neq x}}\,
t_{xy}\,\frac{\flr_x}{\flr_x+\flr_y} \,-\,
\sum_{\latop{\scriptstyle y\in\xst}{\scriptstyle y\neq x}}\,
v_{xy} \,=\, 0,\qquad\forall x\in\xst.\kern-10pt
\end{equation}
However, since $y\in\ist\setminus\xst$ implies on the one hand $\flr_y=0$ and on the other hand $t_{xy}=v_{xy}$, each of the preceding equations is equivalent to a similar one where $y$ varies over the whole of $\ist\setminus\{x\}$:
\begin{equation}
\label{eq:ally}
\funbis_x(\flr,\vaa) \,:=\,
\sum_{\latop{\scriptstyle y\in\ist}{\scriptstyle y\neq x}}\,
t_{xy}\,\frac{\flr_x}{\flr_x+\flr_y} \,-\,
\sum_{\latop{\scriptstyle y\in\ist}{\scriptstyle y\neq x}}\,
v_{xy} \,=\, 0,\qquad\forall x\in\xst.\kern-10pt
\end{equation}
(2)~\label{remarkTwo}Also, it is interesting to see the result of
adding up the equations (\ref{eq:ally}) for all $x$ in some subset $\wst$ of $\xst$. Using the fact that $v_{xy}+v_{yx}=t_{xy}$, one sees that such an addition results in the following equality:
\begin{equation}
\label{eq:sumw}
\sum_{\latop{\scriptstyle x\in\wst}{\scriptstyle y\not\in\wst}}\,
t_{xy}\,\frac{\flr_x}{\flr_x+\flr_y} \,-\,
\sum_{\latop{\scriptstyle x\in\wst}{\scriptstyle y\not\in\wst}}\,
v_{xy} \,=\, 0,\qquad\forall\kern.75pt\wst\sbseteq\xst.\kern-5pt
\end{equation}

\medskip
Let us proceed now with the proof of~(c). In the following, 
$\vaa$ and $\vaabis$\linebreak[3] denote respectively a fixed matrix satisfying the hypotheses of the theorem and a~slight perturbation of it.
As we have done in similar occasions, we systematically use a tilde to distinguish between hom\-olo\-gous objects associated respectively with $\vaa$~and~$\vaabis$; in particular, such a notation will be used in connection with the labels of certain equations.
Our aim is to~show that\, $\flrbis$ approaches $\flr$ \,as\, $\vaabis$ approaches $\vaa$. In this connection we will use the little-o and big-O notations made popular by Edmund Landau (who by the way is the author of a paper on the rating of chess players, namely~\cite{la}, which inspired Zermelo's work).\ensep
This notation refers here to 
functions of $\vaabis$ and their behaviour as $\vaabis$ approaches $\vaa$;
\ensep
if $f$ and $g$ are two such functions,
\ensep
$f=o(g)$\,~means that for every $\epsilon>0$ there exists a~$\delta>0$ such that $\|\vaabis-\vaa\|\le\delta$ implies $\|f(\vaabis)\|\le\epsilon\,\|g(\vaabis)\|$;
\ensep
on the other hand,\linebreak 
$f=O(g)$\, means that there exist $M$ and~$\delta>0$ such that $\|\vaabis-\vaa\|\le\delta$ implies $\|f(\vaabis)\|\le M\,\|g(\vaabis)\|$.

Obviously, if $\vaabis$ is near enough to $\vaa$ then $v_{xy}>0$ implies $\vbis_{xy}>0$.\linebreak 
As a consequence, $x\chrel y$ implies $x\chrelbis y$. In particular, the irreducibility 
of~$\vxx$ entails that $\vaabisxx$~is also irreducible.
Therefore, $\xst$ is entirely contained in~some irreducible component $\xstbis$ of $\ist$ for $\vaabis$.
Besides, $\xstbis$ is a top dominant irreducible component for $\vaabis$; in fact, we have the following chain of implications for $x\in\xst\sbseteq\xstbis$: $y\not\in\xstbis \,\Rightarrow\, y\not\in\xst \,\Rightarrow\, x\chrela y \,\Rightarrow\, x\chrelbis y \,\Rightarrow\, x\chrelabis y$, where we have used successively: the inclusion $\xst\sbseteq\xstbis$, the hypothesis that $\xst$ is top dominant for $\vaa$, the fact that $\vaabis$ is near enough to $\vaa$, and the hypothesis that $y$ does not belong to the irreducible component $\xstbis$.\linebreak 
Now, according to part~(b) and remark~(1) from p.\,\pageref{remarkOne}--\pageref{remarkTwo}, 
$\flrx$ and $\flrxbis$ are~determined respectively by the systems (\ref{eq:yinx}) and (\tref{eq:yinx}), or equivalently by~(\ref{eq:ally}) and (\tref{eq:ally}), whereas $\flrax$ and $\flraxbis$ are both of them equal to zero. So~we~must show that $\flrbisy=o(1)$ for any $y\in\xstbis\setminus\xst$, and that $\flrbisx-\flr_x=o(1)$ for~any $x\in\xst$. The proof is organized in three main steps.

\halfsmallskip
Step~(1).\ensep
\textit{$\flrbisy=O(\flrbisx)$ whenever $v_{xy}>0$}.\ensep
For the moment, we assume $\vaabis$ fixed (near enough to $\vaa$ so that $\vbis_{xy}>0$) and $x,y\in\xstbis$.
Under these hypotheses one can argue as follows: 
Since $\flrxbis$ maximizes $\fbisxx$, the corresponding value of $\fbisxx$ can be bounded from below by any particular value of the same function. On~the other hand, we can bound it from above by the factor $\flrbisx/(\flrbisx+\flrbisy)^{\vbis_{xy}}$. So, we can write
\begin{equation}
\label{eq:desigualtatsstep1}
\kern-10pt
\left(\frac12\right)^{\hskip-3ptN(N-1)} \le \left(\frac12\right)^{\lower4pt\hbox{$\latop{\textstyle\sum}{\vrule width0pt height8pt\scriptstyle p,q\in\xstbiss}$}{\textstyle \tbis_{pq}}} =
\fbisxx(\psi)
\,\le\, \fbisxx(\flrxbis) \,\le\,
\left(\frac{\flrbisx}{\flrbisx+\flrbisy}\right)^{\hskip-3pt\vbis_{xy}}\kern-2pt,\kern-2pt
\end{equation}
where $\psi$ has been taken so that $\psi_q$ has the same value for all $q\in\xstbis$ (and it~vanishes for $q\not\in\xstbis$). The preceding inequality entails that 
\begin{equation}
\label{eq:desigualtatsstep1bis}
\flrbisy \,\le\, \left(2^{\,N(N-1)\,/\,\vbis_{xy}} - 1\right) \,\flrbisx.
\end{equation}
Now, this inequality holds not only for $x,y\in\xstbis$,
but it is also trivially true for~$y\not\in\xstbis$,
since then one has $\flrbisy=0$. 
On the other hand, the case $y\in\xstbis,\ x\not\in\xstbis$ is not possible at all, because the hypothesis that
$\vbis_{xy}>0$ would then contradict the fact that $\xstbis$ is a top dominant irreducible component.
\ensep
Finally, we let $\vaabis$ vary towards $\vaa$. The~desired result is a consequence of~(\ref{eq:desigualtatsstep1bis}) since $\vbis_{xy}$ approaches $v_{xy}>0$.

\halfsmallskip
Step~(2).\ensep
\textit{$\flrbisy=o(\flrbisx)$ for any $x\in\xst$ and $y\not\in\xst$}.\ensep
Again, we will consider first the special case where $v_{xy}>0$.
In this case the result is easily obtained as a consequence of the equality (\tref{eq:sumw}) for $\wst=\xst$:
\begin{equation}
\label{eq:sumx}
\sum_{\latop{\scriptstyle x\in\xst}{\scriptstyle y\not\in\xst}}\,
\tbis_{xy}\,\frac{\flrbisx}{\flrbisx+\flrbisy} \,-\,
\sum_{\latop{\scriptstyle x\in\xst}{\scriptstyle y\not\in\xst}}\,
\vbis_{xy} \,=\, 0.
\end{equation}
In fact, this equality implies that
\begin{equation}
\label{eq:sumxbis}
\sum_{\latop{\scriptstyle x\in\xst}{\scriptstyle y\not\in\xst}}\,
\tbis_{xy}\,\left(1-\frac{\flrbisx}{\flrbisx+\flrbisy}\right) \,=\,
\sum_{\latop{\scriptstyle x\in\xst}{\scriptstyle y\not\in\xst}}\,
\vbis_{yx}.
\end{equation}
Now, it is clear that the right-hand side of this equation is $o(1)$ and that each of the terms of the left-hand side is positive. Since $\tbis_{xy} - v_{xy} = \tbis_{xy} - t_{xy} = o(1)$, the hypothesis that $v_{xy}>0$ allows to conclude that $\flrbisx/(\flrbisx+\flrbisy)$ approaches~$1$, or equivalently, $\flrbisy=o(\flrbisx)$.
\ensep
Let us consider now the case of any $x\in\xst$ and $y\not\in\xst$. Since $\xst$ is top dominant, we know that there exists a path $x_0x_1\dots x_n$ from $x_0=x$ to $x_n=y$ such that $v_{x_ix_{i+1}}>0$ for all $i$. According to step~(1) we have $\flrbissub_{x_{i+1}}=O(\flrbissub_{x_i})$. On the other hand, there must be some $j$ such that $x_j\in\xst$ but $x_{j+1}\not\in\xst$, which has been seen to imply that $\flrbissub_{x_{j+1}}=o(\flrbissub_{x_j})$. By combining these facts one obtains the desired result.

\halfsmallskip
Step~(3).\ensep
\textit{$\flrbisx-\flr_x=o(1)$ for~any $x\in\xst$}.\ensep
Consider the equations (\tref{eq:ally}) for $x\in\xst$ and split the sums in two parts depending on whether $y\in\xst$ or $y\not\in\xst$:
\begin{equation}
\label{eq:arranged}
\sum_{\latop{\scriptstyle y\in\xst}{\scriptstyle y\neq x}}\,
\tbis_{xy}\,\frac{\flrbisx}{\flrbisx+\flrbisy} \,-\,
\sum_{\latop{\scriptstyle y\in\xst}{\scriptstyle y\neq x}}\,
\vbis_{xy} \,=\,
\sum_{y\not\in\xst}\,(\vbis_{xy}-\tbis_{xy}\,\frac{\flrbisx}{\flrbisx+\flrbisy}).
\end{equation}
The last sum is $o(1)$ since step~(2) ensures that $\flrbisy=o(\flrbisx)$ and we also know that $\tbis_{xy}-\vbis_{xy}=\vbis_{yx} = o(1)$ (because $x\in\xst$ and $y\not\in\xst$). So $\flrbis$ satisfies a system of the following form, where $x$ and $y$ vary only within~$\xst$ \,and\, $\widetilde w_{xy}$ is a slight modification of $\vbis_{xy}$ which absorbs the right-hand side of (\ref{eq:arranged}):
\begin{equation}
\label{eq:vwdetall}
\funter_x(\flrbisX,\vaabis,\waabis) \,:=\,
\sum_{\latop{\scriptstyle y\in\xst}{\scriptstyle y\neq x}}\,
\tbis_{xy}\,\frac{\flrbisx}{\flrbisx+\flrbisy} \,-\,
\sum_{\latop{\scriptstyle y\in\xst}{\scriptstyle y\neq x}}\,
\widetilde w_{xy} \,=\, 0,\hskip1.5em\forall x\in\xst.\kern-5pt
\end{equation}
Here, the second argument of $\funter$ refers to the dependence on $\vaabis$ through~$\tbis_{xy}$.
We know that $\tbis_{xy} - t_{xy} = o(1)$ and also that $\widetilde w_{xy} - v_{xy} = (\widetilde w_{xy} - \vbis_{xy}) +\linebreak 
 (\vbis_{xy} - v_{xy}) = o(1)$.
So~we are interested in the preceding equation near the point $(\flrx,\vaa,\vaa)$. Now in this point we have $\funter(\flrx,\vaa,\vaa) = \fun(\flrx,\vaa) = 0$, as~well~as $(\partial\funter_x/\partial\flrbisy)(\flrx,\vaa,\vaa) = (\partial\fun_x/\partial\flr_y)(\flrx,\vaa)$. Therefore, the implicit function theorem can be applied similarly as in Theorem~\ref{st:zermeloirre}, with the  result that $\flrbisX={\cal H}(\vaabis,\waabis)$, where $\cal H$ is a smooth function which satisfies \hbox{${\cal H}(\vaa,\vaa)=\flrx$.} In particular, the continuity of $\cal H$ allows to conclude that $\flrbisX$~approaches $\flrx$, since we know that both $\vaabis$ and $\waabis$ approach $\vaa$.

\halfsmallskip
Finally, by combining the results of steps~(2) and (3) one obtains $\flrbisy\cd=o(1)$ for any $y\not\in\xst$.
\end{proof}

\vskip-4pt
\remarks

1. The convergence of $\flrn$ to $\flr$ is a necessary condition for $\flrn$ being a maximizing sequence but not a sufficient one. The preceding proof shows that a~necessary and sufficient condition is that the ratios $\flrn_y/\flrn_z$ tend to $0$\linebreak 
whenever $y\chrela z$, whereas, if $y\chrels z$, \ie if $y$ and $z$ belong to the same irreducible component $Z$, these ratios approach 
the homologous ones for the unique maximizer of~$\fzz$.

2. If there is not a dominant component then the maximizing sequences can have multiple limit points. However, as we will see in the next section, the projected Llull matrices are always in the hypotheses of Theorem~\ref{st:zermelore}.

\section{The fraction-like rates}

Let us recall from~\secpar{2.9} that the fraction-like rates~$\flr_x$ will be obtained by applying Zermelo's method to the projected Llull matrix~$(\psc_{xy})$.

\bigskip
The next results show that this matrix has a very special structure in connection with irreducibility.



\smallskip
\begin{lemma}\hskip.5em
\label{lema_irreductibilitat_projectada}
The projected Llull matrix $(\psc_{xy})$ has the following properties for any admissible order $\xi$ ($p'$ denotes the immediate successor of $p$ in $\xi$):

\iim{a} If $x\rxi y$ and $\psc_{yx}=0$, then $\psc_{p'p}=0$ for some $p$ such that $x\rxieq p \rxi y$.
 
\iim{b} If $\psc_{p'p}=0$ for some $p$, then $\psc_{yx}=0$ for all $x,y$ such that $x\rxieq p\rxi y$.
 
\iim{c} If $x\rxi y$ and $\psc_{xy}=0$, then $\psc_{ab}=0$ for all $a,b$ such that $x\rxieq a$.

\end{lemma}
\begin{proof}\hskip.5em
Part~(a).\ensep
Assume that $x\rxi y$. Then $\psc_{yx}$ is the left end of the interval $\gamma_{xy}$. Now, since $\gamma_{xy}=\bigcup\,\{\gamma_{pp'}\,|\,x\rxieq p\rxi y\}$, a vanishing left end for~$\gamma_{xy}$ implies the same property for some of the $\gamma_{pp'}$, \ie $\psc_{p'p}=0$.
  
\halfsmallskip
Part~(b).\ensep 
According to Theorem~\ref{st:propertiesOfProjection}.(a), $x\rxieq p\rxi y$ implies the inequalities 
$\psc_{yx}\le \psc_{p'x}\le \psc_{p'p}$. Therefore, 
$\psc_{p'p}=0$ implies  $\psc_{yx}=0$.
  
\halfsmallskip
Part~(c).\ensep
For $x\rxi y$, $\psc_{xy}=0$ means that $\gamma_{xy}=[0,0]$.
This implies that $\gamma_{pp'}=[0,0]$ for all $p$ such that $x\rxieq p\rxi y$. Now, according to Lemma~\ref{st:intervals}, the barycentres of the intervals $\gamma_{qq'}$ decrease or stay the same when $q$ moves towards the bottom. So $\gamma_{qq'}=[0,0]$ for all $q$ such that $x\rxieq q$. As a consequence, we immediately get $\psc_{ab} = 0$ for any $a,b$ such that $x\rxieq a,b$. Furthermore, for $b\rxi x\rxieq a$, part~(a) of Theorem~\ref{st:propertiesOfProjection} gives the following inequalities: 
$\psc_{ab}\le\psc_{ax}$ for $a\ne x$, and $\psc_{ab}\le\psc_{ay}$ for $a=x$, where the right-hand sides are already known to vanish. So $\psc_{ab}$ vanishes also for such~$a$ and~$b$.
\end{proof}

\begin{proposition}\hskip.5em
\label{st:existencia-X}
Let us assume that the projected Llull matrix $(\psc_{xy})$
is not the zero matrix.
Let us consider the set
\begin{equation}
\xst \,=\, \{x\in\ist\mid \psc_{p'p}>0 \text{ for all $p$ such that }p\rxi x\},
\end{equation}
where the right-hand side makes use of an admissible order~$\xi$.
This set has the following properties:

\iim{a} It does not depend on the admissible order $\xi$.\par

\iim{b} $\psc_{xy}>0$ for any $x\in\xst$ and $y\in\ist$.\par

\iim{c} $\psc_{yx}=0$ for any $x\in\xst$ and $y\not\in\xst$.\par

\iim{d} $r_x<r_y$ for any $x\in\xst$ and $y\not\in\xst$.\par

\iim{e} $\xst$~is the top dominant irreducible component of $\ist$ for $(\psc_{xy})$.\par
\end{proposition}












\newcommand\topzero{h}

\begin{proof}\hskip.5em
Statement~(a) will be proved at the end.
The definition of $\xst$ is equivalent to the following one:
\ensep
$\xst=\ist$\, if~$\psc_{p'p}>0$ for any $p$;
\ensep
otherwise, $\xst=\{x\in\ist\,|\,x\rxieq\topzero\}$, where $\topzero$ is the topmost (in~$\xi$) element of $\ist$ which satisfies $\psc_{\topzero'\topzero}=0$.
In particular, $\xst$ reduces to the topmost element of $\ist$\, when~$\psc_{p'p}=0$ for any $p$.

Statement~(b).\ensep
In view of Lemma~\ref{lema_irreductibilitat_projectada}.(a), the definition of $\xst$ implies that $\psc_{yx}>0$ for any $x,y\in\xst$ such that $x\rxi y$.
\ensep
This statement is empty when $\xst$ reduces to a single element $a$, but then we will make use of the fact that $\psc_{aa'}>0$, which is bound to happen because otherwise Lemma~\ref{lema_irreductibilitat_projectada}.(c) would entail that the whole matrix is zero, against our hypothesis.
\ensep
These facts imply statement~(b) by virtue of Theorem~\ref{st:propertiesOfProjection}.(a).

Statement~(c).\ensep
If $\xst=\ist$ there is nothing to prove. Otherwise, if $\topzero$ is the above-men\-tioned topmost  element of $\ist$ which satisfies $\psc_{\topzero'\topzero}=0$, then Lemma~\ref{lema_irreductibilitat_projectada}.(b) ensures that $\psc_{yx}=0$ for any $x,y$ such that $x\rxieq \topzero\rxi y$, \ie any $x\in\xst$ and $y\not\in\xst$.

Statement~(d).\ensep
If $\xst=\ist$ there is nothing to prove.
Otherwise, the result follows from parts~(b) and (c) together with Lemma~\ref{st:obsRosa}.(c).

Statement~(e).\ensep
This is an immediate consequence of (b) and (c).

Statement~(a).\ensep
A top dominant irreducible component is always unique because the relation of dominance between irreducible components is antisymmetric.
\end{proof}

\vskip-4pt
\remarks
 1.\ensep In the complete case, the average ranks $\bar r_x$ defined by equation~(\ref{eq:avranksfromscores})\linebreak[3] are easily seen to satisfy already a property of the same kind as~(d):\ensep if $\xst$ and $\yst$ are two irreducible components of $(v_{xy})$ such that $\xst$ dominates $\yst$, then $\bar r_x < \bar r_y$ for all $x\in\xst$ and $y\in\yst$
\cite[Thm.\,2.5]{mp}.\par

2.\ensep Even in the complete case, Zermelo's rates associated with the original Llull matrix $(v_{xy})$ are not necessarily compatible with the average ranks $\bar r_x$. However, as we will see below, the projected Llull matrices will always enjoy such a compatibility.\par

\bigskip
From now on, $\xst$ denotes the top dominant irreducible component whose existence is established by the preceding proposition. According to  Theorem~\ref{st:zermelore}, the fraction-like rates $\flr_x$ vanish \ifoi $x\in\ist\setminus\xst$ and their values for $x\in\xst$ are determined by
the restriction of $(\psc_{xy})$ to $x,y\in\xst$.
More specifically, the latter are determined by the condition of maximizing the function
\begin{equation}
\label{eq:funcioFP}
F(\flr) \,\,=\, \prod_{\{x,y\}}\,
\frac{{\flr_x}^{\psc_{xy}}\,{\flr_{y}}^{\psc_{yx}}}{(\flr_x+\flr_{y})^{\pto_{xy}}},
\end{equation}
under the restriction
\begin{equation}
\begin{repeated}{eq:fratesa}
\label{eq:fratesaP}
\sum_x\,\flr_x \,\,=\, \nev.
\end{repeated}
\end{equation}
where we will understand that $x$ and $y$ are restricted to $\xst$, and 
$\nev$~denotes the fraction of non-empty votes (\ie $\nev = F/V$ where $F$~is the number of non-empty votes and $V$ is the total number of votes). Moreover, we know that $(\flr_x\,|\,x\in\xst)$ is the solution of the following system of equations besides~(\ref{eq:fratesaP}):
\begin{equation}
\begin{repeated}{eq:frates}
\label{eq:fratesP}
\sum_{y\neq x}\,
\pto_{xy}\,\frac{\flr_x}{\flr_x+\flr_y} \,=\,
\sum_{y\neq x}\,
\psc_{xy}.
\end{repeated}
\end{equation}
where the sums extend to all $y\ne x$ in $\xst$.
\ensep
The next result shows that the resulting fraction-like rates are fully compatible with the rank-like ones except for the vanishing of those outside the top dominant component. 

\smallskip
\begin{theorem}\hskip.5em
\label{st:phis}

\iim{a} $\flr_x > \flr_y \,\Longrightarrow\, r_x \,<\, r_y$.

\iim{b} $r_x \,<\, r_y \,\Longrightarrow\, \hbox{either\, }\flr_x > \flr_y \hbox{ \,or\, }\flr_x = \flr_y = 0$.
\end{theorem}

\begin{proof}\hskip.5em
Let us begin by noticing that both statements hold if $\flr_y=0$, \ie if $y\not\in\xst$. In the case of statement~(a), this is true because of Proposition~\ref{st:existencia-X}.(d). So, we can assume that $\flr_y>0$, \ie $y\in\xst$. But in this case, each one of the hypotheses of the present theorem implies that $\flr_x>0$, \ie $x\in\xst$. In the case of statement~(b), this is true because of
Proposition~\ref{st:existencia-X}.(d) (with $x$ and $y$ interchanged with each other) and the fact that $\xst$ is a top interval for any admissible order.
\ensep
So, from now on we can assume that $x$ and $y$ are both in $\xst$, or, on account of Theorem~\ref{st:zermelore}, that $\xst=\ist$.

\halfsmallskip
Statement~(a): It will be proved by seeing that a simultaneous occurrence of the inequalities $\flr_x > \flr_y$ and $r_x \ge r_y$  would entail a contradiction with the fact that $\flr$ is the unique maximizer of $F(\flr)$. More specifically, we will see that one would have $F(\flrbis)\ge F(\flr)$ where $\flrbis$ is obtained from $\flr$ by interchanging the values of $\flr_x$ and $\flr_y$, that is 
\begin{equation}
\label{eq:intercanvi}
\flrbis_z \,=\,
\begin{cases}
\flr_y, &\text{if }z=x,\\
\flr_x, &\text{if }z=y,\\
\flr_z, &\text{otherwise.}\\
\end{cases}
\end{equation}
In fact, $\flrbis$ differs from $\flr$ only in the components  associated with $x$ and $y$, so that
\begin{equation}
\begin{split}
\label{eq:quocientFs}
\frac{F(\flrbis)}{F(\flr)} \,=\,
&\left(\frac{\flrbisx}{\flr_x}\right)^{\hskip-3pt\psc_{xy}}
\prod_{z\ne x,y}\,
\left(\frac{\flrbisx/(\flrbisx+\flr_z)}{\flr_x/(\flr_x+\flr_z)}
\right)^{\hskip-3pt\psc_{xz}}
\left(\frac{\flr_x+\flr_z}{\flrbisx+\flr_z}\right)^{\hskip-3pt\psc_{zx}}\\
\times
&\left(\frac{\flrbisy}{\flr_y}\right)^{\hskip-3pt\psc_{yx}}
\prod_{z\ne x,y}\,
\left(\frac{\flrbisy/(\flrbisy+\flr_z)}{\flr_y/(\flr_y+\flr_z)}
\right)^{\hskip-3pt\psc_{yz}}
\left(\frac{\flr_y+\flr_z}{\flrbisy+\flr_z}\right)^{\hskip-3pt\psc_{zy}}.
\end{split}
\end{equation}
More particularly, in the case of (\ref{eq:intercanvi}) this expression becomes
\begin{equation}
\label{eq:canviphis-aF}
\frac{F(\flrbis)}{F(\flr)} \,=\,
\left(\frac{\flr_y}{\flr_x}\right)^{\hskip-3pt\psc_{xy}-\psc_{yx}}\,
\prod_{z\ne x,y}\,
\left(\frac{\flr_y/(\flr_y+\flr_z)}{\flr_x/(\flr_x+\flr_z)}
\right)^{\hskip-3pt\psc_{xz}-\psc_{yz}}
\left(\frac{\flr_y+\flr_z}{\flr_x+\flr_z}\right)^{\hskip-3pt\psc_{zy}-\psc_{zx}},
\end{equation}
where all of the bases are strictly less than~$1$, since 
$\flr_x > \flr_y$, and all of the the exponents are non-positive, because of Lemma~\ref{st:obsRosa}.(c). Therefore, the product is greater than or equal to $1$, as claimed.

\newcommand\flro{\omega}
\newcommand\incr{\epsilon}

\halfsmallskip
Statement~(b): Since we are assuming $x,y\in\xst$, it is a matter of proving that $r_x < r_y \,\Rightarrow\, \flr_x > \flr_y$. On the other hand, by making use of the contra\-positive of~(a), the problem reduces to proving that $\flr_x = \flr_y \,\Rightarrow\, r_x = r_y$.

Similarly to above, this implication will be proved by seeing that a simultaneous occurrence of the equality $\flr_x=\flr_y=:\flro$ together with the inequal\-ity $r_x < r_y$ (by symmetry it suffices to consider this one) would entail a contradiction with the fact that $\flr$ is the unique maximizer of $F(\flr)$.\linebreak[3]
 More specifically, here we will see that one would have $F(\flrbis)>F(\flr)$ where $\flrbis$ is obtained from $\flr$ by slightly increasing $\flr_x$ while decreasing $\flr_y$, that is
\begin{equation}
\label{eq:dissociacio}
\flrbisz \,=\,
\begin{cases}
\flro + \incr, &\text{if }z=x,\\
\flro - \incr, &\text{if }z=y,\\
\flr_z, &\text{otherwise.}\\
\end{cases}
\end{equation}
This claim will be proved by checking that
\begin{equation}
\label{eq:depspositiva}
\left.
\frac{\textup{d}\hphantom{\incr}}{\textup{d}\incr}\,
\log \frac{F(\flrbis)}{F(\flr)}
\,\right|_{\epsilon=0} \,>\, 0.
\end{equation}
In fact, (\ref{eq:quocientFs}) entails that
\begin{equation}
\begin{split}
\log\,\frac{F(\flrbis)}{F(\flr)} \,=\,\,\,
&C \,+\, \psc_{xy}\log\flrbisx + \psc_{yx}\log\flrbisy\\
+ &\sum_{z\ne x,y}\,\Big(
\psc_{xz}\log\frac{\flrbisx}{\flrbisx+\flr_z}
+ \psc_{yz}\log\frac{\flrbisy}{\flrbisy+\flr_z}\Big)\\
- &\sum_{z\ne x,y}\,\Big(
\psc_{zy}\log(\flrbisy+\flr_z)
+ \psc_{zx}\log(\flrbisx+\flr_z)\Big),
\end{split}
\end{equation}
where $C$ does not depend on $\incr$. Therefore, in view of (\ref{eq:dissociacio}) we get
\begin{equation}
\begin{split}
\left.
\frac{\textup{d}\hphantom{\incr}}{\textup{d}\incr}\,
\log \frac{F(\flrbis)}{F(\flr)}
\,\right|_{\epsilon=0}
\,=\,\,\,
(\psc_{xy}-\psc_{yx})\,\frac{1}{\flro}
\,\,+\, &\sum_{z\ne x,y}\,
(\psc_{xz}-\psc_{yz})\,\frac{\flr_z}{\flro(\flro+\flr_z)}\\
+\, &\sum_{z\ne x,y}\,
(\psc_{zy}-\psc_{zx})\,\frac{1}{\flro+\flr_z}.
\end{split}
\end{equation}
Now, according to Lemma~\ref{st:obsRosa}.(b,\,c), the assumption that $r_x<r_y$ implies the inequalities $\psc_{xy}>\psc_{yx}$, $\psc_{xz}\ge\psc_{yz}$ and $\psc_{zy}\ge\psc_{zx}$, which result indeed in~(\ref{eq:depspositiva}).
\end{proof}

\renewcommand\upla{\vskip-3pt}
\medskip\upla
The next proposition establishes property~\llpv:

\smallskip\upla
\begin{proposition}\hskip.5em
\label{st:plumpflr}
In the case of plumping votes the fraction-like rates coincide with the fractions of the vote obtained by each option.
\end{proposition}

\vskip-9mm\null
\begin{proof}\hskip.5em
Proposition~\ref{st:plumpps} ensures that the projected scores coincide with the original ones. So we have $\psc_{xy}=\plumpf_x$ and $\pto_{xy}=\plumpf_x+\plumpf_y$. In these circumstances it is obvious that equations (\ref{eq:fratesaP}--\ref{eq:fratesP}) are satisfied if we take $\flr_x=\plumpf_x$. So it suffices to invoke the uniqueness of solution of this system.
\end{proof}

\vskip-11mm\null
\section{Continuity}

We claim that both the rank-like rates $r_x$ and the fraction-like ones $\flr_x$ are continuous functions of the binary scores~$v_{xy}$. The main difficulty in proving this statement lies in the admissible order~$\xi$, which plays a central role in the computations. Since $\xi$ varies in a discrete set, its dependence on the data cannot be continuous at~all. Even so, we claim that the final result is still a~continuous function of the data.

In this connection, one can consider as data the normalized Llull matrix~$(v_{xy})$, its domain of variation being the set $\Omega$ introduced in \secpar{3.3}.
\ensep
Alternatively, one can consider as data the relative frequencies of the possible votes, \ie the coefficients $\alpha_k$ mentioned also in \secpar{3.3}.

\smallskip
\begin{theorem}\hskip.5em
\label{st:continuityThm}
The following objects depend continuously on the Llull matrix~$(v_{xy})$:
the~projected scores~$\psc_{xy}$,
the~rank-like rates~$r_x$, and\,
the~fraction-like rates~$\flr_x$.
\end{theorem}

\begin{proof}\hskip.5em
Let us begin by considering the dependence of the rank-like rates and the fraction-like rates on the projected scores. In the case of the rank-like rates, this dependence is given by formula~(\ref{eq:rrates}), which is not only continuous but even linear (non-homogeneous).\ensep
In the case of the fraction-like rates, their dependence on the projected scores is more involved, but is is still continuous.
In fact, Theorem~\ref{st:zermelore}.(c) ensures such a continuity under the hypothesis that there is a top irreducible component,
which hypothesis is~satisfied by virtue of Proposition~\ref{st:existencia-X}.(e).

So we are left with the problem of showing that the projection $P:(v_{xy})\mapsto(\psc_{xy})$ is continuous. As it has been mentioned above, this is not so clear, since the projected scores are the result of certain operations which are based upon an admissible order~$\xi$  which is determined separately. However, we will see that, on~the one hand, $P$~is continuous as long as $\xi$~remains unchanged, and on the other hand, the results of~\secpar{8,\,9} allow to conclude that $P$~is continuous on the whole of $\Omega$ in spite of the fact that $\xi$ can change. In the following we will use the following notation: for every total order $\xi$, we denote by $\Omega_\xi$ the subset of~$\Omega$ which consists of the Llull matrices for which $\xi$ is an~admissible order, and we denote by $P_\xi$ the restriction of $P$ to $\Omega_\xi$. 

We claim that the mapping $P_\xi$ is continuous for every total order $\xi$.\linebreak 
In~order to check the truth of this statement, one has to go over the different mappings whose composition defines $P_\xi$ (see~\secpar{2.8}), namely:\ensep
$(v_{xy})\mapsto(\isc_{xy})\mapsto(\img_{xy})$,\ensep
$(v_{xy})\mapsto(t_{xy})$,\ensep
$(\img_{xy})\mapsto(\ppmg_{xy})$,\ensep
$\Psi:((\ppmg_{xx'}),(t_{xy}))\mapsto(\ppto_{xy})$,\ensep
and finally $((\ppmg_{xx'}),(\ppto_{xx'}))\mapsto(\psc_{xy})$.
\ensep
Quite a few of these mappings involve the \,$\max$\, and \,$\min$\, operations, which are certainly continuous. For instance, the last mapping above can be written as
$\psc_{xy} = \max\,\{\,(\ppto_{pp'}+\ppmg_{pp'})/2\mid x\rxieq p\rxi y\}$ 
 and $\psc_{yx} = \min\,\{\,(\ppto_{pp'}-\ppmg_{pp'})/2\mid x\rxieq p\rxi y\}$ for $x\rxi y$.
\ensep
Concerning the operator $\Psi$, let us recall that its output is the orthogonal projection of $(t_{xy})$ onto a certain convex set determined by $(\ppmg_{xx'})$; a general result of continuity for such an operation can be found in~\cite{da}.

Finally, the continuity of $P$ (and the fact that it is well-defined) is a~consequence of the following facts (see for instance \cite[\S2-7]{mk}):\ensep
(a)~$\Omega=\bigcup_{\xi}\Omega_\xi$; this is true because of the existence of $\xi$ (Corollary~\ref{st:existenceXiCor}).\ensep
(b)~$\Omega_\xi$ is a closed subset of $\Omega$; this is true because $\Omega_\xi$ is described by a set of non-strict inequalities which concern quantities that are continuous functions of $(v_{xy})$ (namely the inequalities $\img_{xy}\ge 0$ whenever $xy\in\xi$).\ensep
(c)~$\xi$ varies over a finite set.\ensep
(d)~$P_\xi$ coincides with $P_{\eta}$ at $\Omega_\xi\cap\Omega_\eta$, as it is proved in Theorem~\ref{st:independenceOfXi}.
\end{proof}


\smallskip
\begin{corollary}\hskip.5em
\label{st:continuityCor}
The rank-like rates, as well as the fraction-like ones,\linebreak[3] depend continuously on the relative frequency of each possible content of an~individual vote.
\end{corollary}

\begin{proof}\hskip.5em
It suffices to notice that the Llull matrix $(v_{xy})$ is simply the center of gravity of the distribution specified by these relative frequencies (formula~(\ref{eq:cog}) of~\secpar{3.3}).
\end{proof}

\section{Decomposition}

Properties~\llrd\ and~\llfd\ are concerned with having a partition of~$\ist$ in two sets $\xst$ and $\yst$ such that the rates for $x\in\xst$ can be obtained by restricting the attention to $\vxx$,
\ie the $v_{x\bar x}$ with $x,\bar x\in\xst$
(and similarly for $y\in\yst$ in the case of property~\llrd).

More specifically, property~\llrd\ considers the case where the
following equalities are satisfied:
\begin{alignat}{2}
r_x \,&=\, \rbis_x,\qquad &&\hbox{for all $x\in\xst$},
\label{eq:condrfx}\\
r_y \,&=\, \rbis_y \,+\, |X|,\qquad &&\hbox{for all $y\in\yst$},
\label{eq:condrfy}
\end{alignat}
where $\rbis_x$ and $\rbis_y$ denote the rank-like rates which are determined respectively by the matrices $\vxx$ and $\vyy$.
Property~\llrd\ states that in the complete case these equalities are equivalent to having
\begin{equation}
\label{eq:v1}
v_{xy} \,=\, 1
\quad \hbox{(and therefore $v_{yx}=0$)}
\quad \hbox{whenever\, $xy\in\xst\times\yst$}.
\end{equation}

In the following we will continue using a tilde to distinguish
between hom\-olo\-gous objects associated respectively with the
whole matrix~$\vaa$ and with its submatrices $\vxx$ and $\vyy$.

\bigskip
First of all we explore the implications of condition~(\ref{eq:v1}).

\smallskip
\begin{lemma}\hskip.5em
\label{lema1-vxy}
Given a partition $\ist=\xst\cup\yst$ in two disjoint nonempty sets, one has the following implications:
\begin{equation}
\left.
\begin{array}{c} v_{xy}=1 \\[2.5pt]
\forall\,xy\in\xst\times\yst
\end{array}
\right\}
\ \Longrightarrow\ 
\left\{
\begin{array}{c} \img_{xy}\cd=1 \\[2.5pt]
\forall\,xy\in\xst\times\yst
\end{array}
\right\}
\ \Longleftrightarrow\ 
\left\{
\begin{array}{c}
v_{xy}^\pi=1 \\[2.5pt]
\forall\,xy\in\xst\times\yst
\end{array}
\right.
\label{eq:lema1-vxy}
\end{equation}
If the individual votes are complete, or alternatively, if they are
transitive relations, then the converse of the first implication holds too.
\end{lemma}

\begin{proof}\hskip.5em
Assume that $v_{xy}=1$ for all $xy\in\xst\times\yst$. Then $v_{yx}=0$, for all such pairs, which implies that $v_\gamma$ vanishes for any path $\gamma$ which goes from $\yst$ to~$\xst$. This fact, together with the inequality $\isc_{xy}\ge v_{xy}$, entails the following equalities for all $x\in\xst$ and $y\in\yst$: $\isc_{yx}=0$, $\isc_{xy}=1$, and consequently $\img_{xy}=1$.

Assume now that $\img_{xy}=1$ for all $xy\in\xst\times\yst$. Let $\xi$ be an admissible order. As an immediate consequence of the definition, it includes the set $\xst\times\yst$. Let $\last$ be the last element of $\xst$ according to~$\xi$. From the present hypothesis it is clear that $\ppmg_{\last\last'}=1$, \ie $\gamma_{\last\last'}=[0,1]$, which entails that $\gamma_{xy}=[0,1]$, \ie $\psc_{xy}=1$, for every $xy\in\xst\times\yst$.

Assume now that $\psc_{xy}=1$ for all $xy\in\xst\times\yst$. Let $\xi$ be an admissible order. Here too, we are ensured that it includes the set $\xst\times\yst$; this is so by virtue of Theorem~\ref{st:propertiesOfProjection}.(a). Let $\last$ be the last element of $\xst$ according to~$\xi$. From the fact that $\ppmg_{\last\last'}=\pmg_{\last\last'}=1$, one infers that $\img_{xy}=1$ for all $xy\in \xst\times\yst$.

Finally, let us assume again that $\img_{xy}=1$ for all $xy\in\xst\times\yst$. Since\linebreak $\img_{xy}=\isc_{xy}-\isc_{yx}$ and both terms of this difference belong to $[0,1]$, the only possibility is $\isc_{xy}=1$ and $\isc_{yx}=0$, which implies that $v_{yx}=0$.\ensep
In the complete case, this equality is equivalent to $v_{xy}=1$.
\ensep
In the case where the individual votes are transitive relations, one can reach the same conclusion in the following way: The equality $\isc_{xy} = 1$ implies the existence of a path $x_0x_1\dots x_n$ from $x$ to $y$ such that $v_{x_ix_{i+1}} = 1$ for all $i$. But this means that all of the votes include each of the pairs $x_ix_{i+1}$ of this path. So, if they are transitive relations, all~of them include also the pair $xy$, \ie $v_{xy} = 1$.
\end{proof}

\smallskip
\begin{lemma}\hskip.5em
\label{lema3-vxy}
Condition~\textup{(\ref{eq:v1})} implies, for any admissible order, the following equalities:
\begin{alignat}{2}
\label{eq:ppmgxx}
\ppmg_{xx'} &=\, \ppmgbis_{xx'}, \qquad&&\hbox{whenever $x,x'\in\xst$,}\\[2.5pt]
\label{eq:ppmgyy}
\ppmg_{yy'} &=\, \ppmgbis_{yy'}, \qquad&&\hbox{whenever $y,y'\in\yst$,}
\\[2.5pt]
\label{eq:ptoxx}
\pto_{x\bar x} &=\, 1, \qquad&&\hbox{for all $x,\bar x\in\xst$.}
\end{alignat}
\end{lemma}
\begin{proof}\hskip.5em
As we saw in the proof of Lemma~\ref{lema1-vxy}, condition~(\ref{eq:v1}) implies the vanishing of $v_\gamma$ for any path $\gamma$ which goes from $\yst$ to~$\xst$. Besides the conclusions obtained in that lemma, this implies also the following equalities:
\begin{alignat}{3}
\label{eq:iscxx}
\isc_{x\bar x} &= \iscbis_{x\bar x},
\hskip1.75em &\img_{x\bar x} &= \imgbis_{x\bar x},
\hskip1.75em &&\hbox{for all $x,\bar x\in\xst$,}\\[2.5pt]
\label{eq:iscyy}
\isc_{y\bar y} &= \iscbis_{y\bar y},
\hskip1.75em &\img_{y\bar y} &= \imgbis_{y\bar y},
\hskip1.75em &&\hbox{for all $y,\bar y\in\yst$.}
\end{alignat}
Let us fix an admissible order~$\xi$. The second equality of (\ref{eq:lema1-vxy}) ensures not only that $\xi$ includes the set $\xst\times\yst$, but it can also be combined with (\ref{eq:iscxx}) and (\ref{eq:iscyy}) to obtain respectively (\ref{eq:ppmgxx}) and (\ref{eq:ppmgyy}).
On the other hand, the third equality of (\ref{eq:lema1-vxy}) implies that $\pto_{xy}=1$ for all $xy\in\xst\times\yst$, from which the pattern of growth of the projected turnovers ---more specifically, equation~(\ref{eq:tinequalities}.2)--- allows to obtain (\ref{eq:ptoxx}).
\end{proof}

\smallskip
\begin{theorem}\hskip.5em
\label{st:condrfPro}
In the complete case\, one has the following equivalences: 
\textup{(\ref{eq:condrfx})} $\Longleftrightarrow$
\textup{(\ref{eq:condrfy})} $\Longleftrightarrow$
\textup{(\ref{eq:v1})}.
\end{theorem}

\begin{proof}\hskip.5em
Since we are considering the complete case, we can make use of the margin-based procedure (\secpar{2.6}). The proof is organized in two parts:

\halfsmallskip
Part~(a):\ensep (\ref{eq:v1}) $\Longrightarrow$
(\ref{eq:condrfx}) and (\ref{eq:condrfy}).\ensep
As a consequence of the equalities~(\ref{eq:ppmgxx}) and (\ref{eq:ppmgyy}), the margin-based procedure ---more specifically, steps (\ref{eq:cprojection3}) and (\ref{eq:cprojection4})--- results in the following equalities:
\begin{alignat}{2}
\label{eq:pmgxx}
\pmg_{x\bar x} \,&=\, \pmgbis_{x\bar x}, \qquad&&\hbox{for all $x,\bar x\in\xst$,}\\[2.5pt]
\label{eq:pmgyy} \pmg_{y\bar y} \,&=\, \pmgbis_{y\bar y},
\qquad&&\hbox{for all $y,\bar y\in\yst$.}
\end{alignat}
On the other hand, the third equality of (\ref{eq:lema1-vxy}) is equivalent to saying that
\begin{equation}
\hskip1.7em\pmg_{xy} \,=\, 1, \qquad\hskip1.1em\hbox{for all $xy\in\xst\times\yst$.}
\end{equation}
When the projected margins are introduced in (\ref{eq:rratesfrommargins}) these equalities result in (\ref{eq:condrfx}) and (\ref{eq:condrfy}).

\halfsmallskip
Part~(b):\ensep (\ref{eq:condrfx}) $\Rightarrow$ (\ref{eq:v1});\,  (\ref{eq:condrfy}) $\Rightarrow$ (\ref{eq:v1}).\ensep
On account of formula~(\ref{eq:rratesfrommargins}), conditions (\ref{eq:condrfx}) and (\ref{eq:condrfy}) are easily seen to be respectively equivalent to the following equalities:
\begin{alignat}{2}
\sum_{\latop{\scriptstyle y\in\ist}{\scriptstyle y\neq x}} \pmg_{xy}
\,\,&=\,\,
\sum_{\latop{\scriptstyle \bar x\in\xst}{\scriptstyle \bar x\neq x}}
\pmgbis_{x\bar x}
\,+\, |Y|,\qquad &&\hbox{for all $x\in\xst$},
\label{eq:condrfxbis}
\\[2.5pt]
\sum_{\latop{\scriptstyle x\in\ist}{\scriptstyle x\neq y}}\pmg_{yx}
\,\,&=\,\,
\sum_{\latop{\scriptstyle \bar y\in\yst}{\scriptstyle \bar y\neq y}}
\pmgbis_{y\bar y}
\,-\, |X|,\qquad &&\hbox{for all $y\in\yst$}.
\label{eq:condrfybis}
\end{alignat}
Let us add up respectively the equalities (\ref{eq:condrfxbis}) over $x\in\xst$ and the equalities (\ref{eq:condrfybis}) over $y\in\yst$. Since $\pmg_{pq}+\pmg_{qp} = \pmgbis_{pq}+\pmgbis_{qp} = 0$, we obtain
\begin{align}
\sum_{\latop{\scriptstyle x\in\xst}{\scriptstyle y\in\yst}} \pmg_{xy}
\,\,&=\,\,
|X|\,|Y|,
\label{eq:condrfxsum}
\\[2.5pt]
\sum_{\latop{\scriptstyle y\in\yst}{\scriptstyle x\in\xst}}\pmg_{yx}
\,\,&=\,\,
- |X|\,|Y|.
\label{eq:condrfysum}
\end{align}
Since the projected margins belong to $[-1,1]$, the preceding equalities imply respectively
\begin{alignat}{2}
\pmg_{xy} \,&=\, 1, \qquad&&\hbox{for all $x\in\xst$ and $y\in\yst$,}
\label{eq:condrfxall}
\\[2.5pt]
\pmg_{yx} \,&=\, -1, \qquad&&\hbox{for all $x\in\xst$ and $y\in\yst$,}
\label{eq:condrfyall}
\end{alignat}
(which are equivalent to each other since $\pmg_{xy}+\pmg_{yx}=0$).
Finally, either of these equalities implies that $\psc_{xy}=1$ for all $xy\in\xst\times\yst$, from which Lemma~\ref{lema1-vxy} allows to obtain (\ref{eq:v1}).
\end{proof}


\bigskip
The following propositions do not require the votes to be complete, but they require them to be rankings, or, more generally, in the case of Proposition~\ref{st:condrfmitjana}, to be transitive relations.

\smallskip
\begin{lemma}\hskip.5em
\label{lema2-vxy}
In the case of ranking votes, condition~\textup{(\ref{eq:v1})} implies that\linebreak $t_{xy}=1$ \,for any $x\in\xst$ and $y\in\ist$.
\end{lemma}
\begin{proof}\hskip.5em
In fact, even if we are dealing with truncated ranking votes, the rules that we are using for translating them into binary preferences ---namely, rules~(a--d) of \secpar{2.1}--- entail the following implications:\ensep (i)~$v_{xy}=1$ for some $y\in\ist$ implies that $x$~is explicitly mentioned in all of the ranking votes;\ensep and\,~(ii)~$x$~being explicitly mentioned in all of the ranking votes implies that $t_{xy}=1$ for any $y\in\ist$.
\end{proof}

\smallskip
\begin{proposition}
\label{st:condrfRV}
In the case of ranking votes,
condition~\textup{(\ref{eq:v1})} implies~\textup{(\ref{eq:condrfx})}.
\end{proposition}

\begin{proof}\hskip.5em
Let us fix an admissible order. According to Lemma~\ref{lema3-vxy},
we have $\pto_{x\bar x}=1$ for all $x,\bar x\in\xst$. On the other hand, Lemma~\ref{lema2-vxy} ensures that $t_{x\bar x}=1$ for all $x,\bar x\in\xst$, from which it follows that $\ptobis_{x\bar x}=1$ for all $x,\bar x\in\xst$\linebreak 
(since $\ptobis_{x\bar x}$ are the turnovers obtained from the restriction to the matrix $\vxx$, which belongs to the complete case). In particular, we have $\ppto_{xx'}=\pptobis_{xx'}=1$ whenever $x,x'\in\xst$.\ensep
On the other hand, Lemma~\ref{lema3-vxy} ensures also that $\ppmg_{xx'}=\ppmgbis_{xx'}$ whenever $x,x'\in\xst$.
\ensep
These equalities entail that 
$\psc_{x\bar x}=\pscbis_{x\bar x}$ for all $x,\bar x\in\xst$.
\ensep
By Lemma~\ref{lema1-vxy} we know also that $\psc_{xy}=1$ for all $xy\in\xst\times\yst$. 
Therefore,
\begin{equation*}
r_x \,=\,
N-\sum_{{\scriptstyle y\ne x \atop\scriptstyle y\in A }} \psc_{xy}
\,=\,
|X|-\sum_{{\scriptstyle \bar x\ne x  \atop\scriptstyle \bar x\in\xst
}} \pscbis_{x\bar x}
\,=\,
\rbis_x,\qquad \forall x\in\xst.
\end{equation*}
\end{proof}

\smallskip
\begin{proposition}\hskip.5em
\label{st:condrfmitjana}
Assume that the individual votes are transitive relations.
In this case, the equality
\begin{equation}
\sum_{x\in\xst}r_x \,=\, |\xst|(|\xst|+1)/2
\label{eq:rmitjana}
\end{equation}
implies \textup{(\ref{eq:v1})} (with $\yst=\ist\setminus\xst$).
\end{proposition}

\begin{proof}\hskip.5em
Let us introduce formula (\ref{eq:rrates}) for $r_x$ into (\ref{eq:rmitjana}). By using the fact that $\psc_{x\bar x} + \psc_{\bar xx} \le 1$, one obtains that
\begin{equation}
\sum_{\latop{\scriptstyle x\in\xst}{\scriptstyle y\in\yst}} \psc_{xy}
\,\,\ge\,\,
|X|\,|Y|.
\label{eq:condrfmitjanasum}
\end{equation}
The only possible way to satisfy this inequality is having 
$\psc_{xy}=1$ for all $xy\in\xst\times\yst$. Finally, (\ref{eq:v1}) follows by virtue of Lemma~\ref{lema1-vxy} since we are assuming that the individual votes are transitive relations. 
\end{proof}

\smallskip
\begin{corollary}\hskip.5em
\label{st:condrf1} Assume that the votes are rankings. Then $r_x=1$
if and only if \,all voters have put $x$ into first place.
\end{corollary}

\begin{proof}\hskip.5em
It suffices to apply Propositions~\ref{st:condrfRV} and \ref{st:condrfmitjana} with $\xst=\{x\}$.
\end{proof}

\begin{comment}
Si fos cert que en el cas de rankings $\sum_{\{x,y\}
}\pto_{xy}=\sum_{\{x,y\}} t_{xy}$ aleshores podr\'{\i}em demostrar que,
en aquest cas, $r_x=\tilde r_x$ for all $x\in\xst$ implica $v_{xy}=1$
per a tot $xy\in \xst\times\yst$: $$\sum_{x\in
X}r_x=N\,|X|-\sum_{\{x,\bar x\}\subset X}\pto_{x\bar x} -\sum_{x\in
X,y\in\yst} \psc_{xy}=\sum_{x\in\xst}\tilde r_x=|X|^2-\sum_{\{x,\bar
x\}\subset X} \pto_{x\bar x}$$ $$\Leftrightarrow |X|\,|Y|
+\sum_{x\in\xst,y\in\yst}\psc_{yx} -\sum_{x\in\xst,y\in\yst} t_{xy}= 0$$
$$\Leftrightarrow \psc_{yx}-t_{xy}=-1, \; \forall x\in\xst,y\in\yst$$
which implies $t_{xy}=1$ and $\psc_{yx}=0$.... Tampoc surt!!!
\end{comment}

\bigskip
The next theorem establishes property~\llfd.

\smallskip
\begin{theorem}\hskip.5em
\textup{(a)} In the complete case, or alternatively, under the hypothesis that the individual votes are rankings, one has the following implication: Assume that $\xst\sbset\ist$ has the property that $v_{xy}=1$ whenever $x\in\xst$ and $y\in\yst=\ist\setminus\xst$, and that there is no proper subset with the same property. In~that case, the fraction-like rates satisfy $\flr_x=\flrbisx>0$ for all $x\in\xst$ and $\flr_y=0$ for all $y\in\yst$.\ensep
\textup{(b)} In the complete case the converse implication holds too.
\end{theorem}

\begin{proof}\hskip.5em
Statement~(a).\ensep
Let us fix an admissible order $\xi$. By Lemma~\ref{lema1-vxy},
the hypothesis that $v_{xy}=1$ for all $xy\in\xst\times\yst$ implies
the following facts for all $xy\in\xst\times\yst$:\ensep $\img_{xy}=1$,\ensep $xy\in\xi$,\ensep $\psc_{xy}=1$,\ensep $\psc_{yx}=0$.\ensep
On the other hand, we can see that under the present hypothesis one has
\begin{equation}
\label{eq:pscxx}
\psc_{x\bar x}=\pscbis_{x\bar x},\qquad \hbox{for~any $x,\bar x\in\xst$.}
\end{equation}
In~the complete case this follows from Lemma~\ref{lema3-vxy}. Under the alternative hypothesis that the individual votes are rankings, it can be obtained  as in the proof of Proposition~\ref{st:condrfRV} as a consequence of Lemma~\ref{lema3-vxy} and the fact that in this case $\pto_{x\bar x} = \ptobis_{x\bar x} = 1$ for any $x,\bar x\in\xst$.

Now, according to Lemma~\ref{lema_irreductibilitat_projectada}, the matrix $(\psc_{xy})$ has a top dominant irreducible component $\xsth$. Since $\psc_{yx}=0$ for all $xy\in\xst\times\yst$, it is clear that $\xsth\sbseteq\xst$. However, a strict inclusion $\xsth\sbset\xst$ would imply $\psc_{x\xh}=0$ and therefore $\psc_{\xh x}=1$ for any $x\in\xst\setminus\xsth$ and $\xh\in \xsth$. Since we also have $\psc_{xy}=1$ for $x\in\xst$ and $y\not\in\xst$, we would get $\psc_{\xh\yh}=1$ for all $\xh\in\xsth$ and $\yh\not\in\xsth$, which would imply, by Lemma~\ref{lema1-vxy}, that $v_{\xh\yh}=1$ for all such pairs. This would contradict the supposed minimality of $\xst$. So, $\xst$ itself is the top dominant irreducible component of the matrix $(\psc_{xy})$.

By making use of  Theorem~\ref{st:zermelore}, it follows that $\flr_x=\flrbisx>0$ for all $x\in\xst$ and $\flr_y=0$ for all $y\in\yst$. In principle, $\flrbisx$ are here the fraction-like rates computed from the restriction of $(\psc_{xy})$ to the set $\xst$. However, (\ref{eq:pscxx}) allows to view them also as the fraction-like rates computed from the matrix~$(\pscbis_{xy})$, which by definition has been worked out from the restriction of $(v_{xy})$ to $x,y\in\xst$.

\smallskip
Statement~(b).\ensep 
Let us begin by noticing that the hypothesis that $\flr_x>0$ for all $x\in\xst$ and $\flr_y=0$ for all $y\in\yst=\ist\setminus\xst$ implies that $\xst$ is the~top dominant irreducible component of the matrix $(\psc_{xy})$. In fact, otherwise Theorem~\ref{st:zermelore} would imply the existence of some $x\in\xst$ with $\flr_x=0$ or~some $y\in\yst$ with $\flr_y>0$.\ensep
In particular, we have $\psc_{yx}=0$ for~all $xy\in\xst\times\yst$. Because of the completeness assumption, this implies that $\psc_{xy}=1$ and \hbox{---by Lemma~\ref{lema1-vxy}---} $v_{xy}=1$ for all those pairs.\ensep
Finally, let us see that $X$ is minimal for this property: If we had $\xsth\sbset\xst$ satisfying $v_{\xh\yh}=1$ for all $\xh\yh\in\xsth\times\ysth$ with $\ysth=\ist\setminus\xst$, then Lemma~\ref{lema1-vxy} would give $\psc_{\xh\yh}=1$ and therefore $\psc_{\yh\xh}=0$ for all such pairs, so $\xst$ could not be the top dominant irreducible
component of the matrix $(\psc_{xy})$.
\end{proof}

\section{The majority principle}

\smallskip
\begin{theorem}\hskip.5em
\label{st:majPrinciple}
The relation~$\mu(\isc)$ complies with the majority principle:
Let $\ist$ be partitioned in two sets $\xst$ and $\yst$ with the property that $v_{xy} > 1/2$ whenever $x\in\xst$ and $y\in\yst$; in that case, $\mu(\isc)$ includes any pair $xy$ with $x\in\xst$ and $y\in\yst$.
\end{theorem}

\begin{proof}\hskip.5em
Assume that $x\in\xst$ and $y\in\yst$.
Since $\isc_{xy} \ge v_{xy}$, the hypothesis of the theorem entails that $\isc_{xy} > 1/2$.
\ensep
On the other hand, let $\gamma$ be a path from $y$ to $x$ such that $\isc_{yx} = v_\gamma$; since it goes from $\yst$ to $\xst$, this path must contain at least one link $y_iy_{i+1}$ with $y_i\in\yst$ and $y_{i+1}\in\xst$; now, for this link we have $v_{y_iy_{i+1}} \le 1 - v_{y_{i+1}y_i} < 1/2$, which entails that $\isc_{yx} = v_\gamma < 1/2$.
\ensep
Therefore, we get $\isc_{yx} < 1/2 < \isc_{xy}$, \ie $xy\in\mu(\isc)$.
\end{proof}

\smallskip
\begin{corollary}\hskip.5em
\label{st:RvsNuCor}
The social ranking determined by the rank-like rates complies with the majority principle: Let $\ist$ be partitioned in two sets $\xst$ and $\yst$ with the property that $v_{xy} > 1/2$ whenever $x\in\xst$ and $y\in\yst$; in that case, the inequality $r_x < r_y$ holds for any $x\in\xst$ and $y\in\yst$.
\end{corollary}
\begin{proof}\hskip.5em
It follows from Theorem~\ref{st:majPrinciple} by virtue of part~(c) of Theorem~\ref{st:RvsNu}.
\end{proof}

\smallskip
\begin{corollary}\hskip.5em
\label{st:condorcetPrinciple}
In the complete case the social ranking determined by the rank-like rates complies with the Condorcet principle:
If $x$ has the property that $v_{xy} > v_{yx}$ for any $y\neq  x$, then $r_x < r_y$ for any $y\neq x$.
\end{corollary}
\begin{proof}\hskip.5em
In the complete case $v_{xy} > v_{yx}$ implies $v_{xy} > 1/2$. So, it suffices to apply the preceding result with $\xst=\{x\}$ and $\yst=A\setminus\xst$.
\end{proof}

\section{Clone consistency}

The notion of a cluster (of clones) was defined in \S{5} in connection with a binary relation: A subset $\cst\sbseteq\ist$ is said to be a cluster for a relation~$\rho$ when, for any~$x\not\in\cst$, having $ax\in\rho$ for some $a\in\cst$ implies $bx\in\rho$ for any $b\in\cst$, and similarly, having $xa\in\rho$ for some $a\in\cst$ implies $xb\in\rho$ for any $b\in\cst$.

Here we will extend the notion of a cluster in the following way:
$\cst\sbseteq\ist$ is said to be a cluster for a system of binary scores~$(v_{xy})$ when
\begin{equation}
v_{ax} = v_{bx},\quad
v_{xa} = v_{xb},\qquad
\text{whenever $a,b\in\cst$ and $x\not\in\cst$.}
\end{equation}
This definition can be viewed as an extension of the preceding one because of the following obvious fact:
\begin{lemma}\hskip.5em
\label{st:clons0}
$\cst$ is a cluster for a relation $\rho$ \ifoi $\cst$ is a cluster for the corresponding system of binary scores, which is defined in~$(\ref{eq:binrelmatrix})$.
\end{lemma}
\noindent
In particular, the extended notion allows the following results to include the case where the individual votes belong to the general class considered in~\secpar{3.3}.

\bigskip
In this section we will prove the clone consistency property~\llcc:
If a set of options is a~cluster for each of the individual votes, then: (a)~it is a cluster for the social ranking;\, and\, (b)~contracting it to a single option in all of the individual votes has no other effect in the social ranking than getting the same contraction.

\bigskip
In the remainder of this section we assume the following \textbf{standing\linebreak 
hypothesis}:
$$
\text{\textit{$\cst$~is a cluster for all of the individual votes.}}
$$
Since the collective binary scores are obtained by adding up the individual ones (equation~(\ref{eq:cog})), the preceding hypothesis immediately implies that
$$
\text{\textit{$\cst$~is a cluster for the collective binary scores $v_{xy}$.}}
$$
In the following we will see that this property of being a cluster is
maintained throughout the whole procedure which defines the social ranking.

\bigskip
\begin{lemma}\hskip.15em
\label{st:clons2}
Assume that either $x$ or $y$, or both, lie outside $\cst$. In this case\linebreak
\begin{equation*}
\isc_{xy} \,=\, \max\,\{\,v_\gamma\mid\gamma\text{ contains no more than one element of }\cst\,\}
\end{equation*}
\end{lemma}
\begin{proof}\hskip.5em
It suffices to see that any path $\gamma=x_0 \dots x_n$ from $x_0=x$ to~$x_n=y$ which contains more than one element of~$\cst$ can be replaced by~another one~$\gammabis$ which contains only one such element and satisfies $v_{\gammabis}\ge v_\gamma$. Consider first the case where $x,y\not\in\cst$. In this case it will suffice to take $\gammabis=x_0\dots x_{j-1}x_k\dots x_n$, where $j = \min\,\{\,i\mid x_i\in\cst\,\}$ and $k = \max\,\{\,i\mid x_i\in\cst\,\}$, which obviously satisfy $0<j<k<n$. Since $x_{j-1}\not\in\cst$ and $x_j,x_k\in\cst$, we~have $v_{x_{j-1}x_j}=v_{x_{j-1}x_k}$, so that
\begin{equation*}
\begin{split}
v_\gamma
\,&=\, \min\,\left(v_{x_0x_1},\dots,v_{x_{n-1}x_n}\right)\\
\,&\le\, \min\,\left(v_{x_0x_1},\dots,v_{x_{j-1}x_j},
v_{x_k x_{k+1}},\dots,v_{x_{n-1}x_n}\right)\\
\,&=\, \min\,\left(v_{x_0x_1},\dots,v_{x_{j-1}x_k},
v_{x_k x_{k+1}},\dots,v_{x_{n-1}x_n}\right)\\
\,&=\, v_{\gammabis}.
\end{split}
\end{equation*}
The case where $x\not\in\cst$ but $y\in\cst$ can be dealt with in a similar way by~taking $\gammabis=x_0\dots x_{j-1}x_n$, and analogously, in the case where $x\in\cst$ and $y\not\in\cst$ it~suffices to take $\gammabis=x_0x_{k+1}\dots x_n$.
\end{proof}

\smallskip
\begin{proposition}\hskip.5em
\label{st:clons3}
$\cst$~is a cluster for the indirect scores $\isc_{xy}$.
\end{proposition}
\begin{proof}\hskip.5em
Consider $a,b\in\cst$ and $x\not\in\cst$. Let $\gamma=x_0x_1x_2 \dots x_n$ be a path from $a$ to $x$ such that $\isc_{ax}=v_\gamma$. By Lemma~\ref{st:clons2}, we can assume that $a$ is the only element of $\gamma$ that belongs to $\cst$. In particular, $x_1\not\in\cst$, so that $v_{ax_1} = v_{bx_1}$, which allows to write
\begin{equation*}
\begin{split}
\isc_{ax} \,=\, v_\gamma
\,&=\, \min\,\left(v_{ax_1},v_{x_1x_2},\dots,v_{x_{n-1}x}\right)\\
\,&=\, \min\,\left(v_{bx_1},v_{x_1x_2},\dots,v_{x_{n-1}x}\right)\\
\,&\le\, \isc_{bx}.
\end{split}
\end{equation*}
By~interchanging $a$ and $b$, one gets the reverse inequality $\isc_{bx}\le\isc_{ax}$ and there\-fore the equality $\isc_{ax}\cd=\isc_{bx}$. An analogous argument shows that $\isc_{xa}\cd=\isc_{xb}$.
\end{proof}

\smallskip
\begin{proposition}
\label{st:clonsNu}
\!$\cst$\!~is a cluster for the indirect comparison relation $\nu\!=\!\mu(\isc)$.
\end{proposition}
\begin{proof}
This is an immediate consequence of the preceding proposition.
\end{proof}

\medskip
\begin{proposition}\hskip.5em
\label{st:clonsXi}
There exists an admissible order $\xi$
such that $C$ is a~cluster for $\xi$.
\end{proposition}
\begin{proof}
This result is given by Theorem~\ref{st:existenceXiCluster} of p.\,\pageref{st:existenceXiCluster}.
\end{proof}

\begin{comment}
Pending:
Is $\cst$ a cluster for the projected scores $\psc_{xy}$?
\end{comment}

\smallskip
\begin{theorem} 
\label{st:clons4}
$\cst$~is a cluster for the ranking defined by the rank-like rates \textup{(}\ie for the relation $\sigma = \{xy\in\tie\mid r_x<r_y\}$\textup{)}.
\end{theorem}
\begin{proof}\hskip.5em
We must show that, for any $x\not\in\cst$ and any $a,b\in\cst$,
$r_a<r_x$ implies $r_b<r_x$
\,and\, $r_x<r_a$ implies $r_x<r_b$
(from which it follows that $r_a=r_x$ implies $r_b=r_x$).
\ensep
Equivalently, it suffices to show that:\ensep
(a)~$r_a<r_x$ implies $r_b\le r_x$;\ensep
(b)~$r_x<r_a$ implies $r_x\le r_b$;\ensep
and\, (c)~$r_a=r_x$ implies $r_b=r_x$.
\ensep
The~proof will make use of an admissible order~$\xi$
with the property that $\cst$~is a cluster for~$\xi$
(whose existence is ensured by Proposition~\ref{st:clonsXi}).

Parts~(a) and (b) are then a straightforward consequence of
part~(a) of Lemma~\ref{st:obsRosa}.(a).
In fact, by combining this result, and its contrapositive,
with the fact that $\cst$~is a cluster for~$\xi$, we~have the following implications: $r_a<r_x \Rightarrow ax\in\xi \Rightarrow bx\in\xi \Rightarrow r_b\le r_x$, and similarly, $r_x<r_a \Rightarrow xa\in\xi \Rightarrow xb\in\xi \Rightarrow r_x\le r_b$.

Part~(c):\,~$r_a=r_x$ implies $r_b=r_x$ (for $x\not\in\cst$ and $a,b\in\cst$). Since $\xi$~is a total order, we must have either $ax\in\xi$ or $xa\in\xi$; in the following we assume $ax\in\xi$ (the other possibility admits of a similar treatment). In order to deal with this case we will consider the last element of $\cst$ according to $\xi$, which we will denote as $\ell$, and its immediate successor $\ell'$, which does not belong to $\cst$. Since $a\rxieq \ell\rxi \ell'\rxieq x$ 
and $r_a = r_x$, we must have $r_\ell = r_{\ell'}$.
Now, according to part~(b) of Lemma~\ref{st:obsRosa}, $\pmg_{\ell\ell'}=0$; in other words, $\ppmg_{\ell\ell'}=0$. By~the definition of 
$\ppmg_{\ell\ell'}$, this means that there exist $p$ and $q$ with $p\rxieq \ell \rxi q$ such that $\img_{pq}=\isc_{pq}-\isc_{qp}=0.$
Obviously, $q\not\in\cst$, whereas $p$ either belongs to $\cst$ or it precedes all elements of $\cst$. In the latter case, we immediately get $\ppmg_{cc'}=0$ for all $c\in\cst$ (by the definition of~$\ppmg_{cc'}$). If $p\in\cst$, we arrive at the same conclusion thanks to 
Proposition~\ref{st:clons3}, which ensures that $\img_{cq}=\img_{pq}$. 
So,~the intervals $\gamma_{cc'}$ with $c\in C$ are all of them reduced to a point. Since $\cst$ is a cluster for the total order~$\xi$, this implies that $\gamma_{ab}$ is also reduced to the same point (this holds for any $a,b\in\cst$). According to part~(b) of Lemma~\ref{st:obsRosa}, this implies that $r_a=r_b$, as it was claimed.
\end{proof}

\medskip
Finally, we consider the effect of contracting $\cst$ to a single element.\linebreak 
So~we consider a new set~$\istbis = (\ist\setminus\cst)\cup\{\clustit\}$ together with the scores 
$\vbis_{xy}\ (x,y\in\istbis)$
defined by the following equalities, where $p,q\in\ist\setminus\cst$ and $\clone$~is an arbitrary element of~$\cst$:\ensep
$\vbis_{pq} = v_{pq}$, $\vbis_{p\clustit} = v_{p\clone}$ \,and\, 
$\vbis_{\clustit q} = v_{\clone q}$\linebreak 
(the definition is not ambiguous since $\cst$ is a cluster for the scores~$v_{xy}$).
In the following, a tilde is systematically used to distinguish between hom\-olo\-gous objects associated respectively with $(\ist,v)$ and~$(\istbis,\vbis)$. We will also make~use of the following notation: for every $x\in\ist$, $\contr{x}$ denotes the element of~$\istbis$ defined by $\contr{x} = \clustit$ if~$x\in\cst$ and by $\contr{x} = x$ if~$x\not\in\cst$; in~terms of this mapping, the~preceding equalities say simply that $\vbis_{\contr{x}\contr{y}}=v_{xy}$ whenever $\contr{x}\neq \contr{y}$.

\smallskip
\begin{theorem}\hskip.5em
\label{st:clons5}
The ranking $\widetilde\sigma = \{xy\in\tiebis\mid \rbis_x<\rbis_y\}$ coincides with the contraction of $\sigma = \{xy\in\tie\mid r_x<r_y\}$ by the cluster $\cst$.
\end{theorem}
\begin{proof}\hskip.5em
We begin by noticing that the indirect scores $\iscbis_{xy}\ (x,y\in\istbis)$ coincide with those obtained by contraction of the $\isc_{xy}\ (x,y\in\ist)$, \ie $\iscbis_{\contr{x}\contr{y}}=\isc_{xy}$ whenever $\contr{x}\neq \contr{y}$. 
This follows from the analogous equality between the direct scores because of Lemma~\ref{st:clons2}. As a consequence, $\nubis=\mu(\iscbis)$ coincides with the contraction of $\nu=\mu(\isc)$ by $\cst$. From this fact, parts~(a) and~(b) of~Theorem~\ref{st:RvsNu}, allow  to derive that $r_x < r_y$ implies $\rbis_{\contr{x}} \le \rbis_{\contr{y}}$ whenever $\contr{x}\neq \contr{y}$, and that $\rbis_{\contr{x}} < \rbis_{\contr{y}}$ implies $r_x \le r_y$.

In order to~complete the proof, we~must check that $r_x = r_y$ is equivalent to~$\rbis_{\contr{x}} = \rbis_{\contr{y}}$ whenever $\contr{x}\neq \contr{y}$. According to part~(b) of Lemma~\ref{st:obsRosa}, it~suffices to see that $\pmg_{xy} = 0$ is equivalent to~$\pmgbis_{\contr{x}\contr{y}} = 0$ whenever $\contr{x}\neq \contr{y}$. In order to prove this equivalence, we need to look at the way that $\pmg_{xy}$ and~$\pmgbis_{\contr{x}\contr{y}}$ are obtained, which requires  certain admissible orders $\xi$ and $\xibis$; in~this connection, it will be useful that $\xi$ be one of the admissible orders for which $\cst$~is a~cluster (Proposition~\ref{st:clonsXi}), and that $\xibis$ be the corresponding contraction, which is admissible as a consequence of Proposition~\ref{st:clons3}. Now, that proposition entails not only that $\cst$ is a cluster for the indirect margins $\img_{pq}$, but also that their contraction by $\cst$ coincides with the margins of the contracted indirect scores, \ie $\imgbis_{\contr{p}\contr{q}}=\img_{pq}$ whenever $\contr{p}\neq \contr{q}$. Moreover, 
by the definition of~the intermediate projected margins, namely equation~(\ref{eq:projection2}.1), it follows that $\cst$ is~also a cluster for the intermediate projected margins $\ppmg_{pq}$ and that their contraction by $\cst$ coincides with the homologous quantities obtained from the contracted indirect margins, \ie $\ppmgbis_{\contr{p}\contr{q}}=\ppmg_{pq}$ whenever $\contr{p}\neq \contr{q}$. On~the other hand, it is also clear from equation~(\ref{eq:projection2}.1) that the intermediate projected margins behave in the following way:
\begin{equation}
\ppmg_{pq} \,\le\, \ppmg_{ab}\qquad\hbox{ whenever }a\rxieq p\rxi q\rxieq b.
\label{eq:ppmgbehaviour}
\end{equation}

After these remarks, we proceed with 
showing that $\pmg_{xy}=0$ is equivalent to $\pmgbis_{\contr{x}\contr{y}}=0$ whenever $\contr{x}\neq \contr{y}$.
By symmetry, we~can assume that $xy\in\xi$, which entails that $\contr{x}\contr{y}\in\xibis$.
In view of (\ref{eq:projection3}--\ref{eq:projection5}),
the equality $\pmg_{xy}=0$ is equivalent to saying that $\ppmg_{hh'}=0$ for all $h$ such that $x\rxieq h\rxi y$, and~similarly, the  equality $\pmgbis_{\contr{x}\contr{y}}=0$ is~equivalent to $\widetilde m^\sigma_{\eta\eta'}=0$ for all $\eta$ such that $\contr{x}\rxieq\eta\rxi\contr{y}$. By considering a path $x_0x_1\dots x_n$ from $x_0=x$ to $x_n=y$ with $x_ix_{i+1}$ consecutive in $\xi$, it is clear that the problem reduces to proving the following implications, where $\ell$ denotes the last element of $\cst$ by $\xi$, $f$~denotes the first one, and $\anteh$ denotes the element that immediately precedes $h$ in $\xi$:
\ensep
(a)~$\ppmg_{\ell\ell'}=0 \,\Rightarrow\, \ppmgbis_{\clustit\ell'}=0$;
\ensep
(b)~$\ppmgbis_{\clustit\ell'}=0 \,\Rightarrow\, \ppmg_{cc'}=0$ for any $c\in\cst$;
\ensep
(c)~$\ppmg_{\antef f}=0 \,\Rightarrow\, \ppmgbis_{\antef\clustit}=0$;
\ensep and\, 
(d)~$\ppmgbis_{\antef\clustit}=0 \,\Rightarrow\, \ppmg_{\antec c}=0$ for any $c\in\cst$.\linebreak 
Now, (a) and (c) are immediate consequences of the fact that $\ppmgbis_{\contr{p}\contr{q}}=\ppmg_{pq}$ whenever $\contr{p}\neq \contr{q}$.
On the other hand, (b) and (d)~follow from the same equality together with the inequality (\ref{eq:ppmgbehaviour}). In fact, these facts allow us to write $\ppmg_{cc'} \le \ppmg_{c\ell'} = \ppmgbis_{\clustit\ell'}$, which gives~(b), and similarly, $\ppmg_{\antec c} \le \ppmg_{\antef c} = \ppmgbis_{\antef\clustit}$, which gives~(d).
\end{proof}

\begin{comment}
Is $\cst$ a cluster for the fraction-like rates?
\end{comment}

\section{Approval voting}

In approval voting, each voter is asked for a list of approved options, without any expression of preference between them, and each option~$x$ is then rated by the number of approvals for it~\cite{br}. 
In the following we will refer to this number as the \dfc{approval score} of~$x$, and its value relative to $V$ will be denoted by $\av{x}$.

From the point of view of paired comparisons, an individual vote of approval type can be viewed as a truncated ranking where all of the options that appear in it are tied.
In~this section, we will see that the margin-based variant
orders the options exactly in the same way as the approval scores.
In~other words, the method of approval voting agrees with ours under interpretation~(d$'$) of~\secpar{3.2},
\ie under the interpretation that the non-approved options of each individual vote are tied.

Having said that, the preliminary results~\ref{st:av1}--\ref{st:av3} will hold not only under interpretation~(d$'$) but also under interpretation~(d), \ie that there is no information about the preference of the voter between two non-approved options, and also under the analogous interpretation that there is no information about his preference between two approved options.\ensep
Interpretation~(d$'$) does not play an essential role until Theorem~\ref{st:av4}, where we use the fact that it always brings the problem into the complete case.

In the following, $\lambda(\av{})$ denotes the relation defined by
\begin{equation}
xy\in\lambda(\alpha)
\quad\equiv\quad \alpha_x > \alpha_y.
\end{equation}

\smallskip
\begin{proposition}\hskip.5em
\label{st:av1}
In the approval voting situation, 
the following equality holds:
\begin{equation}
v_{xy} - v_{yx} \,=\, \av{x} - \av{y}.
\label{eq:mgav}
\end{equation}
In particular, $\mu(v)=\lambda(\av{})$.
\end{proposition}
\begin{proof}\hskip.5em
Obviously, the possible ballots are in one-to-one correspondence with the subsets $\xst$ of $\ist$. In the following, $v_\xst$ denotes the relative number of votes that approved exactly the set $\xst$. 
With this notation it is obvious that
\begin{equation}
\av{x} \,=\, \sum_{\xst\ni\,x} v_{\xst} \,=\,
\sum_{\latop{\scriptstyle \xst\ni\,x}{\scriptstyle \xst\not\ni\,y}}
v_{\xst} +
\sum_{\latop{\scriptstyle \xst\ni\,x}{\scriptstyle \xst\ni\,y}}
v_{\xst}.
\end{equation}
On the other hand, one has
\begin{equation}
v_{xy} \,=\, 
\sum_{\latop{\scriptstyle \xst\ni\,x}{\scriptstyle \xst\not\ni\,y}}
v_{\xst} \ \bigg( +
{\textstyle\frac12}
\,\sum_{\latop{\scriptstyle \xst\ni\,x}{\scriptstyle \xst\ni\,y}}
v_{\xst} +
{\textstyle\frac12}
\,\sum_{\latop{\scriptstyle \xst\not\ni\,x}{\scriptstyle \xst\not\ni\,y}}
v_{\xst}\bigg),
\label{eq:vxyav}
\end{equation}
where the terms in brackets are present or not depending on which interpretation is used.
Anyway, the preceding expressions, together with the analogous ones where $x$ and $y$ are interchanged with each other, result in the equality~(\ref{eq:mgav})
independently of those alternative interpretations.
\end{proof}

\smallskip
\begin{corollary}\hskip.5em
\label{st:av2}
In the approval voting situation,
a path $x_0x_1\dots x_n$ is contained in $\mu(v)$ (resp.\ $\hat\mu(v)$) \,\ifoi\, the~sequence $\av{x_i}\ (i=0,1,\dots n)$ is decreasing (resp.\ non-increasing).
\end{corollary}

\smallskip
\begin{proposition}\hskip.5em
\label{st:av3}
In the approval voting situation, one has 
$\mu(\ww^\ast )=\lambda(\av{})$.
\end{proposition}

\begin{proof}\hskip.5em
Let us begin by proving that
\begin{equation}
\av{x} > \av{y}
\ensep\Longrightarrow\ensep
\ww^\ast_{xy}> \ww^\ast_{yx}.
\end{equation}
We will argue by contradiction. So, let us assume that $\ww^\ast_{yx}\ge \ww^\ast_{xy}$. According to Proposition~\ref{st:av1}, the hypothesis that $\av{x}>\av{y}$ is equivalent to $v_{xy}> v_{yx}$, which entails that $w_{xy}>0$ (by the definition of $w_{xy}$ together with the strict inequality $v_{xy} >  v_{yx}$). Now, since $\ww^\ast_{xy}\ge w_{xy}$ and we are assuming that $\ww^\ast_{yx}\ge \ww^\ast_{xy}$, it follows that $\ww^\ast_{yx}> 0$. This implies the existence of a path from $y$ to $x$ which is contained in $\hat\mu(v)$ (by the definitions of $\ww^\ast_{yx}$ and $w_{pq}$). Finally, Corollary~\ref{st:av2} produces a contradiction with the present hypothesis that $\av{x}>\av{y}$.

\halfsmallskip
Let us see now that
\begin{equation}
\av{x}=\av{y} \ensep\Longrightarrow\ensep \ww^\ast_{xy} = \ww^\ast_{yx}.
\end{equation}
Again, we will argue by contradiction. So, let us assume that $\ww^\ast_{xy}\neq  \ww^\ast_{yx}$. Obviously, it suffices to consider the case $\ww^\ast_{xy}> \ww^\ast_{yx}$. Now, this inequality implies that $\ww^\ast_{xy}> 0$, which tells us that $\ww^\ast_{xy}=w_\gamma$ for a certain path $\gamma:x_0x_1\dots x_n$ which goes from $x_0 = x$ to $x_n = y$ and is contained in $\hat\mu(v)$. According to~Corollary~\ref{st:av2}, we are ensured that the~sequence $\av{x_i}\ (i=0,1,\dots n)$ is non-increasing. However, the hypothesis that $\av{x}=\av{y}$ leaves no other possibility than $\av{x_i}$ being constant. So, the reverse path $\gamma':x_nx_{n-1}\dots x_1x_0$ is also contained in $\hat\mu(v)$. Besides, Proposition~\ref{st:av1} ensures that $v_{x_{i+1}x_i}=v_{x_ix_{i+1}}$, so~that $w_{\gamma'}=w_\gamma$. Since $\ww^\ast_{yx}\ge w_{\gamma'}$, it follows that $\ww^\ast_{yx}\ge \ww^\ast_{xy}$, which contradicts the hypothesis that $\ww^\ast_{xy}> \ww^\ast_{yx}$.

\halfsmallskip
Finally, one easily checks that the preceding implications entail that  $\sgn\,(\av{x}-\av{y})$ is always equal to $\sgn\,(\ww^\ast_{xy}-\ww^\ast_{yx})$. This is equivalent to the equality of the relations $\lambda(\av{})$ and $\mu(\ww^\ast)$.
\end{proof}

\smallskip
\begin{theorem}\hskip.5em
\label{st:av4}
In the approval voting situation, the margin-based variant results in a full compatibility relation between the rank-like rates $r_x$ and the approval scores $\alpha_x$: \,$r_x < r_y \,\Leftrightarrow\, \alpha_x > \alpha_y$.
\end{theorem}
\begin{proof}\hskip.5em
Recall that the margin-based variant amounts to using interpretation~(d$'$), which always brings the problem into the complete case (when the terms in brackets are included, equation~(\ref{eq:vxyav}) has indeed the property that $v_{xy}+v_{yx}=1$). So we can invoke Theorem~\ref{st:laia}. By combining it with\linebreak[3] Proposition~\ref{st:av3} we see that the inequality $\alpha_x > \alpha_y$ is equivalent to saying that $xy\in\nu$. 
In the following we will keep this equivalence in mind.
The implication $r_x < r_y \,\Rightarrow\, \alpha_x > \alpha_y$ is then an immediate consequence of part~(b) of Theorem~\ref{st:RvsNu}.
\ensep
The converse implication $\alpha_x > \alpha_y \,\Rightarrow\, r_x < r_y$ can be proved in the following way: Let $\xi$ be an admissible order. By definition, it contains $\nu$. So, the inequality $\alpha_x > \alpha_y$ implies $xy\in\xi$. On the other hand, that inequality implies also the existence of a consecutive pair $hh'$ with $x\rxieq h$ and $h'\rxieq y$ such that $\alpha_h > \alpha_{h'}$. As~a~consequence, one has $\alpha_p > \alpha_q$ whenever $p\rxieq h$ and $h'\rxieq q$. So, the sets $\xst=\{p\mid p\rxieq h\}$ and $\yst=\{q\mid h'\rxieq q\}$ are in the hypotheses of part~(c) of~Theorem~\ref{st:RvsNu}, which ensures the desired inequality $r_x < r_y$.
\end{proof}

\vskip-5pt
\remark
So in this case we get a converse of Theorem~\ref{st:RvsNu}.(b).
By following the same arguments as in the preceding proof, one can see that such a converse holds whenever there exists a function $s:\ist\ni x\mapsto s_x\in\bbr$, such that $xy\in\nu \,\Leftrightarrow\, s_x<s_y$.

\vskip5pt
\bigskip
Summing up, the standard approval voting procedure is always in full agreement with the margin-based variant of the CLC~method. In the approval voting situation, this variant amounts to treat all of the candidates which are missing in an approval ballot as equally ‘unpreferred’ (in the same way that all approved candidates are treated as equally preferred). This is quite reasonable if one can assume that the voters are well acquainted with all of the options.

\section{About monotonicity}

In this section we consider the effect of raising a particular option \,$a$\, to a more preferred status in the individual ballots without any change in the preferences about the other options. More generally, we consider the case where the scores~$v_{xy}$ are modified into new values $\vbis_{xy}$ such that
\begin{equation}
\label{eq:mona}
\vbis_{ay} \ge v_{ay},\quad \vbis_{xa} \le v_{xa},\quad \vbis_{xy} = v_{xy},\qquad \forall x,y\neq a.
\end{equation}
In such a situation, one would expect the social rates to behave in the following way, where $y$ is an arbitrary element of $\ist\setminus\{a\}$:
\begin{gather}
\label{eq:rmona}
\rbis_a < r_a,
\\[2.5pt]
\label{eq:rrmona}
r_a < r_y \,\Longrightarrow\, \rbis_a < \rbis_y,\qquad
r_a \le r_y \,\Longrightarrow\, \rbis_a \le \rbis_y,
\end{gather}
where the tilde indicates the objects associated with the modified scores.
Unfortunately, the rating method proposed in this paper does not satisfy these conditions, but generally speaking it satisfies only the following weaker ones:
\begin{gather}
\label{eq:wmona}
r_a < r_y\ \,\Longrightarrow\,\ \rbis_a \le \rbis_y.
\\[2.5pt]
\label{eq:smona}
(r_a < r_y,\ \forall y\neq a)\ \,\Longrightarrow\,\ 
(\rbis_a < \rbis_y,\ \forall y\neq a).
\end{gather}
In particular, (\ref{eq:smona})
is saying that\ensep if $a$ was the only winner for the scores $v_{xy}$, then it is still the only winner for the scores $\vbis_{xy}$.

Let us remark that in the case of ranking votes, situation~(\ref{eq:mona}) includes the following ones: (a)~the option~$a$ is raised to a better position in some of the ranking votes without any change in the preferences between the other options; (b)~the option~$a$ is appended to some ballots which did not previously contain~it; (c)~some ballots are added which plump for option~$a$.
\ensep
However, the third part of~(\ref{eq:mona}) leaves out certain situations which are sometimes considered the matter of other ``monotonicity'' conditions \cite{wo}.

In the terminology of \cite{bo}, property (\ref{eq:wmona}) is saying 
that the method that we are using is ``very weakly monotonic'' as a ranking procedure, whereas property (\ref{eq:smona}) is related to what \cite{bo} calls \hbox{``proper monotonicity''} of a choice procedure. In this connection, it is interesting to remark that the 
method of ranked pairs enjoys the choice\,-\,monotonicity property (\ref{eq:smona}) \cite[p.\,221--222]{t6}, but it lacks the ranking\,-\,monotonicity property (\ref{eq:wmona}). A~profile which exhibits such a failure of the ranking\,-\,monotonicity for the method of ranked pairs is given in \ensep {\small\texttt{http://mat.uab.cat/{\atilde}xmora/CLC\underline{ }calculator/}} (number~9 of ``Example inputs'').

\paragraph{18.1}This section is devoted to giving a proof of properties~(\ref{eq:wmona}) and~(\ref{eq:smona}).

\begin{theorem}\hskip.5em
\label{st:mono}
Assume that $(v_{xy})$ and $(\vbis_{xy})$ are related to each other 
in accordance with \textup{(\ref{eq:mona})}.
In this case, the following properties are satisfied for any $x,y\neq a$:
\begin{gather}
\iscbis_{ay} \ge \isc_{ay},\qquad
\iscbis_{xa} \le \isc_{xa},
\label{eq:mono1}
\\[2.5pt]
\pre{a}(\nubis) \sbseteq \pre{a}(\nu),\qquad
\suc{a}(\nubis) \spseteq \suc{a}(\nu),
\label{eq:mono2}
\\[2.5pt]
\begin{repeated}{eq:wmona}
r_a < r_y \,\Longrightarrow\, \rbis_a \le \rbis_y,
\end{repeated}
\label{eq:mono3}
\end{gather}
where $\nu=\mu(\isc)$ and $\nubis=\mu(\vbis^\ast)$
\end{theorem}

\begin{proof}\hskip.5em
Let us begin by seeing that (\ref{eq:mono3}) will be a consequence of (\ref{eq:mono2}).
In~fact, we have the following chain of implications: $r_a < r_y \,\Rightarrow\, y\in\suc{a}(\nu) \,\Rightarrow\, y\in\suc{a}(\nubis) \,\Rightarrow\, \rbis_a \le \rbis_y$, where the central one is provided by (\ref{eq:mono2}.2) and the other two are guaranteed by Theorem~\ref{st:RvsNu}.

The proof of (\ref{eq:mono1}--\ref{eq:mono2}) is organized in three steps.
In the first one, we~look at the special case where one increases the score of a single pair~$ab$.
After this, we will consider the case where an increase in the score of $ab$ is combined with a decrease in the score of $ba$. Finally, the third step deals with the general situation~(\ref{eq:mona}).

\medskip\noindent
\textsl{Special case~1.} \textit{Assume that}
\begin{equation}
\label{eq:hcas1}
\vbis_{ab}>v_{ab},\quad \vbis_{xy} = v_{xy},\qquad \forall\,xy\neq ab.
\end{equation}
\emph{In this case, the following properties are satisfied:}
\begin{alignat}{2}
\iscbis_{xy}&\ge \isc_{xy},\qquad &\forall x&,y
\label{eq:monot-cas1-1}
\\[2.5pt]
\iscbis_{xa}&=\isc_{xa},\qquad &\forall x&\ne a
\label{eq:monot-cas1-2}
\\[2.5pt]
\iscbis_{by}&=\isc_{by},\qquad &\forall y&\ne b
\label{eq:monot-cas1-3}
\\[2.5pt]
\hbox to0pt{\hss\small$(\ref{eq:mono2})$\hss}\hskip29mm 
\pre{a}(\nubis) &\sbseteq \pre{a}(\nu),\qquad
&\suc{a}(\nubis) &\spseteq \suc{a}(\nu),
\hskip29mm 
\label{eq:mono2R}
\\[2.5pt]
\pre{b}(\nubis) &\spseteq \pre{b}(\nu),\qquad
&\suc{b}(\nubis) &\sbseteq \suc{b}(\nu).
\label{eq:monot-cas1-4}
\end{alignat}

\halfsmallskip
In fact, under the hypothesis~(\ref{eq:hcas1}) it is obvious that $\vbis_\gamma\ge v_\gamma$ and that the strict inequality happens only when the path $\gamma=x_0\dots x_n$ contains the pair~$ab$ and the latter realizes the minimum of the scores $v_{x_ix_{i+1}}$. As a consequence, the indirect scores satisfy the inequality~(\ref{eq:monot-cas1-1}). Furthermore, a strict inequality in~(\ref{eq:monot-cas1-1}) implies that the maximum which defines $\iscbis_{xy}$ is realized by a path $\gamma$ which satisfies $\vbis_\gamma>v_\gamma$ and therefore contains the pair~$ab$.

Now, in order to obtain the indirect score for a pair of the form $xa$ it is useless to consider paths involving $ab$, since such paths contains cycles whose omission results in paths not involving $ab$ and having a better or equal score. So, the maximum which defines $\iscbis_{xa}$ is realized by a path which does not involve $ab$. According to
the last statement of the preceding paragraph, this implies~(\ref{eq:monot-cas1-2}). An entirely analogous argument establishes~(\ref{eq:monot-cas1-3}).

Finally, (\ref{eq:mono2R}) is obtained in the following way: $x\in\pre{a}(\nubis)$ means that $\iscbis_{xa}>\iscbis_{ax}$, from which (\ref{eq:monot-cas1-2}) and (\ref{eq:monot-cas1-1}) allow to derive that $\isc_{xa}=\iscbis_{xa}> \iscbis_{ax}\ge \isc_{ax}$, \ie $x\in\pre{a}(\nu)$.
\ensep
Similarly, $x\in\suc{a}(\nu)$ implies $x\in\suc{a}(\nubis)$ because one has
$\iscbis_{ax}\ge\isc_{ax}>\isc_{xa}=\iscbis_{xa}$.
\ensep
An analogous argument establishes~(\ref{eq:monot-cas1-4}).

\medskip\noindent
\textsl{Special case~2.} \textit{Properties \textup{(\ref{eq:mono1}--\ref{eq:mono2})} are satisfied in the following situation:}
\begin{equation}
\label{eq:hstep2}
\vbis_{ab}\ge v_{ab},\quad \vbis_{ba}\le v_{ba},\quad \vbis_{xy} = v_{xy},\qquad \forall\,xy\neq ab,ba.
\end{equation}

\pagebreak 

\halfsmallskip
This result will be obtained from the preceding one by going through an intermediate Llull matrix $\vbisbis$ defined in the following way
\begin{equation}
\begin{cases}
\vbisbis_{ab}=v_{ab}\\
\vbisbis_{ba}=\vbis_{ba}\\
\vbisbis_{xy}=\vbis_{xy}=v_{xy}, &\forall\,xy\ne ab,ba\\
\end{cases}
\end{equation}

If the hypothesis $\vbis_{ab}\ge v_{ab}$ is satisfied with strict inequality, then $\vbis$ and~$\vbisbis$ are in the hypotheses of the special case~1 (they play respectively the roles of $\vbis$ and $v$). In particular, we get
\begin{equation}\label{eq:lema-monotonia-intermitja-1}
\iscbis_{xy}\ge \iscbisbis_{xy},\quad
\iscbis_{xa}=\iscbisbis_{xa},\quad
\pre{a}(\nubis) \sbseteq \pre{a}(\nubisbis),\quad
\suc{a}(\nubis) \spseteq \suc{a}(\nubisbis).
\end{equation}
On the other hand, if $\vbis_{ab}=v_{ab}$ then $\vbisbis=\vbis$ and the preceding relations hold as equalities. 

Similarly, if the hypothesis $\vbis_{ba}\le v_{ba}$ is satisfied with strict inequality, then $v$ and~$\vbisbis$ are in the hypotheses of the special case~1 with $ab$ replaced by~$ba$ (they play respectively the roles of $\vbis$ and $v$). In particular, we get
\begin{equation}\label{eq:lema-monotonia-intermitja-2}
\isc_{xy}\ge \iscbisbis_{xy},\quad
\isc_{ay}=\iscbisbis_{ay},\quad
\pre{a}(\nu) \spseteq \pre{a}(\nubisbis),\quad
\suc{a}(\nu) \sbseteq \suc{a}(\nubisbis).
\end{equation}
As before, if $\vbis_{ba}=v_{ba}$ then $\vbisbis=v$ and the preceding relations hold as equalities.

Finally, (\ref{eq:mono1}--\ref{eq:mono2}) are obtained by combining (\ref{eq:lema-monotonia-intermitja-1}) and (\ref{eq:lema-monotonia-intermitja-2}):
\begin{gather*}
\iscbis_{ay} \ge \iscbisbis_{ay} = \isc_{ay},\\[2.5pt]
\iscbis_{xa} = \iscbisbis_{xa} \le \isc_{xa},\\[2.5pt]
\pre{a}(\nubis) \sbseteq \pre{a}(\nubisbis) \sbseteq \pre{a}(\nu),\\[2.5pt]
\suc{a}(\nubis) \spseteq \suc{a}(\nubisbis) \spseteq \suc{a}(\nu).
\end{gather*}

\vskip-3pt
\medskip\noindent
\textsl{General case.}\ensep
In the general situation~(\ref{eq:mona}), properties (\ref{eq:mono1}--\ref{eq:mono2})
are a direct consequence of the successive application of the special case~2 to every pair~$ay$.
\end{proof}

\smallskip
\begin{corollary}\hskip.5em
Under the hypothesis of Theorem~\ref{st:mono} one has also
\begin{equation}
\label{eq:monophi}
\flr_a > \flr_y \,\Rightarrow\, \flrbis_a \ge \flrbis_y.
\end{equation}
\end{corollary}
\begin{proof}\hskip.5em
It suffices to combine (\ref{eq:mono3}) with Theorem~\ref{st:phis}.
\end{proof}

\begin{comment}
The fraction-like rates are monotone functions of the projected scores \cite{cg}.
\end{comment}

\medskip
\begin{corollary}[\footnote{We thank Markus Schulze for pointing out this fact.}]\unskip
Under the hypothesis of Theorem~\ref{st:mono} one has also the property \textup{(\ref{eq:smona})}.
\end{corollary}
\begin{proof}\hskip.5em
According to Theorem~\ref{st:RvsNu}.(b), the left-hand side of (\ref{eq:smona}) implies the strict inequality $\isc_{ay} > \isc_{ya}$ for all $y\neq a$. Now, this inequality can be combined with (\ref{eq:mono1}) to derive that $\iscbis_{ay} > \iscbis_{ya}$ for all $y\neq a$. Finally, Theorem~\ref{st:RvsNu}.(c) with $\xst=\{a\}$ and $\yst=\ist\setminus\{a\}$ guarantees that the right-hand side of (\ref{eq:smona}) is satisfied.
\end{proof}


\paragraph{18.2}
The statements (\ref{eq:rmona}) and (\ref{eq:rrmona}) can fail even in the complete case. Next we give an example of it, with 5~options (it seems to be the minimum for the failure of (\ref{eq:rrmona})\,) and 10~voters. The only change from left to right is one inversion in one of the votes; more specifically, the eighth ballot changes from the order ¶d\better ¶b\better ¶c\better ¶a\better ¶e to the new one ¶b\better ¶d\better ¶c\better ¶a\better ¶e. In spite of this change, favourable to ¶b and disadvantageous to ¶d, the rank-like rate of ¶b is worsened  from $2.90$ to $3.00$, whereas that of ¶d is improved from $3.10$ to $3.00$. This contradicts (\ref{eq:rmona}) for~$a = ¶b$, as well as (\ref{eq:rrmona}.1) for~$a = ¶b$ and $y = ¶d\!,\!¶c$,
and also (\ref{eq:rrmona}.2) for~$a = ¶d$ and $y = ¶a\!,\!¶b$ (when one goes from right to left).
However, it complies with~(\ref{eq:wmona}).

\bgroup\small
\initsizesmall

\bigskip
\leftline{\hskip2\parindent
\taula{}%
\qnumx a/b/c/d/e◊%
\qmulticol(10)\hskip3.5pt Ranking votes:
1|2|2|2|2|3|3|4|4|5/
3|1|3|5|5|1|5|\textbf{2}|1|1/
5|3|4|3|3|2|2|3|3|2/
2|5|1|4|4|4|4|\textbf{1}|5|3/
4|4|5|1|1|5|1|5|2|4◊%
\fitaula
\hskip1.75em
\taula{}%
\qnumx a/b/c/d/e◊%
\qmulticol(10)\hskip3.5pt Ranking votes:
1|2|2|2|2|3|3|4|4|5/
3|1|3|5|5|1|5|\textbf{1}|1|1/
5|3|4|3|3|2|2|3|3|2/
2|5|1|4|4|4|4|\textbf{2}|5|3/
4|4|5|1|1|5|1|5|2|4◊%
\fitaula
}

\newcommand\septaules{\hskip2.85cm}
\def\xnumx#1◊{\bloc1\hfil\hfil\sf{\lower1pt\hbox{\Strut}}{$x$}#1◊} %
\def\xhmulticol(#1)#2:#3:#4◊{\bloc{#1}\hfil\hfil\rm{\lower1pt\hbox{\Strut}#2}{#3}{#4}◊}

\bigskip
\leftline{\hskip2\parindent
\taula{}%
\numx a/b/c/d/e◊%
\hmulticol(5)\hskip3.5pt$V_{xy}$:
\sf a|\sf b|\sf c|\sf d|\sf e:
\st|5|5|7|5/
5|\st|7|\textbf{4}|7/
5|3|\st|7|5/
3|\textbf{6}|3|\st|5/
5|3|5|5|\st◊%
\fitaula
\septaules
\taula{}%
\numx a/b/c/d/e◊%
\hmulticol(5)\hskip3.5pt\lower1pt\hbox{$\widetilde V_{xy}$}:
\sf a|\sf b|\sf c|\sf d|\sf e:
\st|5|5|7|5/
5|\st|7|\textbf{5}|7/
5|3|\st|7|5/
3|\textbf{5}|3|\st|5/
5|3|5|5|\st◊%
\fitaula
}

\newpage

\bigskip
\leftline{\hskip2\parindent
\taula{}%
\numx a/b/c/d/e◊%
\hmulticol(5)\hskip3.5pt$V^\ast_{xy}$:
\sf a|\sf b|\sf c|\sf d|\sf e:
\st|6|6|7|6/
5|\st|7|7|7/
5|6|\st|7|6/
5|6|6|\st|6/
5|5|5|5|\st◊%
\fitaula
\septaules
\taula{}%
\xnumx a/b/c/d/e◊%
\hmulticol(5)\hskip3.5pt\lower1pt\hbox{$\widetilde V^\ast_{xy}$}:
\sf a|\sf b|\sf c|\sf d|\sf e:
\st|5|5|7|5/
5|\st|7|7|7/
5|5|\st|7|5/
5|5|5|\st|5/
5|5|5|5|\st◊%
\fitaula
}

\bigskip
\leftline{\hskip2\parindent
\taula{}%
\numx a/b/c/d/e◊%
\hmulticol(5)\hskip3.5pt$M^\mast_{xy}$:
\sf a|\sf b|\sf c|\sf d|\sf e:
\st|1|1|2|1/
\st|\st|1|1|2/
\st|\st|\st|1|1/
\st|\st|\st|\st|1/
\st|\st|\st|\st|\st◊%
\fitaula
\septaules
\taula{}%
\xnumx a/b/c/d/e◊%
\hmulticol(5)\hskip3.5pt\lower1pt\hbox{$\widetilde M^\mast_{xy}$}:
\sf a|\sf b|\sf c|\sf d|\sf e:
\st|0|0|2|0/
\st|\st|2|2|2/
\st|\st|\st|2|0/
\st|\st|\st|\st|0/
\st|\st|\st|\st|\st◊%
\fitaula
}

\bigskip
\leftline{\hskip2\parindent
\taula{}%
\numx a/b/c/d/e◊%
\hmulticol(5)\hskip3.5pt$M^\pi_{xy}$:
\sf a|\sf b|\sf c|\sf d|\sf e:
\st|1|1|1|1/
\st|\st|1|1|1/
\st|\st|\st|1|1/
\st|\st|\st|\st|1/
\st|\st|\st|\st|\st◊%
\hmulticol(1):$r_x$:
2.80/
2.90/
3.00/
3.10/
3.20◊%
\fitaula
\hskip1.8cm
\taula{}%
\xnumx a/b/c/d/e◊%
\hmulticol(5)\hskip3.5pt\lower1pt\hbox{$\widetilde M^\pi_{xy}$}:
\sf a|\sf b|\sf c|\sf d|\sf e:
\st|0|0|0|0/
\st|\st|0|0|0/
\st|\st|\st|0|0/
\st|\st|\st|\st|0/
\st|\st|\st|\st|\st◊%
\xhmulticol(1):$\widetilde r_x$:
3.00/
3.00/
3.00/
3.00/
3.00◊%
\fitaula
}

\egroup

\bigskip
As one can see, the multiple zeroes present in $\widetilde M^\mast_{xy}$ force a complete tie of the rank-like rates $\widetilde r_x$ in spite of the fact that $\widetilde\nu$ is not empty. Not only the latter contains the pair \textsf{bd}, but in fact $\widetilde M^\mast_{\textsf{bd}} = 2 > 1 = M^\mast_{\textsf{bd}}$.


\end{document}